\documentclass[onecolumn,final,sort&compress]{elsart3p}

\usepackage{graphicx}
\usepackage{amsmath,amssymb}
\usepackage{graphics}
\usepackage{graphicx}
\usepackage[latin1]{inputenc}
\usepackage{epsfig}
\usepackage{comment}

\newcommand{\R}{\mathbb{R}}

\newcommand{\Xnewnm}{\tilde x_{nm} }
\newcommand{\Xoldnm}{x_{n}}
\newcommand{\Vnewnm}{\tilde v_{nm} }
\newcommand{\Voldnm}{v_{m}}

\newcommand{\Xvnewnm}{\tilde\xv_{nm} }

\newcommand{\Vvnewnm}{\tilde\vv_{nm} }

\newcommand{\incrXnm}{ v_{m} \, \Delta t}

\newcommand{\incrVnmb}{E_N(x_{n}) \, \Delta t}
\newcommand{\incrXnmVnmb}{ v_{m} \,E_N(x_{n})\, \Delta t^{2}}

\newcommand{\XnewnmIIord}{\hat x_{nm} }
\newcommand{\VnewnmIIord}{\hat v_{nm} }

\newcommand{\EOD}{\end{document}}

\newcommand{\species}{{\tt s}}
\newcommand{\fss}{\fs^{\species}}
\newcommand{\dx}{dx}
\newcommand{\dv}{dv}

\newcommand{\dvv}{d\vv}

\newcommand{\xvt}{\tilde\xv}
\newcommand{\vvt}{\tilde\vv}
\newcommand{\xst}{\tilde\xs}
\newcommand{\vst}{\tilde\vs}

\renewcommand{\cv}{\mathbf{c}}

\newcommand{\vv}{\mathbf{v}}
\newcommand{\wv}{\mathbf{w}}
\newcommand{\xv}{\mathbf{x}}

\newcommand{\Bv}{\mathbf{B}}

\newcommand{\Ev}{\mathbf{E}}

\newcommand{\Xv}{\mathbf{X}}
\newcommand{\Yv}{\mathbf{Y}}

\newcommand{\as}{a}
\newcommand{\bs}{b}
\newcommand{\cs}{c}
\newcommand{\ds}{d}

\newcommand{\fs}{f}
\newcommand{\gs}{g}

\newcommand{\ms}{m}

\newcommand{\qs}{q}

\newcommand{\vs}{v}
\newcommand{\ws}{w}
\newcommand{\xs}{x}

\newcommand{\Bs}{B}
\newcommand{\Cs}{C}

\newcommand{\Es}{E}

\newcommand{\Gs}{G}

\newcommand{\Ps}{P}
\newcommand{\Qs}{Q}

\newcommand{\Vs}{V}

\newcommand{\Xs}{X}

\newcommand{\abs}[1]{\vert#1\vert}

\usepackage{xcolor}

\begin{document}

\begin{frontmatter}
  
  \title{Arbitrary-Order Time-Accurate Semi-Lagrangian Spectral
    Approximations of the Vlasov-Poisson System}

  \author[Unicam] {L. Fatone}
  \author[Unimore,CNR]{, D. Funaro}
  \author[LANLT5] {, and G. Manzini}

  \address[Unicam]{
    Dipartimento di Matematica,  Universit\`a degli Studi di Camerino,
    Italy;
    \emph{e-mail: lorella.fatone@unicam.it}
  }
  \address[Unimore]{
    Dipartimento di Scienze Fisiche, Informatiche e Matematiche,
    Universit\`a degli Studi di Modena \\ e Reggio Emilia,
    Italy;
    \emph{e-mail: daniele.funaro@unimore.it}
  }
  \address[CNR]{
    Istituto di Matematica Applicata e Tecnologie Informatiche,
    Consiglio Nazionale \\delle Ricerche,
    via Ferrata 1,
    27100 Pavia,
    %\emph{e-mail: daniele.funaro@unimore.it}
  }
  \address[LANLT5]{
    Group T-5, 
    Applied Mathematics and Plasma Physics,
    Theoretical Division,\\
    Los Alamos National Laboratory,
    Los Alamos, NM,
    USA;
    \emph{e-mail: gmanzini@lanl.gov}
  }

  \maketitle
  \begin{abstract}
    The Vlasov-Poisson system, modeling the evolution of
    non-collisional plasmas in the electrostatic limit, is
    approximated by a Semi-Lagrangian technique.
    Spectral methods of periodic type are implemented through a
    collocation approach.
    Groups of particles are represented by the Fourier Lagrangian
    basis and evolve, for a single timestep, along an high-order
    accurate representation of the local characteristic lines.
    The time-advancing technique is based on Taylor developments that
    can be, in principle, of any order of accuracy, or by coupling the
    phase space discretization with high-order accurate Backward
    Differentiation Formulas (BDF) as in the method-of-lines
    framework.
    At each timestep, particle displacements are reinterpolated and
    expressed in the original basis to guarantee the order of accuracy
    in all the variables at relatively low costs.
    Thus, these techniques combine excellent features of spectral
    approximations with high-order time integration.
    Series of numerical experiments are performed in order to assess
    the real performance.
    In particular, comparisons with standard benchmarks are examined.
  \end{abstract}
  
\end{frontmatter}

\raggedbottom

%%% PAPER
%\input introduction.tex

%%%%%%%%%%%%%%%%%%%%%%%%%%%%%%%%%%%%%%%%%%%%%%%%%%%%%%%%%
\section{Introduction}\label{sec1}
%%%%%%%%%%%%%%%%%%%%%%%%%%%%%%%%%%%%%%%%%%%%%%%%%%%%%%%%%%
%%
%% \noindent 
%%
The Vlasov-Poisson system of equations describes the dynamics of a
collisionless plasma of charged particles (electrons and ions), where
the only relevant interaction is driven by a self-consistent
electrostatic field~\cite{Boyd-Sandserson:2003}.
Although the Vlasov-Poisson system is one of the simplest models that
can be considered in plasma physics, its numerical treatment is quite
challenging to the numerical modelers.
In fact, each plasma species is described by a distribution function
that is defined on a high-dimensional phase space.
Since the beginning of numerical plasma simulations in the '60s, a
number of methods have been proposed to the scientific
community and thoroughly investigated.
%%
%% 4 families of methods
We can roughly regroup them in a few big families: Particle-in-Cell
(PIC) methods, Transform methods, Eulerian and Semi-Lagrangian
methods.

The PIC method is very popular in the plasma physics community, as it
is the most widely used method because of its robustness and relative
simplicity~\cite{Birdsall-Langdon:2005}.
There, the evolution of a plasma is described by the motion of a
finite number of macro-particles in the physical space.
These macro-particles are tracked along the characteristics of the
Vlasov equation and their mutual interaction is driven by a
nonlinearly coupled electric field, which solves the Poisson equation.
The right-hand side of the Poisson equation depends on the charges
carried by the macro-particles.
The convergence of the PIC method for the Vlasov-Poisson system was
proved
in~\cite{Cottet-Raviart:1984,Wollman:2000,Wollman-Ozizmir:1996}.
%%
% PIC does not suffer of the curse of dimensionality
The PIC method has been successfully used to simulate the behavior of
collisionless laboratory and space plasmas and provides excellent
results for the modeling of large scale phenomena in one, two or three
space dimensions~\cite{Birdsall-Langdon:2005}.
Also, implicit and energy preserving PIC formulations that are
suitable to long time integration problems are available from the most
recent
literature~\cite{Brackbill:2016,Chen-Chacon:2015,Chen-Chacon-Barnes:2011,Lapenta:2017,Lapenta-2011,Markidis-Lapenta:2011,Taitano-Knoll-Chacon-Chen:2013}.
%%
%% but there are other drawbacks
Nonetheless, PIC codes suffer from intrinsic drawbacks.
As proved in~\cite{Cottet-Raviart:1984}, achieving 
%%``reasonable'' approximation standards requires
high numerical resolution in multidimensional plasma physics
simulations may require a huge number of particles, thus making
such simulations infeasible even with the most powerful
supercomputers currently available.
Since only a relatively limited number of particles can be considered
in practical calculations, the method is used in a suboptimal way and
tends to be intrinsically noisy.
Although research has been carried out to reduce PIC
noise~\cite{Myers-Colella-VanStraalen:2017}, the method remains effective
mainly for problems with a low noise-to-signal ratio, and where the
physics is not driven by fine phase space structures.

\medskip
%% Transform methods
Based on the seminal paper~\cite{Grad:1949}, an alternative approach,
called the Transform method, was developed at the end of the '60s, which uses
a spectral decomposition of the distribution function and leads to a
truncated set of moment equations for the expansion
coefficients~\cite{Armstrong:1970}.
To this end, Hermite basis functions are used for unbounded domains,
Legendre basis functions for bounded domains, and Fourier basis
functions for periodic domains, see,
e.g.,~\cite{Manzini:2016:JCP:journal,Parker-Dellar:2015,Klimas:1983,Vencels-Delzanno-BoPeng-Laure-Markidis:2015,Vencels:2016:JPCS:journal}.
These techniques can outperform
PIC~\cite{Camporeale-Delzanno-Vergen-Moulton:2015,Camporeale-Delzanno-Lapenta-Daughton:2006}
in Vlasov-Poisson simulations.
Moreover, they can be extended in an almost straightforward way to
multidimensional simulations of more complex models, like
Vlasov-Maxwell~\cite{Delzanno:2015}.
Convergence of various formulations of these methods was shown
in~\cite{Gajewski-Zacharias:1977,Manzini-Funaro-Delzanno:2017:journal}.
Transform methods offer a few indisputable advantages.
First of all, they may be extremely accurate since they
are based on a spectral approximations of the differential operators.
Furthermore, physically meaningful discrete invariants (such as total
number of particles, momentum and total energy) can be built directly
from the expansion
coefficients~\cite{Schumer-Holloway:1998,Holloway:1996}.
The existence of such discrete invariants implies better stability
properties in long-time integration problems.
However, despite their good properties their implementation may be computational
demanding.
As a matter of fact, they suffer of the ``curse of dimensionality'' (i.e., a bad
scaling of the computational complexity with the number of
dimensions), when multidimensional basis functions are built by tensor
product of one-dimensional ones.

%% Eulerian-Vlasov and semi-Lagrangian method
%%
\medskip
An alternative to PIC and Transform methods is offered by the class of
Eulerian and Semi-Lagrangian methods, which discretize the Vlasov
equation on a grid of the phase space.
Common approaches for the implementation are: Finite Volume
Methods~\cite{Filbet:2001,Banks-Hittinger:2010}, Discontinuous
Galerkin~\cite{Ayuso-Carrillo-Shu:2011,Ayuso-Carrillo-Shu:2012,Heath-Gamba-Morrison-Michler:2012},
finite difference methods based on ENO and WENO
polynomial reconstructions~\cite{Christlieb-Guo-Morton-Qiu:2014}, or
propagation of the solution along the characteristics in an operator
splitting framework~\cite{Arber-Vann:2002,%
  Carrillo-Vecil:2007,%
  Cheng-Knorr:1976,%
  Filbet-Sonnendrucker-Bertrand:2001,%
  Filbet-Sonnendrucker:2003,%
  Sonnendrucker-Roche-Bertrand-Ghizzo:1999,%
  Crouseilles-Respaud-Sonnendrucker:2009}.
Semi-Lagrangian methods were first developed for meteorological
applications in the early '90s
\cite{Bermejo:1991,Bermejo:1995,Staniforth-Cote:1991}.
The aim was to take advantage of both Lagrangian and Eulerian
approaches.
Indeed, these methods allow for a relatively accurate description of the phase
space using a fixed mesh and propagating the values of the
distribution function along the characteristics curves forward or
backward in time.
High-dimensionality is typically addressed by a splitting operator
strategy in order to advance the solution in time.
Such a splitting makes it possible to approximate a multi-dimensional
time-dependent problem by a sequence of one-dimensional problems.
For the one-dimensional Vlasov-Poisson system, the splitting
reformulates the Vlasov equation in two advection subproblems that
advance the distribution function in space and velocity independently.
High-order approximations are described
in~\cite{Qiu-Christlieb:2010}.

\medskip
%%
%% WHAT WE ARE DOING
%%
In this paper, we propose  Semi-Lagrangian methods that provide the
spectral accuracy of Transform methods.
In particular, space and velocity representations are discretized
using a spectral collocation approach and the approximation of
the distribution function is advanced in time by following backward
the characteristic curves.
Furthermore, we do not resort to any time splitting of the Vlasov
equation and the desired order of accuracy in time, e.g.,
$\mathcal{O}(\Delta t^2)$ or even higher, is attained by using 
well calibrated representations of the characteristic
curves. The major advantage of our approach is to combine, in a simple
and natural way, spectral accuracy with on purpose time discretization techniques,
in principle of any order of convergence.
The formulation of the method is the same for any space and velocity
dimension, provided we adopt a multi-index notation.
Finally, although we do not address these topics in the present paper,
we note that an efficient implementation is possible by resorting to
standard libraries such as the Discrete Fast Fourier Transform
(DFT)]~\cite{Brigham:1988}.
Moreover, we remark that unsplit algorithms, like the ones that we
propose in this work, are more suited to task parallelization on
multicore processors in comparison to split algorithms in the standard
Semi-Lagrangian approaches.

\medskip
The paper is organized as follows.
In Section~\ref{sec:continuous:Model}, we present the continuous model.
In Section~\ref{sec:phase-space discretization}, we introduce the
spectral approximation in the phase space.
In Section~\ref{sec:time:discretization}, we present a
Semi-Lagrangian scheme based on a first-order accurate
approximation of the characteristic curves, making use of a
suitable Taylor expansion.
In Section~\ref{sec:more:advanced:time:discretizations}, we derive
more refined time discretization schemes, built in the framework of
the method-of-lines, applying second-order and third-order multi-step
Backward Differentiation Formula (BDF).
%%$
To show the flexibility of our approach, we also present a single-step
second-order approximation in time.
In Section~\ref{sec:conservation:properties}, we investigate the
conservation properties of the method and we show that the number of
particles is always an exact invariant of the method, regardless of
the order of the time discretization. Within a spectral accurate error,
this is also true for momenta.
Concerning the total energy, this is conserved up to an approximation error
that depends on the accuracy of the time discretization.
In Section~\ref{sec:numerical:experiments}, we show the predicted
convergence rate in time by using a manufactured solution.
Furthermore, we assess the performance of the method on standard
benchmark problems as the two stream instability and the Landau
damping.
In Section~\ref{sec:conclusions}, we present our final remarks and
conclusions.

%\input method.tex

%%%%%%%%%%%%%%%%%%%%%%%%%%%%%%%%%%%%%%%%%%%%%%%%%%%%%%%%%
\section{The continuous model}
\label{sec:continuous:Model}
%%%%%%%%%%%%%%%%%%%%%%%%%%%%%%%%%%%%%%%%%%%%%%%%%%%%%%%%%%

%% multidimensional formulation 
%%
\subsection{Multidimensional multispecies formulation}
The distribution functions $\fss(t,\xv,\vv)$,
$\species=1,2,\ldots,n^{\species}$, solving the Vlasov-Poisson system
describe the statistical evolution of a collection of collisionless
charged particles of $n^{\species}$ distinct species, subject to
mutual electrostatic interactions~\cite{Boyd-Sandserson:2003}.
From a physical viewpoint, each $\fss(t,\xv,\vv)\dx\dv$ represents
\emph{the probability of finding particles of species $\species$ in an
  element of volume $\dx\dv$, at time $t$ and point $(\xv,\vv)$ in the
  phase space
  $\Omega=\Omega_x\times\Omega_v$}, where
$\Omega_{x}\subseteq\R^{3}$, $\Omega_{v}\subseteq\R^{3}$.
The 3D-3V Vlasov equation for the $\species$-th species with mass
$\ms^{\species}$ and electric charge $\qs^{\species}$ reads as:
\begin{equation}
  \frac{\partial\fss}{\partial t} + \vv\cdot\nabla_{\xv}\fss
  +\frac{\qs^{\species}}{\ms^{\species}}\Ev\cdot\nabla_{\vv}\fs = 0, 
  \quad 
  t\in (0,T],\,\,
  \xv\in\Omega_{x},\,\,
  \vv\in\Omega_{v},
  \label{eq:3D:V}
\end{equation}
where $\Ev(t,\xv)$ represents the electric field.
The initial condition for $\fss$ is given by a function
$\bar{f}^{\species}$, so that
\begin{equation}
  \fss(0,\xv,\vv)=\bar{\fs}^{\species}(\xv,\vv),
  \quad\species=1,\ldots,n^{\species},\, \, \xv\in\Omega_{x},\,\,\vv\in\Omega_{v}.
  \label{eq:3D:Vci}
\end{equation}
The coupling with the self-consistent electric field $\Ev(t,\xv)$ is
taken into account through the divergence equation:
\begin{equation}
  \epsilon_{0}\big(\nabla\cdot\Ev\big)(t,\xv) 
  = \qs^{\species}\sum_{\species=1}^{n^{\species}}\rho^{\species}(t,\xv)
  = \qs^{\species}\sum_{\species=1}^{n^{\species}}\int_{\Omega_{\vv}}\fss(t,\xv,\vv)\dvv , 
  \qquad t\in[0,T],\,\,\xv\in\Omega_{x},
  \label{eq:3D:EMaxwell}
\end{equation}
where $\rho^{\species}(t,\xv)$ is the charge density of species $\species$.
In~\eqref{eq:3D:EMaxwell} $\epsilon_{0}$ is the dielectric vacuum
permittivity and $\rho(t,\xv)$, is the total charge density.
We refer the reader interested in the theoretical analysis of the
Vlasov-Poisson model and the properties of its solutions to
\cite{Bouchut-Golse-Pulvirenti:2000,Glassey:1996,Dolbeault:2002}.

\subsection{1D-1V formulation of the Vlasov-Poisson system}
To ease the presentation of the numerical scheme, we consider the
1D-1V Vlasov-Poisson formulation for the electron-ion coupled system.
Consistently, we restrict the domain to $\Omega_x\subseteq\R$ and
$\Omega_v\subseteq\R$.
Since positive ions (protons) are much heavier than electrons, we may
assume that they do not move, so that their density distribution
function is constant over $\Omega_{x}$.
Without altering the generality of the exposition, we can set
$\qs=-1$, $\ms=1$, $\epsilon_{0}=1$.
By dropping out the label $\species$, we only have one distribution function
$\fs$ for the electron species, so that the corresponding Vlasov
equation and initial condition read as:
\begin{align}
  &
  \frac{\partial\fs}{\partial t} + \vs\frac{\partial\fs}{\partial\xs} 
  -\Es(t,\xs)\ \frac{\partial\fs}{\partial\vs}=0, \qquad t\in (0,T],\,\,\xs\in\Omega_{x},\,\,
  \vs\in\Omega_{v},
  \label{eq:1D1V:V}\\[3mm]
  & \fs(0,\xs,\vs)=\bar{\fs}(\xs,\vs),\quad\xs\in\Omega_{x},\,\,\vs\in\Omega_{v},
  \label{eq:1D1V:Vci}
\end{align}
where the coupled electric field $\Es$ verifies the equation:
\begin{equation}
  \label{eq:1D1V:EMaxwell}
  \frac{\partial\Es}{\partial\xs}(t,\xs) = 1-\rho(t,\xs),
  \quad t\in [0,T],\,\,\xs\in\Omega_{x}.
\end{equation}
We recall that $\rho(t,\xs)$ is the electron charge density defined
by:
\begin{equation}\label{eq:1D1V:charge-density}
  \rho(t,\xs)= \int_{\Omega_{v}}\fs(t,\xs,\vs)\dv.
  %% \quad t\in [0,T],\,\,\xs\in\Omega_{x}.
\end{equation}
%%
% Equations~\eqref{eq:1D1V:Vci} prescribe the initial conditions for the
% particle distribution functions $\fs$, while the electric field $\Es$
% is now given at every time instant by
% solving~\eqref{eq:1D1V:EMaxwell}.
%% 
We assume the constraints (charge conservation):
% \footnote{\textbf{Marco: Deve essere
%     $\int_{\Omega_x}(1-\rho(t,\xs))\dx=0$ perch\'e il plasma \`e
%     globalmente neutro, quindi l'integrale in spazio della densit\`a
%     di carica, che \`e il rhs della Eq.~\eqref{eq:1D1V:EMaxwell} deve
%     fare zero, da cui segue~\eqref{eq:1D1V:intrhoconst}.} }
%% 
\begin{equation} 
  \label{eq:1D1V:intrhoconst}
  \int_{\Omega_{x}}\Es(t,\xs)dx = 0,
  \qquad \textrm{which~implies~that}\qquad 
  \int_{\Omega_{x}}\rho(t,\xs)dx = \abs{\Omega_x}, %\quad t\in [0,T],
\end{equation}
where $\abs{\Omega_x}$ measures the extension of $\Omega_x$.
%%
%% By differentiation, equation~\eqref{eq:1D1V:EMaxwell} is usually
By taking 
\begin{equation} 
  \Es(t,\xs)=- \frac{\partial \Phi}{\partial\xs}(t,\xs), 
\end{equation}
equation~\eqref{eq:1D1V:EMaxwell} can be transformed into the Poisson
equation for the potential field $\Phi(t,\xs)$:
\begin{equation} 
  \label{eq:1D1V:EPoisson}
  -\frac{\partial^{2}\Phi}{\partial\xs^{2}}(t,\xs) = 1-\rho(t,\xs).
  %% \qquad\textrm{with}~
  %% \Es(t,\xs)=- \frac{\partial \Phi}{\partial\xs}(t,\xs), 
  %% \quad t\in [0,T],\,\,\xs\in\Omega_{x}
\end{equation}
% %%
% Let $(\fs,\Es)$ be the solution of the Vlasov-Poisson problem defined
% by equations~\eqref{eq:1D1V:V}, \eqref{eq:1D1V:Vci},
% \eqref{eq:1D1V:EMaxwell}, \eqref{eq:1D1V:charge-density}.
%% 
As far as boundary constraints in $x$ and $v$ are concerned, we will
assume a periodic boundary condition for the Poisson equation and
either periodic or homogeneous Dirichlet boundary conditions for the
Vlasov equation.

%\subsection{Properties}
%%
%It is well-known (see.....) that in 

%%
In the continuum setting, the total number of plasma particles is
preserved.
Hence, from a straightforward calculation and
using~\eqref{eq:1D1V:intrhoconst} it follows that:
\begin{equation} 
  \label{eq:1D1V:massConserv}
  \frac{d}{dt}\int_{\Omega}\fs(t,\xs,\vs)\,dx\,dv = 0. %,\qquad t\in [0,T].
\end{equation}
Moreover, the distribution function $\fs$ solving the Vlasov-Poisson
system satisfies the so-called $L^{p}$-stability property for $p\ge
1$:
\begin{equation}
  \label{eq:1D1V:L2Conserv}
  \frac{d}{dt}\,\|f(t,\cdot,\cdot)\|^{p}_{L^{p}} =  
  \frac{d}{dt}\,\int_{\Omega} |\fs(t,\xs,\vs)|^{p}\,dx\,dv=0, %%,\qquad t\in [0,T].
\end{equation}
which holds for any $t\in[0,T].$
In particular, we will be concerned with $p=2$.
In this case, \eqref{eq:1D1V:L2Conserv} implies the $L^{2}$-stability
of the method~\cite{Glassey:1996} (sometimes called also \emph{energy
  stability} in the literature).

Finally, we consider the total energy of the system defined by:
\begin{equation}
  \label{eq:1D1V:totEnergy}
  \mathcal{E}(t) 
  = \frac12 \int_{\Omega}\fs(t,\xs,\vs)\,{|\vs|^{2}}\,dx\,dv\, 
  + \,\frac12\int_{\Omega_{x}}|\Es(t,\xs)|^{2 }\,dx,
  % ,\qquad t\in [0,T],
\end{equation}
where the first term represents the kinetic energy 
% at time $t\in [0,T]$ 
and the second one the potential energy.
% at time $t\in [0,T]$.
%%
The Vlasov-Poisson model is characterized by the exact conservation of
the energy, i.e.:
\begin{equation}
  \label{eq:1D1V:energyConserv}
  \frac{d}{dt}\,\mathcal{E} (t)=0.  %,\qquad t\in [0,T].
\end{equation}
If the electric field is smooth enough, for a ``sufficiently small''
$\delta>0$, the local system of characteristics associated
with~\eqref{eq:1D1V:V} is given by the phase space curves
$(\Xs(\tau),\Vs(\tau))$ solving 
\begin{equation}
  \label{eq:1D1V:char}
  \frac{d\Xs}{d \tau} = -\Vs(\tau),\qquad
  \frac{d\Vs}{d \tau} =  \Es(\tau,\Xs(\tau)), \qquad \tau \in ]t-\delta ,t+\delta [,
\end{equation}
with the condition that $(\Xs(t),\Vs(t))=(\xs,\vs)$ when $\tau=t$.
Under suitable regularity assumptions, there exists a unique solution
of the Vlasov-Poisson problem \eqref{eq:1D1V:V}, \eqref{eq:1D1V:Vci},
\eqref{eq:1D1V:EMaxwell} and \eqref{eq:1D1V:charge-density},
see~\cite{Glassey:1996}, which can formally be expressed by
propagating the initial condition~\eqref{eq:1D1V:Vci} along the
characteristic curves that solve~\eqref{eq:1D1V:char}.
Therefore, for every $t\in(0,T]$ we have that
\begin{equation}
  \label{eq:1D1V:solChar}
  \fs(t,\xs,\vs) = \bar{\fs}(\Xs(t),\Vs(t)).
  %% \qquad t\in[0,T],\,\,\xs\in\Omega_{x},\,\,\vs\in\Omega_{v}.
\end{equation}
By using a first-order approximation of the characteristic curves given
by:
\begin{equation}
  \label{eq:1D1V:char1}
  \Xs(\tau) = \xs - \vs (\tau-t),\qquad
  \Vs(\tau) = \vs + \Es(t,\xs)(\tau-t) ,
\end{equation}
the Vlasov equation is satisfied up to an error that decays as
$(\tau-t)$, for $\tau$ tending to $t$.
To achieve a higher order of convergence, we need a more accurate
approximation of the characteristic curves, such as, for example, the
one given by setting:
\begin{align}
  \Xs(\tau) &= \xs - \vs (\tau-t) - \frac12\Es(t,x) (\tau-t)^2, \nonumber\\[3mm]%\quad s\in[0,T],\nonumber\\[3mm]
  \Vs(\tau) &= \vs + \Es(t,\xs)(\tau-t) 
  -\frac12 \left( 
    \frac{\partial\Es}{\partial t}(t,\xs) 
    + \vs\frac{\partial\Es}{\partial\xs}(t,\xs)
  \right) (\tau-t)^2. %,\qquad s\in[0,T].
  \label{eq:1D1V:char2}
\end{align}
By direct substitution in~\eqref{eq:1D1V:V}, the Vlasov equation is
satisfied at every point $(t,\xs,\vs)$ up to the quadratic remainder
$(\tau-t)^2$ for $\tau$ tending to $t$.
Of course, \eqref{eq:1D1V:char2} can be replaced by other more
accurate expansions leading to a high-order remainder term
proportional to $(\tau-t)^{S}$ for some integer $S>2$.
Without exhibiting the explicit formulas, which look rather involved,
we point out this property as a possible extension for further
generalizations.

In view of the expression above, it is also convenient to write the
time derivative of the electric field $\Es$ by arguing as follows. 
We evaluate the time derivative of $\rho$ in
\eqref{eq:1D1V:charge-density} and use the Vlasov-Poisson equation:
\begin{align}
  \frac{\partial\rho}{\partial t}(t,\xs)
  = -\int_{\Omega_{v}}\vs\frac{\partial\fs}{\partial x}(t,\xs,\vs)\,\dv 
  +   \Es(t,\xs)\int_{\Omega_{v}}\frac{\partial\fs}{\partial\vs} (t,\xs,\vs)\dv
  = -\int_{\Omega_{v}}\vs\frac{\partial\fs}{\partial x}(t,\xs,\vs)\,\dv,
  \label{eq:1D1V:chargeder}
\end{align}
%
% \begin{align}
%   \frac{\partial\rho}{\partial t}(t,\xs)
%   &= -\int_{\Omega_{v}}\vs\frac{\partial\fs}{\partial x}(t,\xs,\vs)\,\dv 
%   +   \Es(t,\xs)\int_{\Omega_{v}}\frac{\partial\fs}{\partial\vs} (t,\xs,\vs)\dv\nonumber\\[0.5em]
%   &= -\int_{\Omega_{v}}\vs\frac{\partial\fs}{\partial x}(t,\xs,\vs)\,\dv,\qquad t\in [0,T],\,\,\xs\in\Omega_{x},
%   \label{eq:1D1V:chargeder}
% \end{align}
%% 
where we observe that the integral of $\partial\fs/\partial\vs$ is
zero for a periodic function or in presence of homogeneous Dirichlet
conditions.
Translated in terms of $\Es$, the above equation implies the Amp\`ere
equation, which reads as:
\begin{equation}
  \frac{\partial\Es}{\partial t}(t,\xs) +
  \int_{\Omega_{v}}\vs\fs(t,\xs,\vs)\dv = C_A,
  \label{eq:1D1V:eder}
\end{equation}
(after an integration with respect to $\xs$).
Finally, in order to preserve the conditions in
\eqref{eq:1D1V:intrhoconst}, we must set $C_A=0$ in
\eqref{eq:1D1V:eder}.

%%%%%%%%%%%%%%%%%%%%%%%%%%%%%%%%%%%%%%%%%%%%%%%%%%%%%%%%%%
\section{Phase-space discretization}
\label{sec:phase-space discretization}
%%%%%%%%%%%%%%%%%%%%%%%%%%%%%%%%%%%%%%%%%%%%%%%%%%%%%%%%%%
%
We propose a Semi-Lagrangian method to find numerical approximations
to the self-consistent solutions of the 1D-1V Vlasov-Poisson problem
defined by equations \eqref{eq:1D1V:V}, \eqref{eq:1D1V:Vci},
\eqref{eq:1D1V:EMaxwell}  and \eqref{eq:1D1V:charge-density}.
%%
%% In this section we introduce the discretization in the one-dimensional
%% variables $\xs$ and $\vs$.
%%
The extension to higher-dimensional problems, e.g., the 3D-3V case, is
straightforward and is discussed at the end of this section.
Instead, in the subsequent sections, we will analyze suitable
time discretization techniques.
In view of imposing periodic boundary conditions, we start by
considering the domain:
\begin{equation}
  \Omega = \Omega_{x}\times\Omega_{v} = [0,2\pi[\times[0,2\pi[.
  \label{eq:1D1V:phaseSpaceReale}
\end{equation}
A function $\fs$ defined in $\Omega$ is requested to be periodic in
both $\xs$ and $\vs$.
This means that for any integer $s\ge 0$ we must have:
\begin{align}\label{eq:1D1V:PeriodicBC}
  \frac{\partial^{s}\fs}{\partial\xs^{s}}(0,\vs) = \frac{\partial^{s}\fs}{\partial\xs^{s}} (2\pi,\vs),\qquad\textrm{for~every~}\vs\in\Omega_{v},
\end{align}
and
\begin{align}
  \frac{\partial^{s}\fs}{\partial\vs^{s}}(\xs,0) = \frac{\partial^{s}\fs}{\partial\vs^{s}} (\xs,2\pi),\qquad\textrm{for~every~}\xs\in\Omega_{x},
\end{align}
where, as usual, the zero-th order derivative of the function (i.e.,
when $s=0$) is the given function itself.

\noindent
Given two positive integers $N$ and $M$, we consider the
equispaced points in $[0,2\pi[$:
\begin{equation}
  \xs_{i} = \frac{2\pi}{N}\,i ,\,\,\,i=0,1,\ldots, N-1, 
  \qquad \vs_{j} =\frac{2\pi }{M}\,j,\,\,\,j=0,1,\ldots, M-1.
\label{eq:1D1V:collpoints}
\end{equation}
Hereafter, if not otherwise indicated, we will always use the indices
$i$ and $n$ running from $0$ to $N-1$ to label the grid points along
the $x$-direction, and $j$ and $m$ running from $0$ to $M-1$ to label
the grid points along the $v$-direction.

Then, we introduce the Fourier Lagrangian
% ~\footnote{\textbf{Marco: si
%     dice Lagrangian oppure Lagrange? A me suona pi\'u corretto
%     Lagrange, per\`o non so quale dei due sia quello giusto (o se sono
%     giusti entrambi). Daniele: io userei Lagrangian}} 
basis functions for the $\xs$ and $\vs$ variables with respect to
the nodes \eqref{eq:1D1V:collpoints}, that is:
\begin{align}
  \Bs_{i}^{(N)}(\xs) &= \frac{1}{N} \sin\left(\frac{N(\xs-\xs_{i})}{2}\right)\,\cot\left(\frac{\xs-\xs_{i}}{2}\right),\label{eq:1D1V:Basis1}\\[3mm]
  \Bs_{j}^{(M)}(\vs) &= \frac{1}{M} \sin\left(\frac{M(\vs-\vs_{j})}{2}\right)\,\cot\left(\frac{\vs-\vs_{j}}{2}\right).\label{eq:1D1V:Basis2}
\end{align}
It is known that
\begin{equation}
  \Bs_{i}^{(N)} (\xs_{n}) = \delta_{in} \qquad\textrm{and}\qquad\Bs_{j}^{(M)}(\vs_{m}) = \delta_{jm},
  %%\Bs_{i}^{(N)} (\xs_{n}) = \delta_{in},\quad i,n=0,1,\ldots,N-1,\qquad\Bs_{j}^{(M)}(\vs_{m}) = \delta_{jm},\quad j,m=0,1,\ldots, M-1,
  \label{eq:1D1V:Kron}
\end{equation}
where $\delta_{ij}$ is the usual Kronecker symbol.

Furthermore, we define the discrete spaces:
\begin{align}
  \Xv_N    &= \textrm{span}\Big\{ \Bs_{i}^{(N)} \Big\}_{i=0,1,\ldots, N-1}
  \qquad\textrm{and}\qquad
  %%\nonumber\\[3mm]
  %%\qquad\xs\in\Omega_{x},\nonumber\\[3mm]
  \Yv_{N,M} = \textrm{span}\Big\{ \Bs_{i}^{(N)}\Bs_{j}^{(M)} \Big\}_{i=0,1,\ldots, N-1\atop j=0,1,\ldots, M-1}.
  %%\qquad\xs\in\Omega_{x},\,\vs\in\Omega_{v}.
  \label{eq:1D1V:spaces}
\end{align}
In this way, any function $\fs_{N,M}$ that belongs to $\Yv_{N,M}$ can
be decomposed as:
\begin{equation}
  \fs_{N,M}(\xs,\vs) = \sum_{i=0}^{N-1}\,\sum_{j=0}^{M-1} \cs_{ij}\,\Bs_{i}^{(N)}(\xs)\,\Bs_{j}^{(M)}(\vs),
  %%\qquad\xs\in\Omega_{x},\,\,\vs\in\Omega_{v},
  \label{eq:1D1V:fbasis}
\end{equation}
where the coefficients of the decomposition are given by:
\begin{equation}
  \cs_{ij} = \fs_{N,M} (\xs_{i},\vs_{j}).
  %%\quad i=0,1,\ldots,N-1,\,\,j=0,1,\ldots, M-1.
  \label{eq:1D1V:coeff}
\end{equation}

For what follows, it will be useful to have the expression of the
derivatives of the basis functions.
For instance, one has:
\begin{equation}
  \frac{\partial\Bs_{i}^{(N)}}{\partial\xs}(\xs_{n}) 
  = \ds_{ni}^{(N,1)} 
  =
  \begin{cases}
    0                                                            & \mbox{if~$i=n$},\\[2mm]
    \frac12 (-1)^{i+n}\,\cot\left(\frac{\xs_{n}-\xs_{i}}{2}\right) & \mbox{if~$i\neq n$},
  \end{cases}
  \label{eq:1D1V:Der1x}
\end{equation}
and
\begin{equation}
  \frac{\partial^{2}\Bs_{i}^{(N)}}{\partial\xs^{2}}(\xs_{n}) 
  = \ds_{ni}^{(N,2)}
  =
  \begin{cases}
    -\frac{N^{2}}{12} - \frac16                                            & \mbox{if~$i=n$},\\[4mm]
    -\frac12\,\frac{(-1)^{i+n}}{\sin^{2}\left(\frac{x_{n}-x_{i}}{2}\right)} & \mbox{if~$i\neq n$}.
  \end{cases}
  \label{eq:1D1V:Der2x}
\end{equation}
More generally, $\ds_{ni}^{(N,s)}$ will denote the $s$-th derivative
of $\Bs_{i}^{(N)}$ evaluated at point $\xs_n$, which is given by:
\begin{equation}
  \frac{\partial^{s}\Bs_{i}^{(N)}}{\partial\xs^{s}}(\xs_{n}) = \ds_{ni}^{(N,s)}.
  \label{eq:1D1V:Der1v}
\end{equation}
Analogously we can define:
\begin{equation}
  \frac{\partial^{s}\Bs_{j}^{(M)}}{\partial\vs^{s}}(\vs_{m}) = \ds_{mj}^{(M,s)},
  \label{eq:1D1V:Der2v}
\end{equation}
where $\ds_{mj}^{(M,1)}$, $\ds_{mj}^{(M,2)}$ in \eqref{eq:1D1V:Der2v}
are obtained by replacing the nodes $\xs_{i}$ with the nodes $\vs_{j}$
in \eqref{eq:1D1V:Der1x} and \eqref{eq:1D1V:Der2x} and setting up the
indices accordingly.
As a special case we set: $\ds_{ni}^{(N,0)}=\delta_{ni}$,
$\ds_{mj}^{(M,0)}=\delta_{mj}$.
Moreover, it is easy to prove that there exists a constant $C$,
independent of $N$, such that:
\begin{equation}
  \abs{ \ds_{ni}^{(N,1)} } \leq CN.
  \label{eq:1D1V:stimad}
\end{equation}
This estimate will be useful in the next section for studying the
stability conditions in the time-marching schemes.

Furthermore, we remind that the following Gaussian quadrature formula:
\begin{equation}\label{eq:1D1V:fdq}
  \displaystyle
  \frac{1}{2\pi}\,\int_{0}^{2\pi}\phi(\xs)\,d\xs\simeq\frac1N\sum_{i=0}^{N-1}\phi(\xs_{i}),
\end{equation}
which can be applied to any $\phi\in\Cs[0,2\pi)$, is exact for every
$\phi\in\textrm{span}\Big\{1,\big\{\sin nx,\cos nx\big\}_{n=1,\ldots,N-1},\sin Nx\Big\}$.
For more details
see~\cite[Section~2.1.2]{Canuto-Hussaini-Quarteroni-Zhang:1988} and
\cite[Section~2.1.2]{Sheng-Tao-Wang:2011}.
%% {\textcolor{red}{anche Shen, Tao, Wang, sempre section 2.1.2}.

In truth, given an integer $s \ge 0$, the derivative of order $s+1$ is
trivially obtained by applying the first derivative matrix to the
point-values of the $s$-th derivative of a trigonometric polynomial.
Such an operation can be performed by the \emph{fast Fourier
  transform} algorithm, with an excellent cost reduction when the
degree is relatively high and a power of 2, and very efficient
implementations exist in freely available and commercial software
libraries.

It is clear that, with little modifications, we can handle Lagrangian
basis of nonperiodic type. 
Among these, the most representative ones are constructed on Legendre
or Chebyshev algebraic polynomials, or Hermite functions (i.e.,
Hermite polynomials multiplied by a Gaussian function).
In some preliminary tests, we observed that each one of these cases
presents peculiar behavior in applications.
A comparison between the different approaches would be too lengthy for
the aims of the present paper. 
Therefore, we prefer to examine more deeply these extensions in a
future analysis.

Now, consider the one-dimensional function $\Es_N\in\Xv_N$. 
% is given.
%%
Given $\Delta t>0$, by taking $\tau=t -\Delta t$ in formula \eqref{eq:1D1V:char1}, we define the new
set of points $\{(\Xnewnm,\Vnewnm)\}_{n,m}$ where
\begin{align}
  \Xnewnm &= \Xoldnm-\incrXnm      ,    \label{eq:1D1V:puntimossi}\\[2mm]
  \Vnewnm &= \Voldnm+E_N(x_n)\Delta t,  \label{eq:1D1V:1charnm}
\end{align}
where we recall that index $n$ is running through the range $[0,N-1]$
and index $m$ through the range $[0,M-1]$.
To evaluate a function $\fs_{N,M}\in\Yv_{N,M}$ at the new points
$(\Xnewnm, \Vnewnm)$ through the coefficients in
\eqref{eq:1D1V:coeff}, we use the Taylor expansion.
%%
% If a function $\fs_{N,M}\in\Yv_{N,M}$ is known through its
% coefficients \eqref{eq:1D1V:coeff}, we would like to evaluate it at
% the new points $(\Xnewnm, \Vnewnm)$.
% %%
% To this end we use a Taylor expansions.
%%
%%We recall that, for a sufficiently smooth function $\Psi$, one has:
%%
For a sufficiently smooth function $\Psi$, we have that
\begin{align}\label{eq:1D1V:Taylor}
  &\Psi (\xs-\vs\Delta t,\vs + \Es_N(\xs) \Delta t) 
  = \Psi(\xs,\vs) - \vs\Delta t\,\frac{\partial\Psi}{\partial\xs}(\xs,\vs) 
  + \Es_N(x) \Delta t\,\frac{\partial\Psi}{\partial\vs}(\xs,\vs) \nonumber \\[3mm]
  &+\frac12\,(\vs\Delta t)^{2}\,\frac{\partial^{2}\Psi}{\partial\xs^{2}}(\xs,\vs)  
  -\vs\Es_N(\xs)\Delta t^{2}\,\frac{\partial^{2}\Psi}{\partial\xs\,\partial\vs}(\xs,\vs)
  +\frac12\,(\Es_N(x)\Delta t)^{2}\,\frac{\partial^{2}\Psi}{\partial\vs^{2}}(\xs,\vs) + \ldots .
\end{align}
%% 
%% where  the usual little-$o$ notation  is used.\\
%% 
Applying \eqref{eq:1D1V:Taylor} to $
\Psi(\xs,\vs)=\Bs_{i}^{(N)}(\xs)\,\Bs_{j}^{(M)}(\vs)$, when
$(\xs,\vs)=(\Xnewnm,\Vnewnm)$, %% $n=0,1,\ldots, N-1$, $ m=0,1,\ldots, M-1$,
is defined in \eqref{eq:1D1V:puntimossi},
% and approximating  up to the second order 
we obtain:
\begin{align}\label{eq:1D1V:eqclou2}
  &
  \Bs_{i}^{(N)}(\Xnewnm)\,\Bs_{j}^{(M)}(\Vnewnm) =
  \Bs_{i}^{(N)}(\Xoldnm)\,\Bs_{j}^{(M)}(\Voldnm) - 
  \incrXnm\, 
  \left[
    \frac{\partial \Bs_{i}^{(N)}}{\partial x}(\Xoldnm)
  \right]\,\Bs_{j}^{(M)}(\Voldnm)\nonumber\\[3mm]
  &\hskip1truecm		
  + \incrVnmb\,\Bs_{i}^{(N)}(\Xoldnm)\,\left[\frac{\partial\Bs_{j}^{(M)}}{\partial\vs}(\Voldnm)\right] +
  \frac12\,(\incrXnm)^{2}\,\left[\frac{\partial^{2}\Bs_{i}^{(N)}}{\partial\xs^{2}}(\Xoldnm)\right]\,\Bs_{j}^{(M)}(\Voldnm) \nonumber\\[3mm]
  &\hskip1truecm		 
  - \incrXnmVnmb\,\left[\frac{\partial\Bs_{i}^{(N)}}{\partial\xs}(\Xoldnm)\right]\,\left[\frac{\partial\Bs_{j}^{(M)}}{\partial\vs}(\Voldnm)\right]\nonumber
  \\[3mm]
  &\hskip1truecm
  + \frac12\,(\incrVnmb)^{2}\,\Bs_{i}^{(N)}(\Xoldnm)\,
  \left[\frac{\partial^{2}\Bs_{j}^{(M)}}{\partial\vs^{2}}(\Voldnm)\right] + \ldots .
  %\nonumber
  %\\[3mm]
  %&\hskip4truecm\quad n=0,1,\ldots, N-1,\,\,m=0,1,\ldots, M-1.
\end{align}
Using \eqref{eq:1D1V:Kron}, \eqref{eq:1D1V:Der1v} and
\eqref{eq:1D1V:Der2v}, we can rewrite \eqref{eq:1D1V:eqclou2} as:
\begin{align}\label{eq:1D1V:eqclou3}
  & \Bs_{i}^{(N)}(\Xnewnm)\,\Bs_{j}^{(M)}(\Vnewnm) =
  \delta_{in}\,\delta_{jm}  - 
  \incrXnm\,\delta_{jm}\,\ds_{ni}^{(N,1)} +	
  \incrVnmb\,\delta_{in}\,\ds_{mj}^{(M,1)}\nonumber\\[3mm]
  & \hskip1truecm	
  + \frac12\,(\incrXnm)^{2}\,\delta_{jm}\,\ds_{ni}^{(N,2)} 
  - \incrXnmVnmb\,\ds_{ni}^{(N,1)}\,\ds_{mj}^{(M,1)}\nonumber\\[3mm]
  & \hskip1truecm	
  + \frac12\,(\incrVnmb)^{2}\,\delta_{in}\,\ds_{mj}^{(M,2)} + \ldots .
  %\nonumber
  %\\[3mm]
  %& \hskip4truecm\quad n=0,1,\ldots, N-1,\,\,m=0,1,\ldots, M-1.
\end{align}
Substituting \eqref{eq:1D1V:eqclou3} in \eqref{eq:1D1V:fbasis}, we
obtain:
\begin{align}\label{eq:1D1V:eqclouDef}
  \fs_{N,M}(\Xnewnm,\Vnewnm) 
  & = \sum_{i=0}^{N-1}\,\sum_{j=0}^{M-1}\cs_{ij}\,\Bs_{i}^{(N)}(\Xnewnm)\,\Bs_{j}^{(M)}(\Vnewnm)\nonumber\\[3mm]
  & = \sum_{i=0}^{N-1}\,\sum_{j=0}^{M-1}\cs_{ij} 
  \Big(
  \delta_{in}\,\delta_{jm} -  
  \incrXnm \,\delta_{jm}\,\ds_{ni}^{(N,1)} +	
  \incrVnmb\,\delta_{in}\,\ds_{mj}^{(M,1)}\nonumber\\[3mm]
  & \hskip0.56truecm % \phantom{\sum_{i=0}^{N-1}\,\sum_{j=0}^{M-1}\cs_{ij}\Big(}\quad
  + \frac12\,(\incrXnm)^{2}\,\delta_{jm}\,\ds_{ni}^{(N,2)}
  - \incrXnmVnmb\,\ds_{ni}^{(N,1)}\,\ds_{mj}^{(M,1)}\nonumber\\[3mm]
  & \hskip0.56truecm % \phantom{\sum_{i=0}^{N-1}\,\sum_{j=0}^{M-1}\cs_{ij}\Big(}\quad
  + \frac12\,(\incrVnmb)^{2}\,\delta_{in}\,\ds_{mj}^{(M,2)} + \ldots.
  \Big)\nonumber
  \\[3mm]
  & = \cs_{nm} +\Delta t\left[ -\vs_m\sum_{i=0}^{N-1}\ds_{ni}^{(N,1)}\cs_{im} 
    + \Es_N(x_n)\sum_{j=0}^{M-1}\ds_{mj}^{(M,1)}\cs_{nj}\right]\nonumber\\
  & \hskip0.56truecm
  + \frac{(\Delta t)^2}{2}\left[ 
    \vs_m^2\sum_{i=0}^{N-1}\ds_{ni}^{(N,2)}\cs_{im} 
    - 2\vs_m\Es_N(x_n)\sum_{i=0}^{N-1}\sum_{j=0}^{M-1}\ds_{ni}^{(N,1)}\ds_{mj}^{(M,1)}\cs_{ij} 
  \right.
  \nonumber\\
  & \phantom{= +\frac{(\Delta t)^2}{2} }\hskip0.5truecm
  \left.
    + (\Es_N(\xs_n))^2\sum_{j=0}^{M-1}\ds_{mj}^{(M,2)}\cs_{nj}
  \right] + \ldots .
  %,\nonumber
  %\\[3mm]
  %&	
  %n=0,1,\ldots, N-1,\,\,m=0,1,\ldots, M-1.
\end{align}
In compact form we can write:
\begin{align}\label{eq:1D1V:coeffTaylor}
  \displaystyle
  & \fs_{N,M} (\Xnewnm,\Vnewnm) = \cs_{nm} + \sum_{s=1}^\infty \sum_{r=0}^s\frac{(-1)^s}{r! (s-r)!}
  \left(\mathcal{I}_{nm}^{\,r}\mathcal{J}_{nm}^{s-r}\sum_{i=0}^{N-1}\sum_{j=0}^{M-1}\ds_{ni}^{(N,r)}\ds_{mj}^{(M,s-r)}\cs_{ij}\right),
  %\nonumber
  %\\[3mm]
  %& \hskip4truecm	
  %\quad n=0,1,\ldots, N-1,\,\,m=0,1,\ldots, M-1,
\end{align}
where we set ${\cal I}_{nm} = x_n -\Xnewnm$ and ${\cal J}_{nm} = v_m
-\Vnewnm$. %% (for $n=0,1,\ldots, N-1$, $ m=0,1,\ldots, M-1$).
%%
% It is clear that one can truncate the summation with respect to $s$ at
% any integer $S$, so that obtaining a reminder of order $(\Delta
% t)^{S+1}$. 
% %%
% Note also that the differentiation in the variables $x$ and $v$ can be
% computed exactly, through the multiplication of the corresponding
% derivative matrices. 
% %%
% Therefore, if we assume that in \eqref{eq:1D1V:coeffTaylor} the integer $s$
% can range from 1 to infinity, no approximation has been introduced so
% far.
%%
Finally, we truncate the summation with respect to $s$ at the integer
$S\ge1$ to have a remainder term of order $(\Delta t)^{S+1}$.
The differentiation in the variables $x$ and $v$ can be computed
exactly by multiplying the corresponding derivative matrices.
Therefore, no approximation is introduced if we assume that the
integer $s$ can range from 1 to infinity in
\eqref{eq:1D1V:coeffTaylor}.

\subsection{Three-dimensional extension}

The three-dimensional extension of~\eqref{eq:1D1V:coeffTaylor} is
straightforward by using the multi-index notation.
To this end, we consider all indices $n,m,i,j,s,r$
in~\eqref{eq:1D1V:coeffTaylor} as multi-indices of order three.
More precisely, $n$ is the triplet of nonnegative integers$(n_1,n_2,n_3)$ and $\abs{n}=n_1+n_2+n_3$ is the order of $n$.
The position vector is thus given by $\xv=(\xs^{1},\xs^{2},\xs^{3})$,
and, a similar notation holds for the velocity position vector
$\vv=(\vs^{1},\vs^{2},\vs^{3})$.
A space vector subindexed by $n$ has to be interpreted as the grid
point $\xv_{n}=(\xs^{1}_{n_1},\xs^{2}_{n_2},\xs^{3}_{n_3})$; a
velocity vector subindexed by $m$ has to be interpreted as the grid
point $\vv_{m}=(\vs^{1}_{m_1},\vs^{2}_{m_2},\vs^{3}_{m_3})$.
Consistently, we also have the double-subindexed vectors
$\xvt_{nm}=(\xst^1_{n_1m_1},\xst^2_{n_2m_2},\xst^3_{n_3m_3})$ and
$\vvt_{nm}=(\vst^1_{n_1m_1},\vst^2_{n_2m_2},\vst^3_{n_3m_3})$.
We use the standard notation
$(\wv)^r=(\ws^1)^{r_1}(\ws^2)^{r_2}(\ws^3)^{r_3}$ for any given
three-dimensional vector $\wv=(\ws^1,\ws^2,\ws^3)$ and multi-index
$r=(r_1,r_2,r_3)$, and we denote the partial derivatives of order
$\abs{r}$ of a generic function $\gs(\xv)$ determined by the multi-index
$r$ as:
\begin{align*}
  \frac{\partial^{\abs{r}}}{\partial\xv^{r}}\gs(\xv) =
  \frac{\partial^{r_1}}{\partial\xs^{1,r_1}} 
  \frac{\partial^{r_2}}{\partial\xs^{2,r_2}} 
  \frac{\partial^{r_3}}{\partial\xs^{3,r_3}} \gs(\xv).
\end{align*}
A similar relation holds for the partial derivatives along $\vv$.
Finally, the three-dimensional basis functions are given by the tensor
product of the one-dimensional basis functions:
\begin{align*}
  \Bv_{i}^{(N)}(\xv) 
  = \Bs_{i_1}^{(N)}(\xs^1) \Bs_{i_2}^{(N)}(\xs^2) \Bs_{i_3}^{(N)}(\xs^3) ,
  \qquad i_1,i_2,i_3=0,\ldots,N-1.
\end{align*}

Now, the three-dimensional version of
equation~\eqref{eq:1D1V:coeffTaylor} becomes:
\begin{multline}
  \displaystyle
  \fs_{N,M} (\Xvnewnm,\Vvnewnm) = \cs_{nm} + \sum_{\abs{s}=1}^{\infty} \sum_{\abs{r}=0}^{\abs{s}} \frac{(-1)^{\abs{s}}}{\abs{r}!\abs{s-r}!}
  \left( \big({\boldsymbol{\cal I}}_{nm}\big)^{\,r}\,\big({\boldsymbol{\cal J}}_{nm}\big)^{s-r}\sum_{\abs{i}=0}^{N-1}\sum_{\abs{j}=0}^{M-1}\mathbf{d}_{ni}^{(N,r)}\mathbf{d}_{mj}^{(M,s-r)}\cs_{ij}\right),\nonumber
  \\[3mm]
  \hskip4truecm	
  \quad \abs{n}=0,1,\ldots, N-1,\,\,\abs{m}=0,1,\ldots, M-1,
  %\label{eq:3D3V:coeffTaylor}
\end{multline}
where we set ${\boldsymbol{\cal I}}_{nm}=\xv_n-\Xvnewnm$,
$({\boldsymbol{\cal I}}_{nm})^{\,r}=\big(\xv_n-\Xvnewnm)^{r}$,
%% $=(\xs_n^{1}-\Xnewnm^{1})^{r_1}\,(\xs_n^{2}-\Xnewnm^{2})^{r_2}\,(\xs_n^{3}-\Xnewnm^{1})^{r_3}$,
%%
%and, similarly, 
${\boldsymbol{\cal J}}_{nm}=\vv_m-\Vvnewnm$,
$({\boldsymbol{\cal J}}_{nm})^{s-r}=(\vv_m-\Vvnewnm)^{s-r}$;
%%
%%$= \big(\vs_m^1-\Vnewnm^1\big)^{s_1-r_1}\,\big(\vs_m^2-\Vnewnm^2\big)^{s_2-r_2}\,\big(\vs_m^3-\Vnewnm^3\big)^{s_3-r_3}$;
%%
%% $\abs{n}=0,1,\ldots, N-1$, $\abs{m}=0,1,\ldots, M-1$,
the partial derivatives of the three-dimensional basis functions are
given by
\begin{align}
  \mathbf{d}_{ni}^{(N,s)} = \frac{\partial^{\abs{s}}\Bv_{i}^{(N)}}{\partial\xv^{s}}(\xv_{n})
  \qquad\textrm{and}\qquad
  \mathbf{d}_{mj}^{(M,s)} = \frac{\partial^{\abs{s}}\Bv_{j}^{(M)}}{\partial\vv^{s}}(\vv_{m}).
  \label{eq:3D3V:Der}
\end{align}
All considerations at the end of the previous section are still true
here.

%%%%%%%%%%%%%%%%%%%%%%%%%%%%%%%%%%%%%%%%%%%%%%%%%%%%%%%%%%
\section{Time discretization}
\label{sec:time:discretization}
%%%%%%%%%%%%%%%%%%%%%%%%%%%%%%%%%%%%%%%%%%%%%%%%%%%%%%%%%%
Given the time instants $t^{k}=k \Delta t= k\,{T}\slash{K}$ for any
integer $k=0,1,\ldots,K$, we consider here the full approximation of
the solution fields $(f,E)$ of the 1D-1V Vlasov-Poisson problem
\eqref{eq:1D1V:V}, \eqref{eq:1D1V:Vci}, \eqref{eq:1D1V:EMaxwell},
\eqref{eq:1D1V:charge-density}:
\begin{align}\label{eq:1D1V:fapproxt}
  \left(
    \fs_{N,M}^{(k)}(\xs,\vs),\,\Es^{(k)}_N(\xs) 
  \right) \simeq 
  \left(
    \fs(t^{k},\xs,\vs),\,\Es(t^{k},\xs) 
  \right) ,
  %,\nonumber \\[3mm]& \hskip5truecm  k=0,1,\ldots,K,\,\,
  \qquad\xs\in\Omega_{x},\,\,\vs\in\Omega_{v},
\end{align}
where the function $\fs_{N,M}^{(k)}$ belongs to $\Yv_{N,M}$ and the
function $\Es_N^{(k)}$ belongs to $\Xv_N$.
Taking into account \eqref{eq:1D1V:charge-density}, we define:
\begin{equation}\label{eq:1D1V:charge-densityapproxt}
  \displaystyle
  \rho^{(k)}_N(\xs)  =
  \int_{\Omega_{v}}\fs^{(k)}_{N,M}(\xs,\vs)\,\dv
  \simeq \rho(t^{(k)},\xs).
  %%\qquad k=0,1,\ldots,K,\,\,\xs\in\Omega_{x}.
\end{equation}

At any timestep $k$, 
%$k=0,1,\ldots, K$, 
we evaluate $\fs^{(k)}_{N,M}$ in the following way:
\begin{align}\label{eq:1D1V:fbasisk}
  \fs^{(k)}_{N,M} (\xs,\vs) 
  = \sum_{i=0}^{N-1}\,\sum_{j=0}^{M-1} \cs^{(k)}_{ij}\,\Bs_{i}^{(N)}(\xs)\,\Bs_{j}^{(M)}(\vs) , 
  %\qquad k=0,1,\ldots, K,\,\,
  %\qquad\xs\in\Omega_{x},\,\,\vs\in\Omega_{v},
\end{align}
where
\begin{equation}\label{eq:1D1V:coeffK}
  \cs^{(k)}_{ij} = \fs^{(k)}_{N,M} (\xs_{i},\vs_{j}).
  %% , \qquad k=0,1,\ldots,K,\,\, 
 % \qquad i=0,1,\ldots, N-1,\,\,j=0,1,\ldots, M-1.
\end{equation}
In particular, at time $t=0$, we use the initial condition for $f$
(see equation \eqref{eq:1D1V:Vci}) by setting
\begin{equation}\label{eq:1D1V:coeff0}
  \cs^{(0)}_{ij} = \fs(0,\xs_{i},\vs_{j})=\bar{\fs}(\xs_{i},\vs_{j}).
  %% , \qquad  i=0,1,\ldots, N-1,\,\,j=0,1,\ldots, M-1.
\end{equation}

% For an initial distribution $\Es_N^{(0)}$ of the electric field (that
% can be a suitable projection in $\Xv_N$ of the initial datum
% $\Es$  in \eqref{eq:1D1V:EMaxwell} for $t=0$), we advance at step $k$ according
% to the following scheme.
%% 
If we suppose that $\Es_N^{(k)}$ is given at step $k$, we first define
(take $\tau=t-\Delta t$ in \eqref{eq:1D1V:char1}):
\begin{align}
  \Xnewnm &= \Xoldnm - \incrXnm,  \nonumber\\[2mm]
  \Vnewnm &= \Voldnm + \Es_N^{(k)}(\xs_n)\Delta t.
  %% ,  \nonumber\\[2mm]
  %% &\hskip2truecm \quad k=0,1,\ldots,K,\,\,n=0,1,\ldots, N-1,\,\,m=0,1,\ldots, M-1. 
  \label{eq:1D1V:2charnm}
\end{align}

Since the solution $\fs$ of the Vlasov-Poisson system is expected to
be constant along the characteristics, the most straightforward method
is obtained by advancing the coefficients of $f_{N,M}\simeq f$ as
follows
\begin{align}
  \cs^{(k+1)}_{nm} 
  = \fs^{(k)}_{N,M}(\Xnewnm,\Vnewnm) 
  = \sum_{i=0}^{N-1}\,\sum_{j=0}^{M-1}\cs^{(k)}_{ij}\,\Bs_{i}^{(N)}(\Xnewnm)\,\Bs_{j}^{(M)}(\Vnewnm),
  \label{eq:1D1V:favanz}
\end{align}
%%
% \begin{align}
%   \cs^{(k+1)}_{nm} 
%   &= \fs^{(k)}_{N,M}(\Xnewnm,\Vnewnm) 
%   = \sum_{i=0}^{N-1}\,\sum_{j=0}^{M-1}\cs^{(k)}_{ij}\,\Bs_{i}^{(N)}(\Xnewnm)\,\Bs_{j}^{(M)}(\Vnewnm),\nonumber\\[3mm]
%   &\hskip2truecm	
%   k=0,1,\ldots, K-1,\,\,n=0,1,\ldots, N-1,\,\,m=0,1,\ldots, M-1,
%   \label{eq:1D1V:favanz}
% \end{align}
where we used representation~\eqref{eq:1D1V:fbasisk}.
This states that %, for any time-step $k$, $k=0,1,\ldots, K$,
the value of $\fs_{N,M}^{(k+1)}$, at the grid points and timestep
$(k+1)\Delta t$, is assumed to be equal to the previous value at time
$k\Delta t$, recovered by going backwards along the characteristics.
Technically, in \eqref{eq:1D1V:2charnm} we should use
$\Es_N^{(k+1)}(\xs_n)$ instead of $\Es_N^{(k)}(\xs_n)$, thus arriving
at an implicit method.
However, the distance between these two quantities is of the order of
$\Delta t$, so that the replacement has no practical effects on the
accuracy of the first-order method.
For higher order schemes, things must be treated more carefully.

Between each step $k$ and the successive one, we need to update the
electric field.
This can be done as suggested here below.

Let $t^{k}$ be fixed.
%%Let $t^{k}$, $k=0,1,\ldots, K$, be fixed.
%% 
Using the Gaussian quadrature formula \eqref{eq:1D1V:fdq} in
\eqref{eq:1D1V:charge-densityapproxt} and~\eqref{eq:1D1V:coeffK} we
write:
\begin{equation}\label{eq:1D1V:charge-densityapproxt2}
  \displaystyle
  \rho^{(k)}_N(\xs_{i})
  = \frac{2\pi}{M}\,\sum_{j=0}^{M-1}\,\fs^{(k)}_{N,M}(\xs_{i},\vs_{j})
  = \frac{2\pi}{M}\,\sum_{j=0}^{M-1}\,\cs_{ij}^{(k)} . 
  %\qquad i=0,1,\ldots, N-1.
\end{equation}
Indeed, it is possible to compute $\rho^{(k)}_N(\xs)$
%% $\xs\in\Omega_{x}=[0,2\pi[$, 
by using the Fourier series:
\begin{equation}\label{eq:1D1V:FourierSeries}
  \displaystyle
  \rho^{(k)}_N(\xs) 
  = 1 + \sum_{n=1}^{N/2} 
  \left[ 
    \hat{\as}_{n}^{(k)}\,\cos(n\xs) + \hat{\bs}_{n}^{(k)}\,\sin(n\xs) 
  \right],
  %% \quad\xs\in\Omega_{x}=[0,2\pi[,
\end{equation}
where the discrete Fourier coefficients $\hat{\as}_{n}^{(k)}$ and
$\hat{\bs}_{n}^{(k)}$ are determined, for $n=1,2,\ldots,N/2$, by the
following formulas:
\begin{align}\label{eq:1D1V:coeffFourierSeries}
  & \hat{\as}_{n}^{(k)}
  = \frac{1}{\pi}\int_0^{2\pi}\rho^{(k)}_N(\xs)\cos(n\xs)\,\dx 
  \simeq \frac{2}{N}\sum_{l=0}^{N-1}\rho^{(k)}_N(\xs_{l})   
  \cos\left(\frac{2nl}{N}\,\pi \right),\nonumber\\
  & \hat{\bs}_{n}^{(k)} 
  = \frac{1}{\pi}\int_0^{2\pi}\rho^{(k)}_N(\xs)\sin(n\xs)\,\dx
  \simeq \frac{2}{N}\sum_{l=0}^{N-1}\rho^{(k)}_N(\xs_{l})
  \sin\left(\frac{2nl}{N}\,\pi\right). 
\end{align}
Actually, for $n$ strictly smaller than $N/2$, the symbol ``$\simeq$''
can be replaced by the symbol ``$=$''.
%\footnote{
%  \marco{
%    Anche a me sembra che ci sia
%      l'uguale. Se ho capito bene, questo succede perch\'e regola di
%      integrazione~\eqref{eq:1D1V:fdq} che permette di definire la
%      $\rho^{(k)}_N(\xs_{s})$ integrando la $\fs^{(k)}_{N,M}(\xs,\vs)$
%      in $v$ \`e esatta per le funzioni $\cos$ e $\sin$ con cui si
%      definiscono le basis functions $\Bs_{i}^{(N)}(\xs)$ e
%      $\Bs_{i}^{(M)}(\vs)$. Tuttavia, non capisco bene perch\'e il
%      caso $k=N/2$ sarebbe diverso. Comunque possiamo lasciare il
%      simbolo $\simeq$ e mettere un breve commento. Vicino
%      alla~\eqref{eq:1D1V:fdq} bisognerebbe, forse, estendere un
%      pochino il commento.
%  }
%}
%% 
%  \marco{
%    Note that equation~\eqref{eq:1D1V:coeffFourierSeries} is an equality
%    for $k=N/2$ due to the exactness of the integration
%    rule~\eqref{eq:1D1V:fdq} on $\Xv_N$ and $\Yv_{N,M}$.}\\

Using equation \eqref{eq:1D1V:FourierSeries} and equation
\eqref{eq:1D1V:EMaxwell} at $t=t^{k}$, 
%% $k=0,1,\ldots, K,$
we conclude that:
\begin{equation}
  \Es^{(k)}_N(\xs) = 
  -\sum_{n=1}^{N/2} \frac{1}{n}
  \left[ 
    \hat{\as}_{n}^{(k)}\,\sin(n\xs) - \hat{\bs}_{n}^{(k)}\,\cos(n\xs)
  \right],
  %% \quad x\in\Omega_{x} = [0,2\pi[,
  \label{eq:1D1V:FourierSeriesE}
\end{equation}
which satisfies (as requested in \eqref{eq:1D1V:intrhoconst}):
\begin{equation}\label{eq:1D1V:mediaEnulla}
  \int_{0}^{2\pi}\,\Es^{(k)}_N(\xs)\,\dx = 0.
  %% , \quad\xs\in\Omega_{x}=[0,2\pi[.
\end{equation}
Finally, from \eqref{eq:1D1V:coeffFourierSeries}, using a
  standard trigonometric formula
  and~\eqref{eq:1D1V:charge-densityapproxt2}, we find that:
\begin{align}\label{eq:1D1V:FourierSeriesEnodi}
  &\Es^{(k)}_N(\xs_i)
  = -\sum_{n=1}^{N/2}\frac{1}{n}\left[ \hat{\as}_{n}^{(k)}\,\sin\left( \frac{2ni}{N}\,\pi \right) - \hat{\bs}_{n}^{(k)}\,\cos\left( \frac{2ni}{N}\,\pi\right)\right]
  \nonumber\\
  \hskip1truecm
  & \simeq \frac{2}{N}\sum_{n=1}^{N/2}\frac{1}{n}\sum_{s=0}^{N-1} \rho^{(k)}_N(\xs_s)
  \left[ 
    \sin\left( \frac{2sn}{N}\,\pi \right) 
    \cos\left( \frac{2in}{N}\,\pi \right) -
    \sin\left( \frac{2in}{N}\,\pi \right) 
    \cos\left( \frac{2sn}{N}\,\pi \right) 
  \right]
  \nonumber\\
  \hskip1truecm
  & = \frac{2}{N}\sum_{n=1}^{N/2}\frac{1}{n}\sum_{s=0}^{N-1} \rho^{(k)}_N(\xs_s) \sin\left( \frac{2(s-i)n}{N}\,\pi \right) \nonumber\\
  \hskip1truecm 
  & = \frac{4\pi}{NM}\sum_{n=1}^{N/2}\frac{1}{n}\sum_{s=0}^{N-1}\sum_{j=0}^{M-1}\cs^{(k)}_{ij} \sin\left( \frac{2(s-i)n}{N}\,\pi \right).
  %,\quad  i=0,1,\ldots, N-1.
\end{align}

By computing the direction of the characteristic lines according to
\eqref{eq:1D1V:2charnm}, the scheme turns out to be only first-order
accurate in $\Delta t$.
Consequently, it is sufficient to stop the development
\eqref{eq:1D1V:coeffTaylor} at $s=1$.
In this way, \eqref{eq:1D1V:favanz} is replaced by:
\begin{align}
  \cs^{(k+1)}_{nm} 
  = \cs_{nm}^{(k)} + \Delta t\,\Phi_{nm}^{(k)},
  \label{eq:1D1V:favanza}
\end{align}
% \begin{align}
%   \cs^{(k+1)}_{nm} 
%   &= \cs_{nm}^{(k)} + \Delta t\,\Phi_{nm}^{(k)} \ ,\nonumber\\[3mm]
%   &\hskip1truecm	
%   \quad k=0,1,\ldots,K-1,\,\,n=0,1,\ldots, N-1,\,\,m=0,1,\ldots, M-1,
%   \label{eq:1D1V:favanza}
% \end{align}
where
\begin{align}\label{eq:1D1V:favanzadef}
  \Phi_{nm}^{(k)} 
  =
  - \vs_m\sum_{i=0}^{N-1}\ds_{ni}^{(N,1)}\cs_{im}^{(k)} 
  + \Es_N^{(k)}(\xs_n)\sum_{j=0}^{M-1}\ds_{mj}^{(M,1)}\cs_{nj}^{(k)}.
\end{align}
%%
% \begin{align}\label{eq:1D1V:favanzadef}
%   \Phi_{nm}^{(k)} 
%   &=
%   - \vs_m\sum_{i=0}^{N-1}\ds_{ni}^{(N,1)}\cs_{im}^{(k)} 
%   + \Es_N^{(k)}(\xs_n)\sum_{j=0}^{M-1}\ds_{mj}^{(M,1)}\cs_{nj}^{(k)} \nonumber\\
%   & \hskip1truecm	
%   \quad k=0,1,\ldots, K-1,\,\,n=0,1,\ldots, N-1,\,\,m=0,1,\ldots, M-1.
% \end{align}

%% In view of solving the non-homogeneous 1D-1V Vlasov-Poisson equation:

Consider a sufficiently regular function $\gs(t,\xs,\vs)$, which
is defined on $\Omega$ for every
$t\in[0,T]$.
To solve the non-homogeneous Vlasov equation:
%%% 
\begin{equation}\label{eq:1D1V:vlasovno}
  \displaystyle
  \frac{\partial\fs}{\partial t} + \vs\ \frac{\partial\fs}{\partial\xs} 
  - \Es(t,\xs)\ \frac{\partial\fs}{\partial\vs} = \gs, 
  %% \quad t\in (0,T],\,\,\xs \in \Omega_{x},\,\, v\in \Omega_{v},
\end{equation}
%% 
%% where $\gs(t,\xs,\vs)$, $t\in(0,T]$, $\xs\in\Omega_{x}$,
%% $\vs\in\Omega_{v}$, is a given sufficiently regular right-hand side,
we modify \eqref{eq:1D1V:favanza} as follows:
\begin{align}\label{eq:1D1V:favanzano}
  %% & 
  \cs^{(k+1)}_{nm} =
  \cs_{nm}^{(k)} + \Delta t\,\Phi_{nm}^{(k)} + \Delta t\,\gs(t^k,\xs_n,\vs_m),
%   ,\nonumber\\[3mm]
%   &\hskip1truecm	
%   \quad k=0,1,\ldots, K-1,\,\,n=0,1,\ldots, N-1,\,\,m=0,1,\ldots, M-1.
\end{align}
where $\Phi_{nm}^{(k)}$
%% , $k=0,1,\ldots, K-1$, $ n=0,1,\ldots, N-1$, $m=0,1,\ldots, M-1$,
is the same as in \eqref{eq:1D1V:favanzadef}.
This is basically a forward Euler iteration. 

As expected from an explicit method, the parameter $\Delta t$ must
satisfy a suitable CFL condition, which is easily obtained by
requiring that the point $(\tilde{\xs}_{nm}, \tilde{\vs}_{nm})$ falls
inside the box $]x_{n-1}, x_{n+1}[\times ]v_{m-1}, v_{m+1}[$.  
From \eqref{eq:1D1V:2charnm}, a sufficient restriction is given by:
\begin{equation}\label{eq:1D1V:CFL}
  \Delta t 
  \leq 2\pi \Big( N \max_m \vert v_m\vert + M \max_n \vert E^{(k)}_N(x_n)\vert\Big)^{-1}.
\end{equation}
By inequality \eqref{eq:1D1V:stimad}, this ensures that the term
$\Delta t \, \Phi_{nm}^{(k)}$ in \eqref{eq:1D1V:favanza} is of the same
order of magnitude as $\cs_{nm}^{(k)}$.  
%%
% We would like to use more potentialities of the expansion
% \eqref{eq:1D1V:coeffTaylor}, so that we will analyze more accurate
% time-marching schemes in the next section.

We will better use the potentialities of
expansion~\eqref{eq:1D1V:coeffTaylor} in the next section to design
more accurate time-marching schemes.

%%%%%%%%%%%%%%%%%%%%%%%%%%%%%%%%%%%%%%%%%%%%%%%%%%%%%%%%%%
\section{More advanced time discretizations}
\label{sec:more:advanced:time:discretizations}
%%%%%%%%%%%%%%%%%%%%%%%%%%%%%%%%%%%%%%%%%%%%%%%%%%%%%%%%%%
%% 
A straightforward way to increase the time accuracy is to use a
higher-order time-marching scheme.
To this end, we consider the second-order accurate two-step explicit
Backward Differentiation Formula (BDF).
With the notation in \eqref{eq:1D1V:favanz},
\eqref{eq:1D1V:favanzadef} and \eqref{eq:1D1V:favanzano},  given the time instants $t^{k}=k \Delta t= k\,{T}\slash{K}$,  $k=0,1,\ldots,K$, we have:
\begin{align}\label{eq:1D1V:favanzabdf2}
  %& 
  \fs_{N,M}^{(k+1)}(\xs_n,\vs_m) 
  = {{\frac43}} \fs_{N,M}^{(k)}  (\tilde{\xs}_{nm}, \tilde{\vs}_{nm}) 
  - {{\frac13}} \fs_{N,M}^{(k-1)}(\tilde{\tilde{\xs}}_{nm}, \tilde{\tilde{\vs}}_{nm}) 
  + {{\frac23}} \Delta t\,\gs(t^{k+1},\xs_n,\vs_m), 
%   \nonumber\\[3mm]
%   & \hskip1truecm	
%   \quad k=1,2,\ldots, K-1,\,\,n=0,1,\ldots, N-1,\,\,m=0,1,\ldots, M-1,
\end{align}
where, based on \eqref{eq:1D1V:2charnm},
$(\tilde{\xs}_{nm},\tilde{\vs}_{nm})$ 
%% $n=0,1,\ldots,N-1$, $m=0,1,\ldots,M-1$,
%%
is the point obtained from $(\xs_n,\vs_m)$ going
back of one step $\Delta t$ along the characteristic lines.
Similarly, 
%% for $n=0,1,\ldots,N-1$, $m=0,1,\ldots,M-1$,
the point $(\tilde{\tilde{\xs}}_{nm},\tilde{\tilde{\vs}}_{nm})$ is
obtained by going two steps back along the characteristic lines
(replace $\Delta t$ with $2\Delta t$ in \eqref{eq:1D1V:2charnm}).
Note that if $\gs=0$, it turns out that $\fs_{N,M}$ is constant along
the characteristic lines.

%% By approximating at the first order the above values as follows:
%%
The first-order accurate approximation of the above values for any
integer $k=1,2,\ldots,K-1$ reads as
\begin{align}\label{eq:1D1V:appf}
  & \fs_{N,M}^{(k)}  (\tilde{\xs}_{nm},\tilde{\vs}_{nm}) \simeq \cs_{nm}^{(k)} + \Delta t\,\Phi_{nm}^{(k)}, \nonumber\\[3mm]
  & \fs_{N,M}^{(k-1)}(\tilde{\tilde{\xs}}_{nm},\tilde{\tilde{\vs}}_{nm}) \simeq \cs_{nm}^{(k-1)} + 2\Delta t\,\Phi_{nm}^{(k-1)},  
  % \nonumber\\[3mm]
%   & \hskip1truecm	
%   \quad k=1,2,\ldots, K-1,\,\,n=0,1,\ldots, N-1,\,\,m=0,1,\ldots, M-1,
\end{align}
and, in terms of the coefficients, we end up with the scheme:
\begin{align}\label{eq:1D1V:favanzabdf2co}
  & \cs_{nm}^{(k+1)} =
  {\frac43} \Big( \cs_{nm}^{(k)}   +  \Delta t\,\Phi_{nm}^{(k)}  \Big) -
  {\frac13} \Big( \cs_{nm}^{(k-1)} + 2\Delta t\,\Phi_{nm}^{(k-1)} \Big) + {\frac23}\Delta t\,\gs(t^{k+1},\xs_n,\vs_m) \nonumber\\
  &\hskip1truecm = {{\frac43}} \cs_{nm}^{(k)} - {{\frac13}} \cs_{nm}^{(k-1)} 
  + {{\frac23}}\Delta t \left[  
    -\vs_m\sum_{i=0}^{N-1} \ds_{ni}^{(N,1)} (2\cs_{im}^{(k)} - \cs_{im}^{(k-1)}) \right. \nonumber\\
  & \hskip1truecm\qquad\left. + \Es_N^{(k)}(x_n)\sum_{j=0}^{M-1} \ds_{mj}^{(M,1)} (2\cs_{nj}^{(k)} - \cs_{nj}^{(k-1)})\right]+
  {\frac23}\Delta t\,\gs(t^{k+1},\xs_n,\vs_m).
%   ,  \nonumber\\[3mm]
%   & \hskip1truecm	
%   \quad k=1,2,\ldots, K-1,\,\,n=0,1,\ldots, N-1,\,\,m=0,1,\ldots, M-1.
\end{align}
This method is second-order accurate in $\Delta t$ as will be shown by the numerical experiments
of Section~\ref{sec:numerical:experiments}.
% we are going to check with some numerical experiments in
% section~\ref{sec:numerical:experiments}.

In the same fashion, a third-order BDF scheme is obtained by setting:
\begin{align}\label{eq:1D1V:favanzabdf3co}
  & \cs_{nm}^{(k+1)} = \frac{18}{11} \Big( \cs_{nm}^{(k)} + \Delta t\,\Phi_{nm}^{(k)}\Big) - \frac{9}{11} \Big( \cs_{nm}^{(k-1)} + 2\Delta t\,\Phi_{nm}^{(k-1)} \Big)\nonumber\\[3mm]
  & \hskip1truecm\qquad + \frac{2}{11} \Big( \cs_{nm}^{(k-2)} + 3\Delta t\,\Phi_{nm}^{(k-2)} \Big) + \frac{6}{11} \Delta t\,\gs(t^{k+1},\xs_n,\vs_m), 
%   \nonumber\\[3mm]
%   & \hskip1truecm	
%   \quad k=2,3,\ldots, K-1,\,\,n=0,1,\ldots, N-1,\,\,m=0,1,\ldots, M-1.
\end{align}
where, now, the time index $k$ ranges from $2$ to $K-1$.

The further question is to see if it is possible to propose an
explicit one-step second-order scheme. 
The problem is delicate, since it is not enough to consider the
quadratic terms of the expansion in \eqref{eq:1D1V:eqclouDef}.
It is also necessary to work with a better representation of the
characteristic lines, such as that in \eqref{eq:1D1V:char2}, where, we
set $\tau=t-\Delta t$.
This time for $ k=0,1,\ldots, K$, we propose:
\begin{align}
  & \Xnewnm = \Xoldnm - \incrXnm -\frac12 \Es_N^{(k+1)}(\xs_n) \Delta t^2 ,  \nonumber\\[3mm]
  & \Vnewnm = \Voldnm + \Es_N^{(k+1)}(\xs_n) \Delta t - \frac12
  \left( 
    \frac{\partial\Es_N^{(k+1)}}{\partial t}(\xs_n) + \vs_m \frac{\partial\Es_N^{(k+1)}}{\partial\xs}(\xs_n)
  \right)\Delta t^2 , 
  \label{eq:1D1V:2charnmho}
 % & \hskip2truecm k=0,1,\ldots, K,\,\,n=0,1,\ldots, N-1,\,\,m=0,1,\ldots, M-1 ,
\end{align}
that corresponds to an implicit method. 
We apply the correction:
\begin{equation}
  \Es_N^{(k+1)} \simeq \Es_N^{(k)} + \frac{\partial\Es_N^{(k)}}
  {\partial t} \Delta t.
  %% , \quad k=0,1,\ldots,K,\label{eq:1D1V:corre}
\end{equation}
Thus, up to errors of the second order, we can modify \eqref{eq:1D1V:2charnmho}
as follows:
\begin{align}
  & \XnewnmIIord = \Xoldnm -\incrXnm -\frac12 \Es_N^{(k)}(\xs_n) \Delta t^2 = \xs_n - \hat{\mathcal{I}}_{nm}, \nonumber\\[2mm]
  & \VnewnmIIord = \Voldnm + \Es_N^{(k)}(\xs_n) \Delta t 
  + \frac12 \left( 
    \frac{\partial\Es_N^{(k)}}{\partial t}(\xs_n) - \vs_m \frac{\partial\Es_N^{(k)}}{\partial\xs}(\xs_n)
  \right) \Delta t^2 
  = \vs_m - \hat{\mathcal{J}}_{nm},
%   \nonumber\\[3mm]
%   & \hskip2truecm k=0,1,\ldots, K,\,\,n=0,1,\ldots, N-1,\,\,m=0,1,\ldots, M-1 ,
  \label{eq:1D1V:2charnmho2}
\end{align}
where, for brevity of notation, we introduced the two quantities
$\hat{\mathcal{I}}_{nm}$ and $\hat{\mathcal{J}}_{nm}$.
%% $ n=0,1,\ldots, N-1$, $m=0,1,\ldots, M-1.$
%% 
The partial derivative of $\Es_N^{(k)}$ with respect to $\xs$ is
available and recoverable from $\rho_N^{(k)}$ (see
\eqref{eq:1D1V:charge-densityapproxt}).
Regarding the time derivative, we can recall \eqref{eq:1D1V:eder} and set:
\begin{equation}\label{eq:1D1V:detapprox}
  \displaystyle
  \frac{\partial\Es_N^{(k)}}{\partial t}(\xs_n) \simeq \int_{\Omega_{v}} \vs\fs^{(k)}_{N,M}(\xs_n,\vs)\,\dv.
  %% , \quad k=0,1,\ldots, K,\,\,n=0,1,\ldots, N-1.
\end{equation}
Successively, the integral on the right-hand side is approximated by
quadrature. 
Once the point $(\Xnewnm , \Vnewnm)$ has been localized with
sufficient detail, one can apply the correction of the coefficients as
suggested by~\eqref{eq:1D1V:eqclouDef} thus neglecting the terms of
order higher than $\Delta t^2$.
In the new situation we have (see also \eqref{eq:1D1V:coeffTaylor} for $s=2$):
\begin{align}\label{eq:1D1V:eqcloumod}
  & \cs^{(k+1)}_{nm} = \cs^{(k)}_{nm} 
  - \hat{\mathcal{I}}_{nm}\sum_{i=0}^{N-1}\ds_{ni}^{(N,1)}\cs^{(k)}_{im} 
  - \hat{\mathcal{J}}_{nm}\sum_{j=0}^{M-1}\ds_{mj}^{(M,1)}\cs^{(k)}_{nj} \nonumber\\
  & \hskip0.5truecm	
  \displaystyle
  + {{\frac12}}\hat{\mathcal{I}}_{nm}^2\sum_{i=0}^{N-1}\ds_{ni}^{(N,2)}\cs^{(k)}_{im} 
  + \hat{\mathcal{I}}_{nm}\hat{\mathcal{J}}_{nm}\sum_{i=0}^{N-1}\sum_{j=0}^{M-1}\ds_{ni}^{(N,1)}\ds_{mj}^{(M,1)}\cs^{(k)}_{ij} 
  + {{\frac12}}\hat{\mathcal{J}}_{nm}^2\sum_{j=0}^{M-1}\ds_{mj}^{(M,2)}\cs^{(k)}_{nj}.
%   , \nonumber
%   \\[3mm]
%   &		\hskip1truecm	
%   k=0,1,\ldots, K-1,\,\,  n=0,1,\ldots, N-1,\,\,m=0,1,\ldots, M-1.
\end{align}

For the non-homogeneous equation \eqref{eq:1D1V:vlasovno}, suitable
adjustments are required to preserve the quadratic convergence.
Indeed, in order to handle the right-hand side, we suggest to use the
trapezoidal rule by defining:
\begin{align}\label{eq:1D1V:gtrapezi}
  %&
  \Delta g^{(k)}_{nm}=
  \frac{\Delta t}{2} \gs(t^k,\Xnewnm,\Vnewnm) + \frac{\Delta t}{2} \gs(t^{k+1},\xs_n,\vs_m) , 
  % ,\nonumber\\[3mm]
  % &\hskip1truecm k=0,1,\ldots,K-1,\,\,n=0,1,\ldots,N-1,\,\,m=0,1,\ldots,M-1.
\end{align}
which is an approximation of the average value of
  $\gs(t,X(t),Y(t))$ for $t\in[t^k,t^{k+1}]$ when moving along the characteristic lines  that solve  \eqref{eq:1D1V:char}.
The term $\Delta g^{(k)}_{nm}$ should be added to the right-hand side of
\eqref{eq:1D1V:eqcloumod}.

%
%Indeed, in order to handle the right-hand side, we suggest to use the
%trapezoidal rule to define
%%%
%\begin{align}\label{eq:1D1V:gtrapezi}
%  %&
%  {\color{red}
%    \Delta t\overline{\gs}^k_{nm} =
%  }
%  \frac{\Delta t}{2} \gs(t^k,\Xnewnm,\Vnewnm) + \frac{\Delta t}{2} \gs(t^{k+1},\xs_n,\vs_m)
%  \qquad t\in[t^{k},t^{k+1}],
%  % ,\nonumber\\[3mm] 
%  % &\hskip1truecm k=0,1,\ldots,K-1,\,\,n=0,1,\ldots,N-1,\,\,m=0,1,\ldots,M-1.
%\end{align}
%%%
%{\color{red}
%which is an approximation of the average value of
%  $\gs(t,X(t),Y(t))$ for $t\in[t^k,t^{k+1}]$.
%}
%%
%Term $\overline{\gs}^k_{nm}$ must be added to the right-hand side of
%\eqref{eq:1D1V:eqcloumod}.

Moreover, $g$ is also involved in the expression \eqref{eq:1D1V:eder}, that
must be rewritten as:
\begin{equation}\label{eq:1D1V:ederg}
  \displaystyle
  \frac{\partial\Es}{\partial t}(t,\xs) =
  \int_{\Omega_{v}} [\vs\fs(t,\xs,\vs) - \Gs(t,\xs,\vs)]\,\dv  , 
\end{equation}
where $\Gs$ is a primitive of the given function $\gs$ with respect to
the variable $\xs$, i.e.: $\partial\Gs/\partial\xs = \gs$.

In all the schemes proposed in this work, a CFL condition of stability
must be imposed on $\Delta t$. This is equivalent to the one shown in
\eqref{eq:1D1V:CFL}.  
We recall once again that all the space derivatives may be computed
with the help of the DFT, with a considerable time saving for $N$ and
$M$ large.
The methods proposed are the starting point to develop, within a
similar framework, more accurate schemes, in principal of any order.

%\input conservation.tex

%%%%%%%%%%%%%%%%%%%%%%%%%%%%%%%%%%%%%%%%%%%%%%%%%%%%%%%%%%
\section{Conservation Properties}
\label{sec:conservation:properties}
%%%%%%%%%%%%%%%%%%%%%%%%%%%%%%%%%%%%%%%%%%%%%%%%%%%%%%%%%%

%\textcolor{red}{}

The discrete counterpart of \eqref{eq:1D1V:massConserv} (i.e., number
of particles/mass/charge conservation) can be proven for the scheme
\eqref{eq:1D1V:favanza} - \eqref{eq:1D1V:favanzadef}.
This is the most basic quantity to be preserved, so that the check of
this relation is quite important from the physics viewpoint.
As in the previous sections  let  $t^{k}=k \Delta t= k\,{T}\slash{K}$,  $k=0,1,\ldots,K$. We start by defining:
\begin{equation}\label{eq:1D1V:massadiscreta}
  \Qs^{(k)}_{N,M} 
  = \frac{2\pi}{N}\frac{2\pi}{M}\sum_{n=0}^{N-1}\sum_{m=0}^{M-1}c^{(k)}_{nm} 
  = \int_\Omega f^{(k)}_{N,M}(x,v) dxdv \approx \int_\Omega f(t^k,x,v)dxdv,
\end{equation}
where we recalled the quadrature formula \eqref{eq:1D1V:fdq}.
The correspondence of the two integrals in
\eqref{eq:1D1V:massadiscreta} is true up to an error that is
spectrally accurate, due to the excellent properties of Gaussian
quadrature.  
By using \eqref{eq:1D1V:massadiscreta} for the timestep $k+1$ and
\eqref{eq:1D1V:favanza} we find that 
\begin{align}
  \Qs^{(k+1)}_{N,M}
  = \frac{2\pi}{N}\frac{2\pi}{M}\sum_{n=0}^{N-1}\sum_{m=0}^{M-1}c^{(k+1)}_{nm}
  = \frac{2\pi}{N}\frac{2\pi}{M}\sum_{n=0}^{N-1}\sum_{m=0}^{M-1}\Big(\cs_{nm}^{(k)} + \Delta t\,\Phi_{nm}^{(k)}\Big) 
  = \Qs^{(k)}_{N,M}+\Delta \Qs^{(k)}_{N,M},
\end{align}
where
\begin{align}
  %% ---
  &\Delta \Qs^{(k)}_{N,M} 
  = \Delta t \frac{2\pi}{N}\frac{2\pi}{M}\sum_{n=0}^{N-1}\sum_{m=0}^{M-1}\Phi^{(k)}_{nm}\nonumber\\
  &\hspace{5mm}= 
  -\Delta t \frac{2\pi}{M}\sum_{m=0}^{M-1} v_m\left[ \frac{2\pi}{N}\sum_{n=0}^{N-1}\frac{\partial f^{(k)}_{N,M}}{\partial x}(x_n, v_m) \right]
  +\Delta t \frac{2\pi}{N}\sum_{n=0}^{N-1} E^{(k)}_N(x_n)\left[ \frac{2\pi}{M}\sum_{m=0}^{M-1}
    \frac{\partial f^{(k)}_{N,M}}{\partial v}(x_n, v_m)\right]\nonumber\\
  %% ---
  &\hspace{5mm}=
  -\Delta t \frac{2\pi}{M}\sum_{m=0}^{M-1} v_m\left[ \int_{\Omega_x}\frac{\partial f^{(k)}_{N,M}}{\partial x}(x, v_m)dx\right]
  +\Delta t \frac{2\pi}{N}\sum_{n=0}^{N-1} E^{(k)}_N(x_n)\left[ \int_{\Omega_v}\frac{\partial f^{(k)}_{N,M}}{\partial v}(x_n, v)dv\right]\nonumber\\
  %% --
  &\hspace{5mm}= 0.
  \label{eq:1D1V:massadiscretacons}
\end{align}
Here, we may note that the two integrals are zero as a consequence of
the boundary conditions (periodic or homogeneous Dirichlet).
This shows that the quantity in \eqref{eq:1D1V:massadiscreta} does not
change from $k$ to $k+1$.
The same property holds for the schemes \eqref{eq:1D1V:favanzabdf2co}
and \eqref{eq:1D1V:favanzabdf3co}.
The proof follows after recognizing that, for $g=0$, the sum of the
coefficients on the right-hand side is equal to 1.
In fact, for \eqref{eq:1D1V:favanzabdf2co} one has: $(4/3)-(1/3)=1$,
and for \eqref{eq:1D1V:favanzabdf3co} one has:
$(18/11)-(9/11)+(2/11)=1$.

Concerning the scheme \eqref{eq:1D1V:eqcloumod}, the conservation of
$\Qs^{(k)}_{N,M}$ is also recovered, but one has to be a bit more
careful in the analysis.
As a matter of fact, there are terms containing second derivatives in
$x$ and $v$, multiplying $(\Delta t)^2$.
With the same arguments followed to recover
\eqref{eq:1D1V:massadiscretacons}, these parts can be transformed in
integrals by Gaussian quadrature.
Their contribution is zero if appropriate boundary conditions are
assumed.
For instance, in the periodic case, all the derivatives are matching
across the point $2\pi$ (see \eqref{eq:1D1V:PeriodicBC}), therefore we
have perfect mass conservation (i.e., the discrete version of it).
With homogeneous Dirichlet boundary conditions, we have no elements to
argue that the integral contribution of the second derivatives
must be zero (because the first derivatives in 0 and $2\pi$ are not
necessarily equal), so that mass conservation is achieved up to an
error proportional to $(\Delta t)^2$.
Nevertheless, if an exponential decay of $f$ is assumed near the
boundary (as it is commonly accepted concerning the variable $v$), the
first derivatives will also decay in the same way, and the integral
contribution of the second derivatives can be again neglected.
In the experiments of the next sections, we assume full periodicity in
the direction $x$, while, in the variable $v$, we will work with
functions exhibiting an exponential decay.
Therefore, up to possible negligible effects developing at the
boundaries, mass conservation is ensured.

Similar considerations can be made regarding the conservation in time
of other quantities, such as the momentum $\int_\Omega v f(t,x,v) dx
dv$, which in the discrete case is defined at time $t^{k}$,  $k=0,1,\ldots,K$,  in the
following way:
\begin{equation}\label{eq:1D1V:momentodiscreto}
  \Ps^{(k)}_{N,M} 
  = \frac{2\pi}{N}\frac{2\pi}{M}\sum_{n=0}^{N-1}\sum_{m=0}^{M-1}v_{m}c^{(k)}_{nm} 
  \approx \int_\Omega v f^{(k)}_{N,M}(x,v) dxdv.
\end{equation}
Here, it has to be noticed that the function $v$ is not a
trigonometric polynomial, so that it is not possible to use the
quadrature formula \eqref{eq:1D1V:fdq} in a straightforward way.
On the other hand, $v$ can be substituted by its projection (in the
$L^2(\Omega)$ norm) on the finite dimensional space ${\bf Y}_{N,M}$
(see \eqref{eq:1D1V:fdq}) up to an error that decays spectrally.
This procedure may however generate a Gibb's phenomenon across the
points of $\Omega$ with $v=2\pi$, where $vf$ is discontinuous.
The trouble can be fixed by supposing that the function $f$ decays as
an exponential (with respect to the variable $v$) near the boundary.
In the end, with assumptions that may be considered standard in
applications, the conservation of momentum can be achieved up to
negligible errors.

A discussion can also be made regarding the discrete version of
\eqref{eq:1D1V:totEnergy} at time $t^{k}$,  $k=0,1,\ldots,K$, i.e.:
\begin{equation}\label{eq:1D1V:energiadiscreta}
  {\cal E}(t^k) \approx {\cal E}^{(k)}_{N,M}
  = \frac12\left( 
    \frac{2\pi}{N}\frac{2\pi}{M}\sum_{n=0}^{N-1}\sum_{m=0}^{M-1}v^2_mc^{(k)}_{nm} +
    \frac{2\pi}{N}\sum_{n=0}^{N-1}\Big[ E^{(k)}_N (x_n)\Big]^2
  \right).
\end{equation}
The theoretical analysis now becomes more involved, since the above
quantity is quadratic.
We expect however that conservation at each step is achieved up to an
error that is at most proportional to $(\Delta t)^S$, where $S$ is the
order of the scheme used.
Exact conservation cannot be expected in this case, due to the fact
that all the time-advancing schemes we consider in this paper are of
explicit type.
Energy conservation is usually a prerogative of implicit schemes (see,
e.g., the Crank-Nicholson method).

Finally, we spend a few words on the treatment of the term in
\eqref{eq:1D1V:totEnergy}.
As already observed above, the function $v^2$ is not a trigonometric
polynomial, therefore in the theoretical analysis we need to replace
it with a suitable projection.
In order to avoid possible Gibb's phenomena at the boundary, we should
rely on the fast decay of the function $f$. On the other hand, these
considerations must also be used in the continuous case, because they
are necessary to give a meaning to the integral $\int_\Omega v^2
f(t,x,v) dx dv$.
In addition, we also point out that there is no proof that the
quantity defined in \eqref{eq:1D1V:totEnergy} is actually a norm,
since it is not guaranteed that, if the discrete quantity  $f_{N,M}\simeq f$ is
positive at time $t=0$, it will remain positive in the subsequent
times.
Anyway, this trouble is frequently present within the framework of any
other type of approximations, unless it is built on purpose to be
sign-preserving (a rather difficult property to achieve).
The possible negativity of $f_{N,M}\simeq f$  has not in general significant
relevance in practical experiments, but makes the theoretical aspects
far more involved.  
For the reasons mentioned above, we omit the details of the study of
energy conservation, because they are rather complicate and out of the
scopes of this paper.
Numerical confirmations of the above statements will be
given in the coming sections.

\raggedbottom
%%%%%%%%%%%%%%%%%%%%%%%%%%%%%%%%%%%%%%%%%%%%%%%%%%%%%%%%%%
\section{Numerical experiments}
\label{sec:numerical:experiments}
%%%%%%%%%%%%%%%%%%%%%%%%%%%%%%%%%%%%%%%%%%%%%%%%%%%%%%%%%%

%% filamentation filtration {Klimas-Farrell:1994}

%%%%%%%%%%%%%%%%%%%%%%%%%%%%%%%%%%%%%%%%%%%%%%%%%%%%%%%%%%
\subsection{Manufactured solution benchmark}
\label{subsec:convergence}
%%%%%%%%%%%%%%%%%%%%%%%%%%%%%%%%%%%%%%%%%%%%%%%%%%%%%%%%%%

The aim of this first test is to assess the convergence rate of our
numerical schemes.
We consider the non-homogeneous Vlasov-Poisson problem
\eqref{eq:1D1V:vlasovno}, \eqref{eq:1D1V:Vci},
\eqref{eq:1D1V:EMaxwell}, \eqref{eq:1D1V:charge-density}, where we set
$\Omega_{x}=[0,2\pi]$, $\Omega_{v}=[-\pi,\pi]$, $T=1$.
%%^
The right-hand side $\gs$ in \eqref{eq:1D1V:vlasovno} is such that
% $t\in[0,T],$ $\xs\in\Omega_{x}$, $\vs\in\Omega_{v}$,
the solution fields $\fs$ and $\Es$ are given by:
\begin{align} 
  &
  \fs(t,\xs,\vs) = \frac{2}{\sqrt{\pi}}  \left[ 1-\cos(2\xs-2\pi t ) \right] \exp (-4 \vs^{2}), 
  %%\quad t\in [0,T],\,\,\xs\in\Omega_{x},\,\,\vs\in\Omega_{v},
  \label{eq:1D1V:exeactf}\\
  &
  \Es(t,\xs) = \frac12 \sin(2\xs-2\pi t).
  %% , \quad t\in[0,T],\,\,\xs\in\Omega_{x}.
  \label{eq:1D1V:exeactE}
\end{align}
We note that both $\fs$ and $\Es$ are $2\pi$-periodic in the variable
$\xs$.
Instead, $\fs$ is not periodic in the variable $\vs$ but we can
effectively approximate it by periodic functions since the Gaussian
function $\exp{(-4\vs^{2})}$ is practically zero at the velocity
boundaries $\vs=\pm\pi$.

%\input tablefE.tex

%!TEX root=paperVlasov.tex
\newcommand{\TABROW}[9]{ #1 &\,\, $#2$ & $ \quad#3$ &\,\, $#4$ &\quad $#5$ &\,\, $#6$ & $\quad #7$ &\,\, $#8$ &\quad $#9$ } %
\newcommand{\HEADER}[9]{ #1 &      #2  &\quad  #3  &\,\,  #4  &\quad  #5  &\,\,  #6  &\quad   #7  &\,\,  #8  &\quad  #9  } %
\renewcommand*{\arraystretch}{0.98}
\begin{table}
  \centering
  \begin{tabular}{lllllllll}
    \hline
    \HEADER{$\Delta t$}{ One-step first-  }{ Rate }{ Second-order  }{ Rate }{ Third-order }{ Rate }{ One-step second-}{ Rate }\\
    \HEADER{          }{ order scheme    }{      }{ BDF method    }{      }{ BDF method  }{      }{ order scheme   }{      }\\
    \HEADER{ }{\eqref{eq:1D1V:favanza}}{ }{\eqref{eq:1D1V:favanzabdf2co}}{ }{\eqref{eq:1D1V:favanzabdf3co}}{ }{\eqref{eq:1D1V:eqcloumod},\eqref{eq:1D1V:gtrapezi},\eqref{eq:1D1V:ederg}}{}\\                 
    \hline
    \TABROW{  0.04      }{ 8.86\,\,10^{-2} }{      }{ 2.78\,\,10^{-2} }{      }{ 4.32\,\,10^{-3} }{      }{ 4.03\,\,10^{-3} }{      } \\
    \TABROW{  0.02      }{ 4.24\,\,10^{-2} }{ 1.06 }{ 6.75,\,10^{-3} }{ 2.04 }{ 5.66\,\,10^{-4} }{ 2.93 }{ 1.01\,\,10^{-3} }{ 2.00 }  \\
    \TABROW{  0.01      }{ 2.07\,\,10^{-2} }{ 1.03 }{ 1.65\,\,10^{-3} }{ 2.03 }{ 7.27\,\,10^{-5} }{ 2.96 }{ 2.51\,\,10^{-4} }{ 2.01 } \\ 
    \TABROW{  0.005     }{ 1.02\,\,10^{-2} }{ 1.02 }{ 4.09\,\,10^{-4} }{ 2.02 }{ 9.25\,\,10^{-6} }{ 2.97 }{ 6.28\,\,10^{-5} }{ 2.00 } \\   
    \TABROW{  0.0025    }{ 5.08\,\,10^{-3} }{ 1.01 }{ 1.02\,\,10^{-4} }{ 2.01 }{ 1.17\,\,10^{-6} }{ 2.99 }{ 1.57\,\,10^{-5} }{ 2.00 } \\
    \TABROW{  0.001325  }{ 2.53\,\,10^{-3} }{ 1.00 }{ 2.53\,\,10^{-5} }{ 2.00 }{ 1.47\,\,10^{-7} }{ 3.00 }{ 3.93\,\,10^{-6} }{ 2.00 } \\
    \hline
  \end{tabular}
  \vspace{.4cm}
  \caption{Relative errors between the exact and the numerical distribution functions in the $L^2(\Omega)$ norm, obtained with different time discretization schemes. 
    The corresponding convergence rate is reported aside.}
  \label{tab1}    
  \vskip1truecm
\end{table}

\begin{table}
  \centering
  \begin{tabular}{lllllllll}
    \hline
    \HEADER{$\Delta t$}{ One-step first-  }{ Rate }{ Second-order  }{ Rate }{ Third-order }{ Rate }{ One-step second-}{ Rate }\\
    \HEADER{          }{ order scheme    }{      }{ BDF method    }{      }{ BDF method  }{      }{ order scheme   }{      }\\
    \HEADER{}{\eqref{eq:1D1V:favanza}}{ }{\eqref{eq:1D1V:favanzabdf2co}}{ }{\eqref{eq:1D1V:favanzabdf3co}}{ }{\eqref{eq:1D1V:eqcloumod},\eqref{eq:1D1V:gtrapezi},\eqref{eq:1D1V:ederg}}{}\\                 
    \hline
    \TABROW{  0.04      }    { 8.18\,\,10^{-2} }{      }   { 3.21\,\,10^{-2} }{      }   { 2.78\,\,10^{-3} }   {      }{ 3.56\,\,10^{-3} }{      }  \\
    \TABROW{  0.02      }    { 4.16\,\,10^{-2} }{ 0.98 }{ 7.95\,\,10^{-3} }{ 2.01 }{ 3.88\,\,10^{-4} }{ 2.84 }{ 8.86\,\,10^{-4} }{ 2.01 }  \\
    \TABROW{  0.01      }    { 2.10\,\,10^{-2} }{ 0.99 }{ 1.97\,\,10^{-3} }{ 2.02 }{ 5.20\,\,10^{-5} }{ 2.90 }{ 2.21\,\,10^{-4} }{ 2.00 } \\ 
    \TABROW{  0.005     }   { 1.05\,\,10^{-2} }{ 0.99 }{ 4.88\,\,10^{-4} }{ 2.01 }{ 6.75\,\,10^{-6} }{ 2.95 }{ 5.52\,\,10^{-5} }{ 2.00 } \\   
    \TABROW{  0.0025    }  { 5.27\,\,10^{-3} }{ 1.00 }{ 1.22\,\,10^{-4} }{ 2.01 }{ 8.59\,\,10^{-7} }{ 2.98 }{ 1.38\,\,10^{-5} }{ 2.00 } \\
    \TABROW{  0.001325  }{ 2.64\,\,10^{-3} }{ 1.00 }{ 3.03\,\,10^{-5} }{ 2.00 }{ 1.08\,\,10^{-7} }{ 3.00 }{ 3.45\,\,10^{-6} }{ 2.00 } \\
    \hline
  \end{tabular}
  \vspace{.4cm}
  \caption{Relative errors between the exact and the numerical electric field in the $L^2(\Omega_x)$ norm, obtained with different time discretization schemes.
    The corresponding convergence rate is reported aside.}
  \label{tab2}    
  \vskip1truecm
\end{table}

Table~\ref{tab1} shows the relative errors and the convergence rates
at the final time $T=1$ between the exact solution
\eqref{eq:1D1V:exeactf} and the numerical solution obtained with the
different schemes proposed in Sections~\ref{sec:time:discretization}
and~\ref{sec:more:advanced:time:discretizations}.
%%
%% \eqref{eq:1D1V:exeactE} risultati per E
%%
These calculations are performed with a fixed number of spectral modes
($N=M=2^{5}$).
We decreased the timestep by halving the initial value $\Delta t=0.04$
at each refinement.
The first column reports the timestep.
The other columns report the relative errors in the $L^2(\Omega)$ norm
and the corresponding convergence rates, when using the various
schemes.
The results of Table~\ref{tab2} pertain to the error of the electric
field.
They confirm the convergence rates shown in Table~\ref{tab1}.
In all these tests we assumed that the time discretization error
dominates the approximation error of the phase space.
Indeed, for the relatively small number of degrees of freedom
$N=M=2^5$, the resolution in $x$ and $v$ is excellent, due to the
convergence properties of the spectral approximations.

%%%%%%%%%%%%%%%%%%%%%%%%%%%%%%%%%%%%%%%%%%%%%%%%%%%%%%%%%%
\subsection{Two-stream instability}\label{sec32}
\begin{figure}
  \centerline{
    \includegraphics[height=6.25cm]{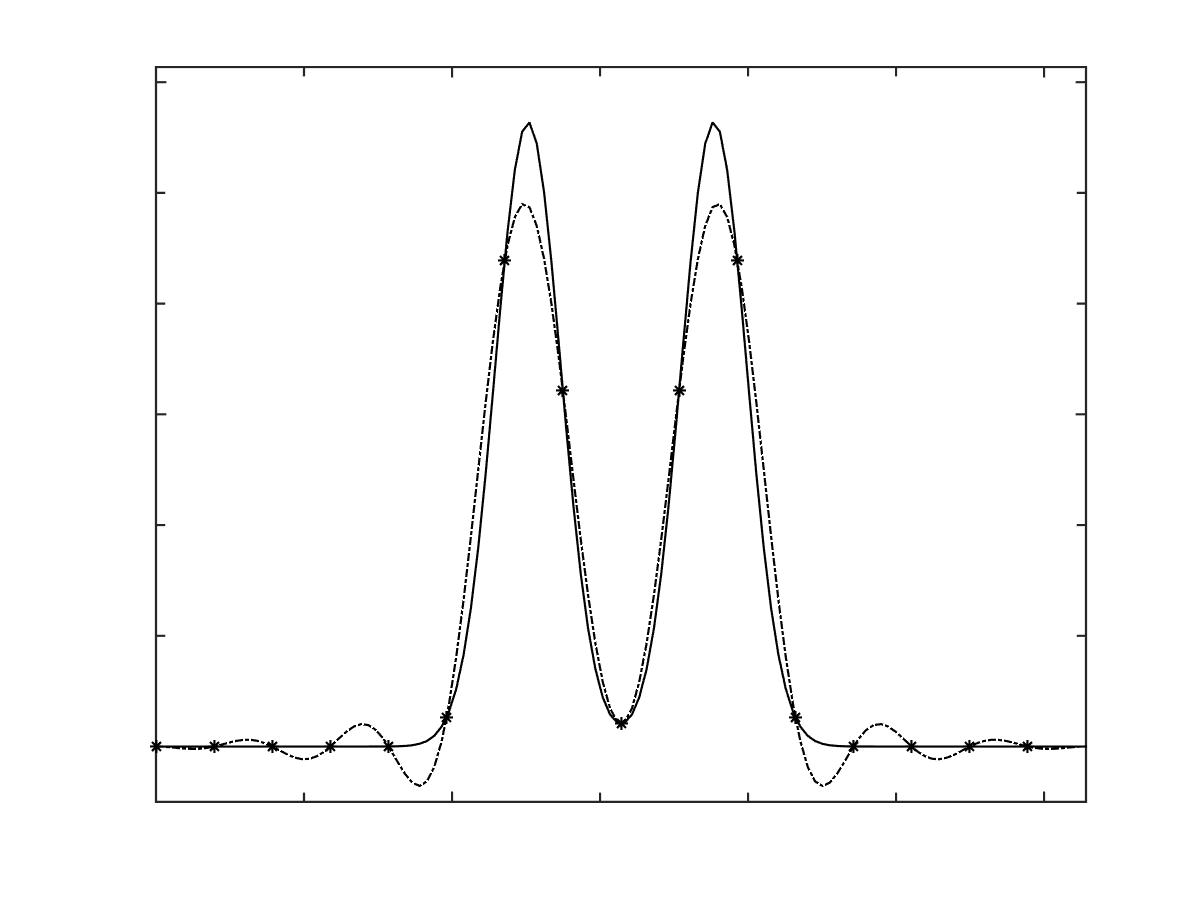}%{./gaussN2^4.jpg}      
    \includegraphics[height=6.25cm]{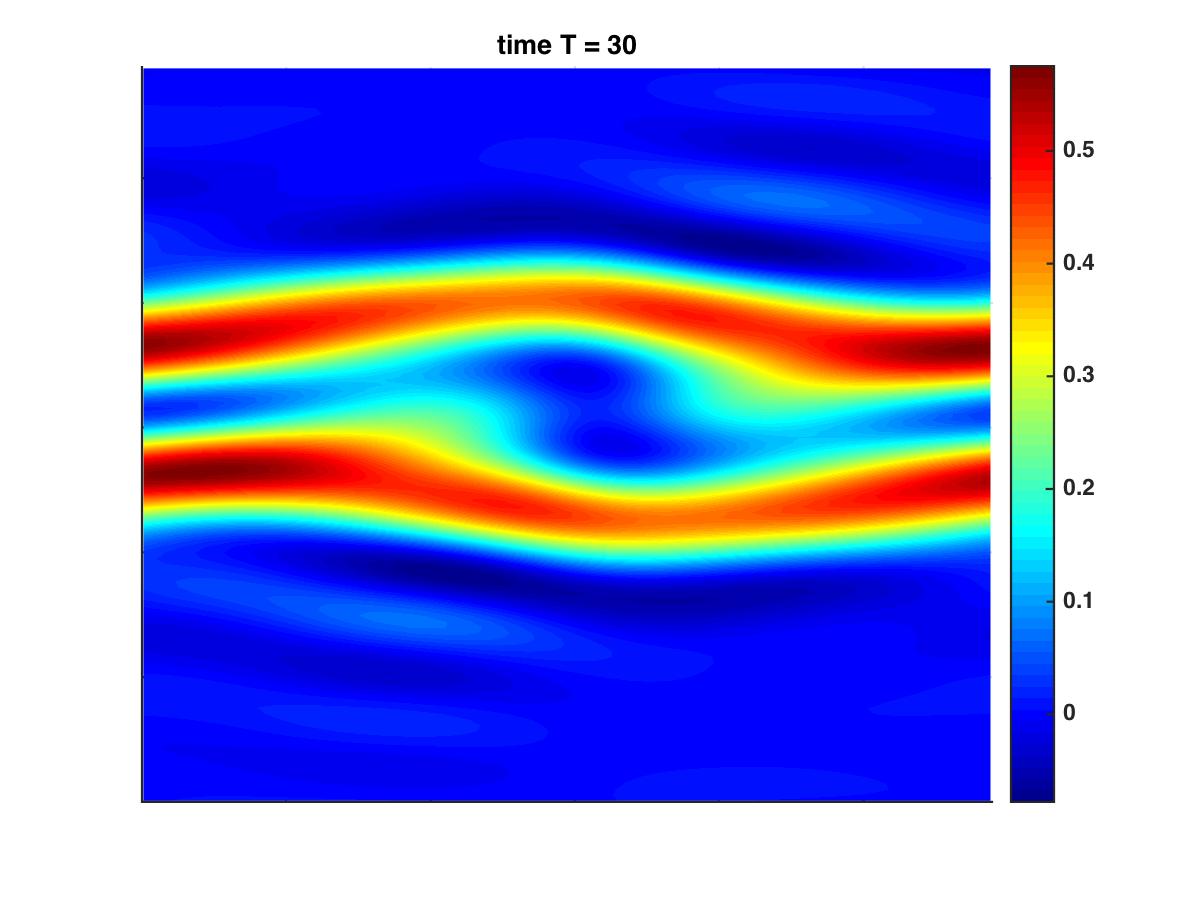}%{./fig1bNM2^4.jpg}
  }
  \centerline{
    \includegraphics[height=6.25cm]{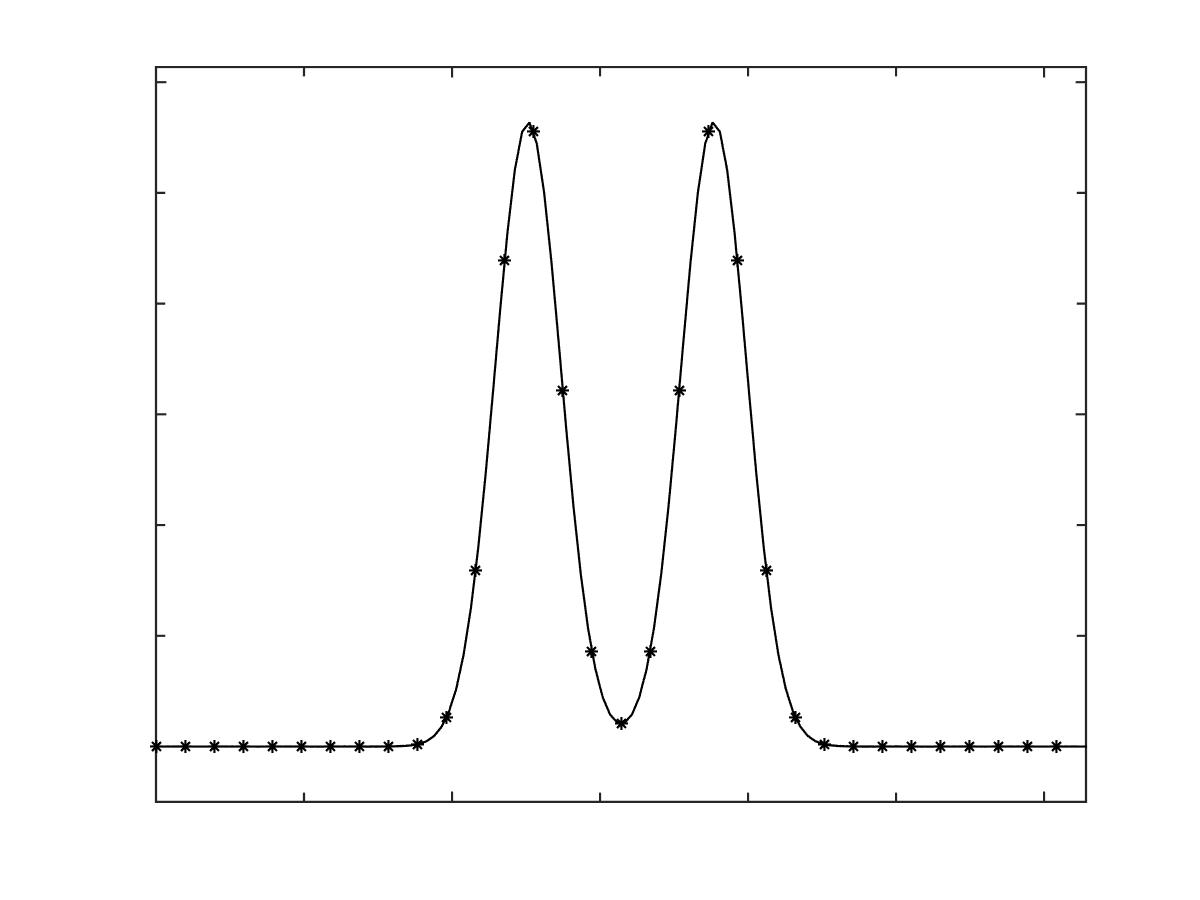}%{./gaussN2^5.jpg} 
    \includegraphics[height=6.25cm]{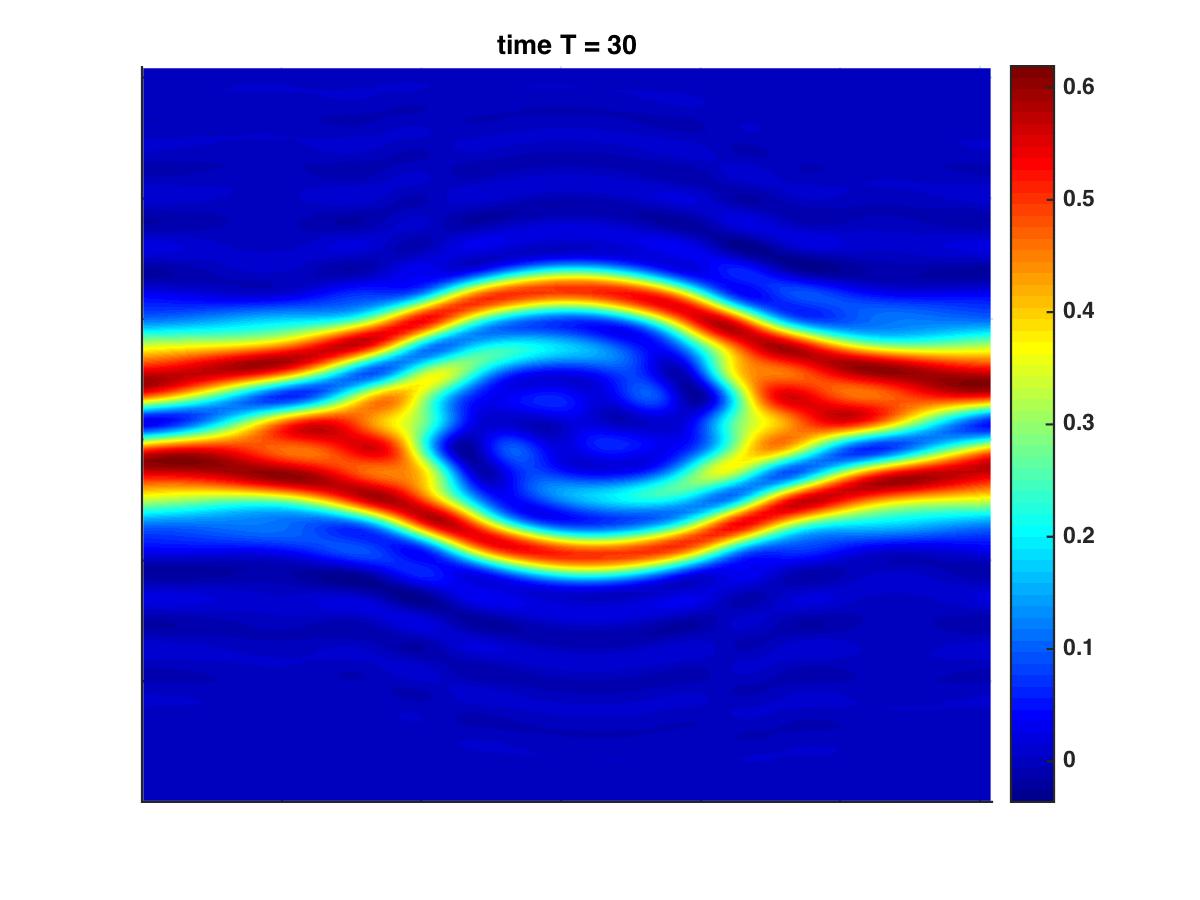}%{./fig1bNM2^5.jpg} 
  }
  \centerline{
    \includegraphics[height=6.25cm]{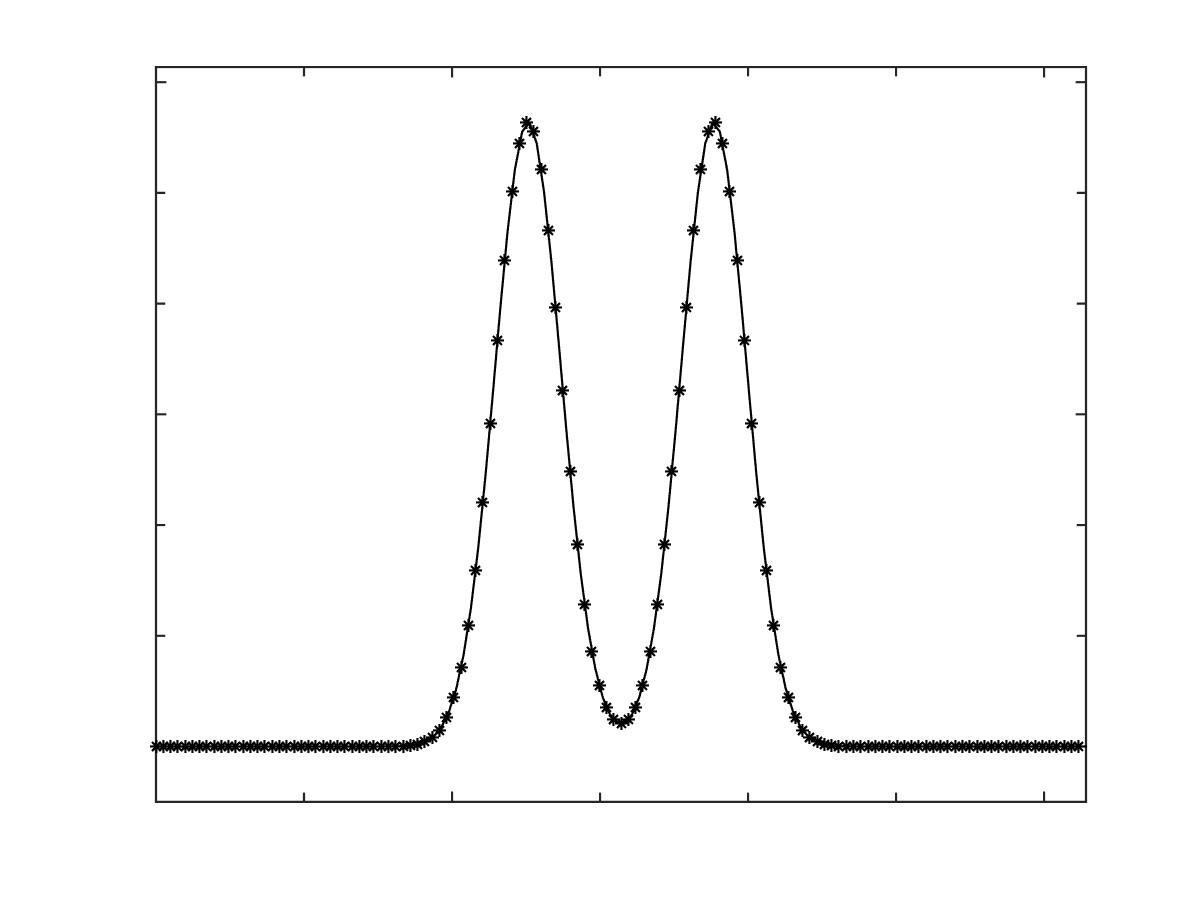}%{./gaussN2^7.jpg} 
    \includegraphics[height=6.25cm]{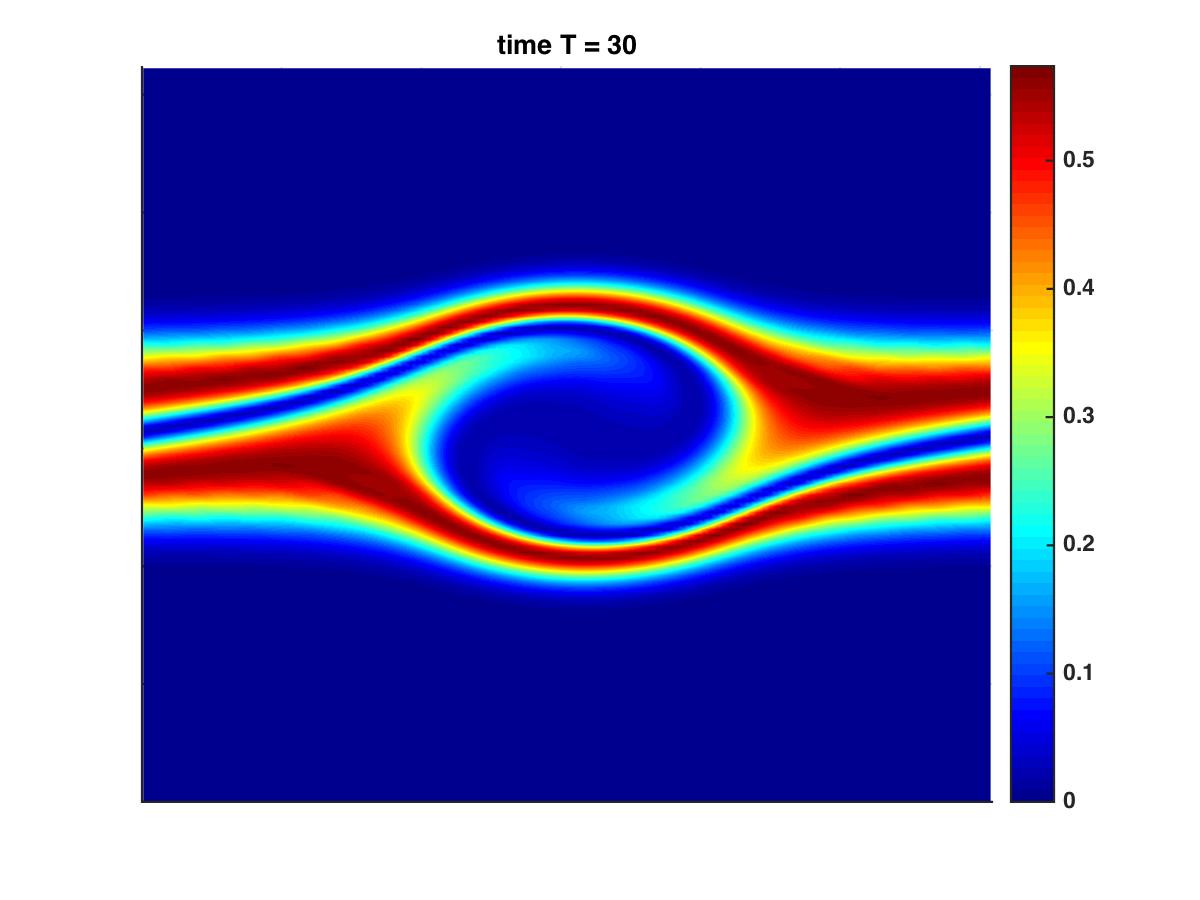}%{./fig1bN2^5M2^7.jpg} 
  }
  \caption{ Two-stream instability test: interpolation with respect to
    $v$ of the initial solution \eqref{eq:1D1V:twostreamf} at time
    $t=0$ and $x=0$ (left plots); approximated distribution function
    in the domain $\Omega=\Omega_x\times\Omega_v$ at time $T=30$
    (right plots).
    The second-order one-step time-marching scheme is implemented with
    $\Delta t=10^{-2}$ and $N=M=2^{4}$ (top), $N=M=2^{5}$ (center) and
    $N=2^{5}$, $M=2^{7}$ (bottom).}
  \label{fig1two-streamApprox}
\end{figure}

\begin{figure}
  \centerline{
    \includegraphics[height=6.25cm]{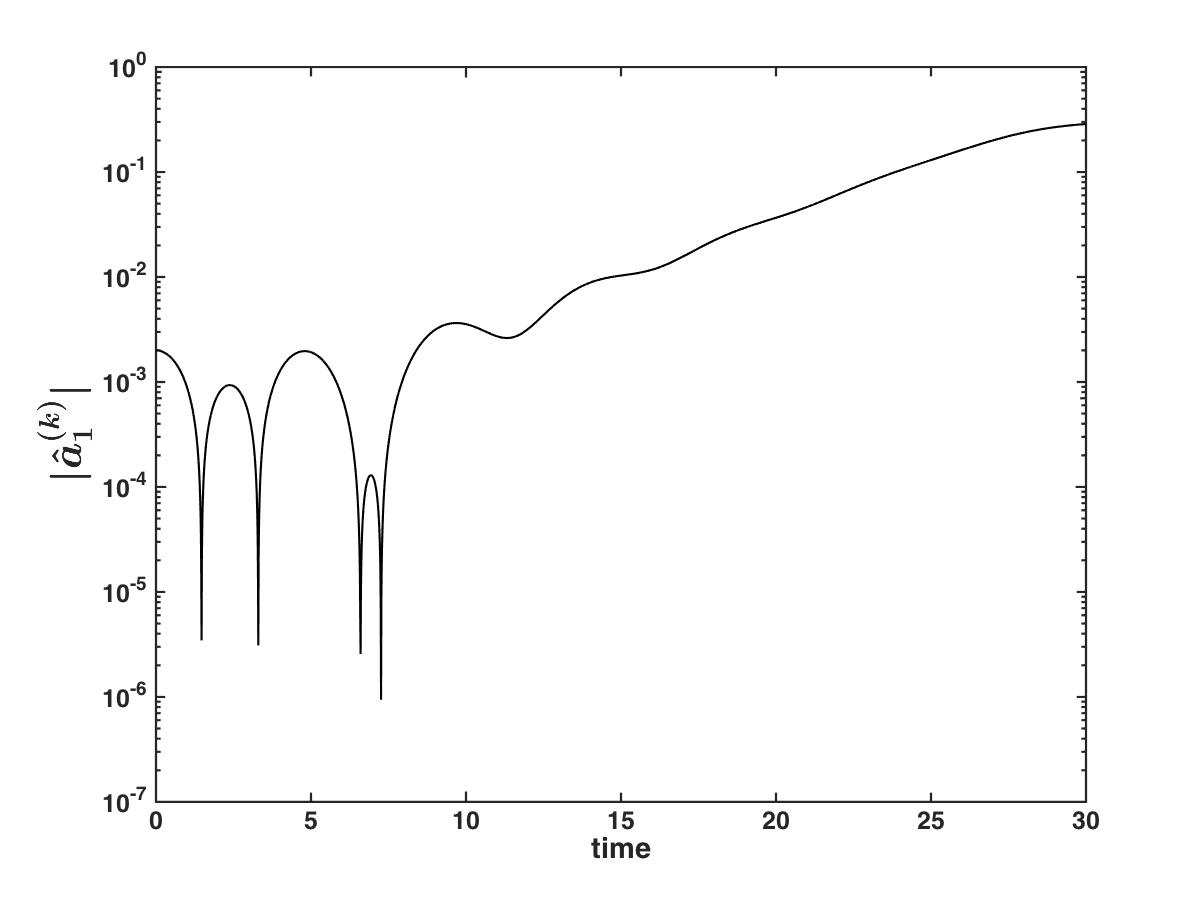}%{./BubboneImodoEN2^5M2^7T30.jpg}
    \includegraphics[height=6.25cm]{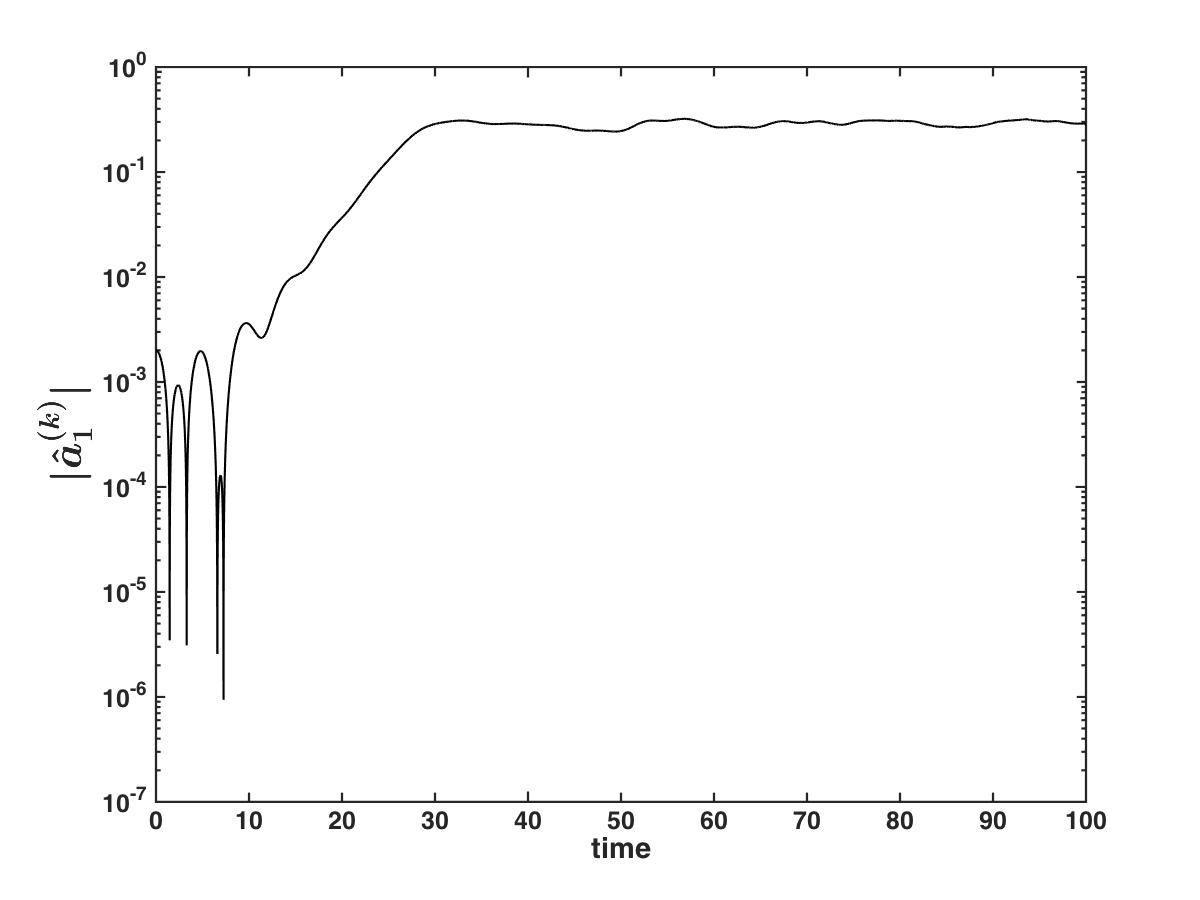}%{./BubboneImodoEN2^5M2^7T100.jpg} 
  }
  \caption{Two-stream instability test: the first Fourier mode versus
    time, $|\hat a_{1}^{(k)}|$,  of the electric field $|E^{(k)}_N|$ in
    \eqref{eq:1D1V:FourierSeriesEnodi}, when using the second-order
    one-step scheme, for $N=2^{5}$, $M=2^{7}$, $\Delta t=10^{-2}$ and
    $T=100$. The plot on the left is an enlargement of the graph in
    the time interval $[0,30]$.}
  \label{fig2two-streamImodoE}
\end{figure}

\begin{figure}
  \centerline{
    \includegraphics[height=6.25cm]{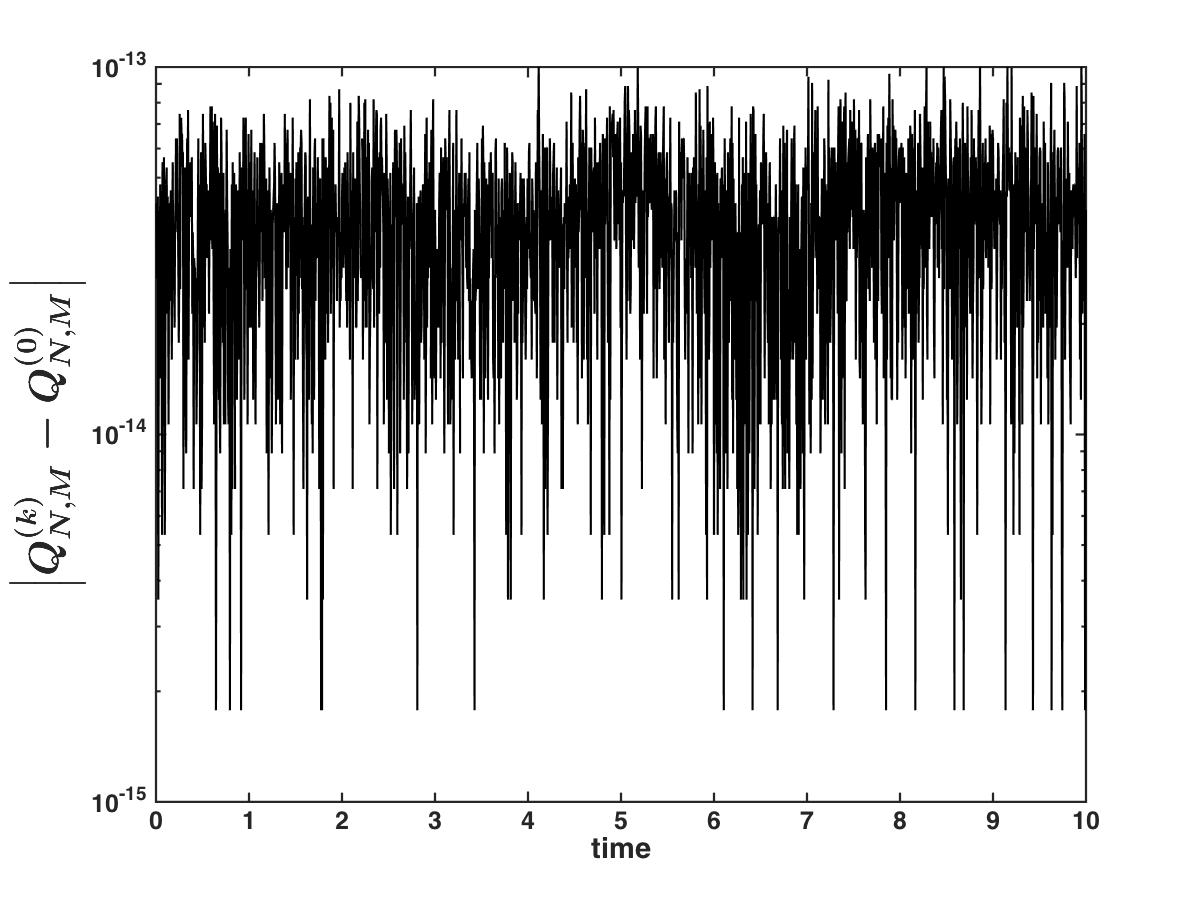}%{./BDF1BubboneMassaT10dt0005.jpg} 
    \includegraphics[height=6.25cm]{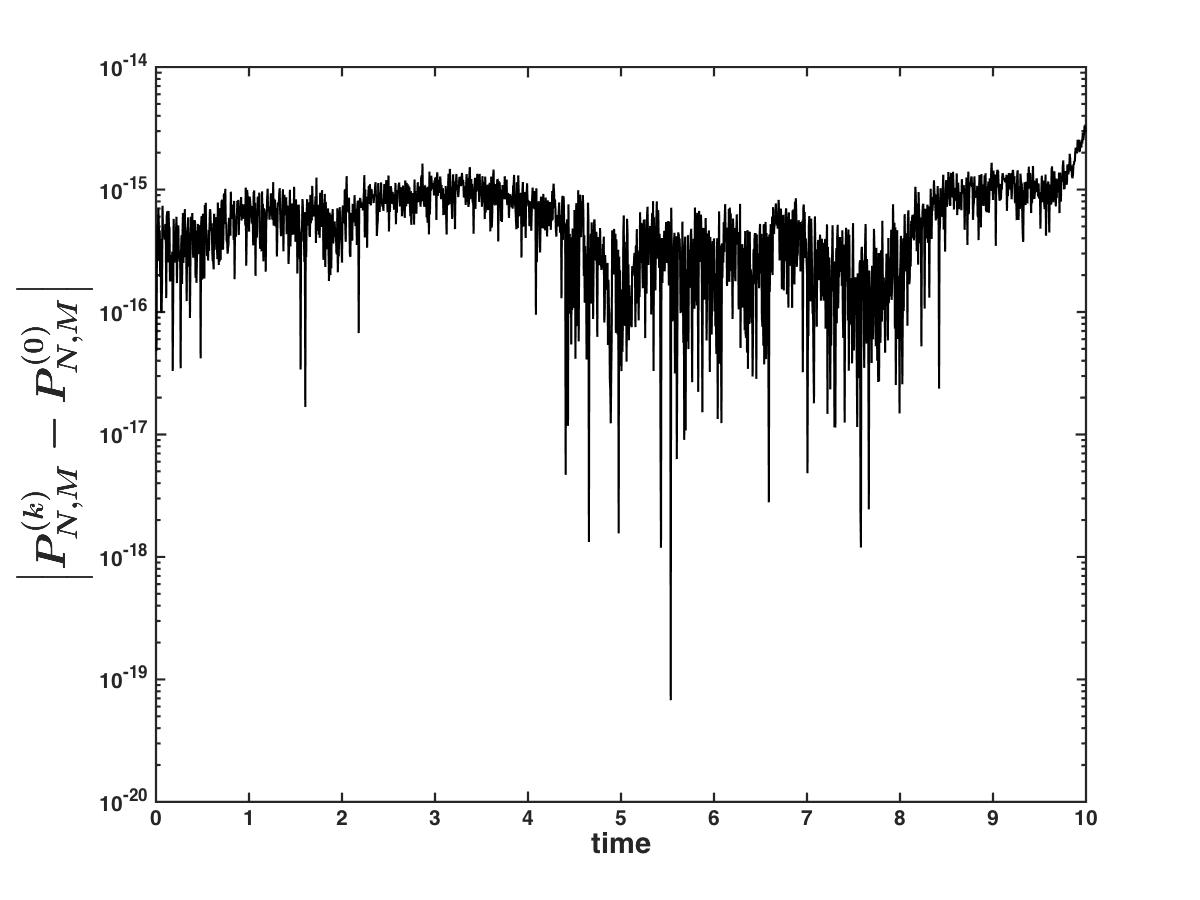}%{./BDF1BubboneMomentoT10dt0005.jpg} 
  }
  \centerline{
    \includegraphics[height=6.25cm]{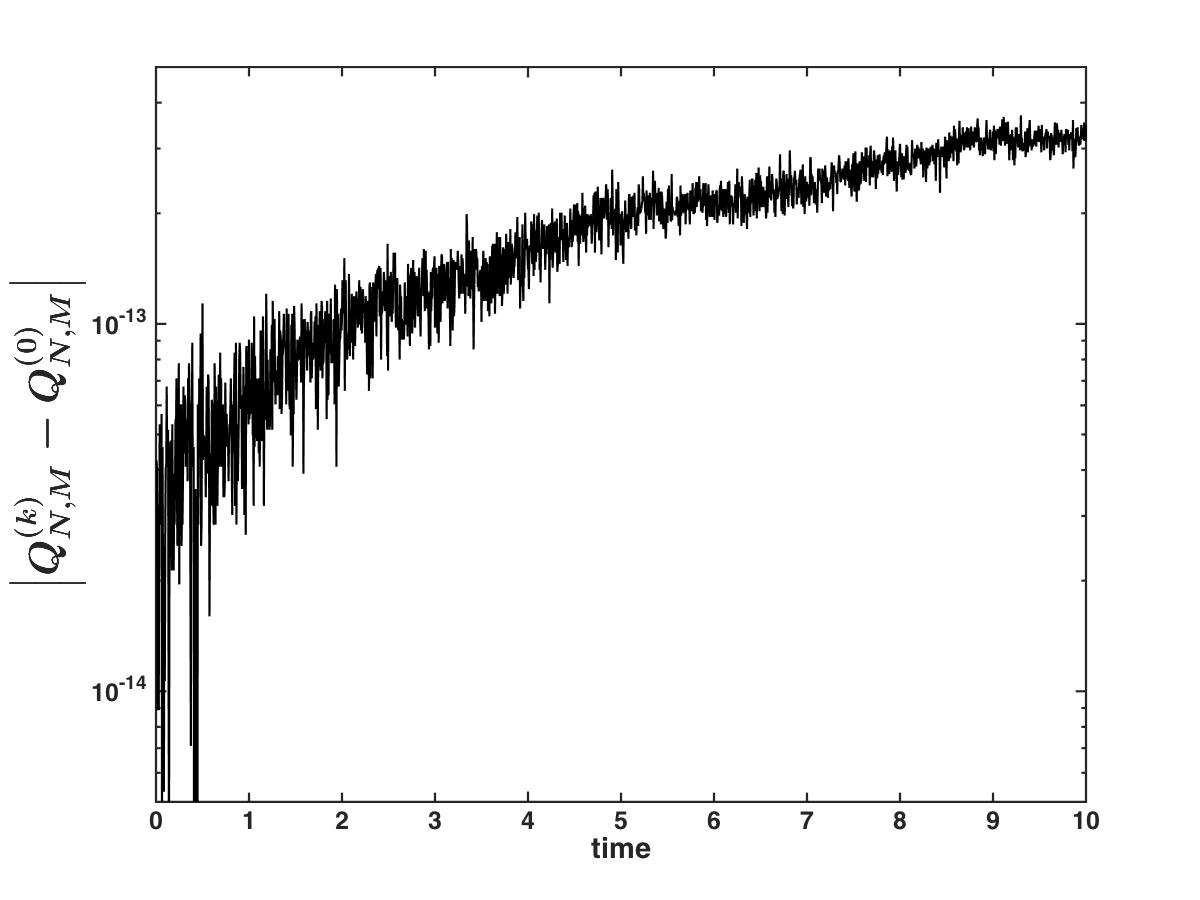}%{./BDF2BubboneMassaT10dt0005.jpg} 
    \includegraphics[height=6.25cm]{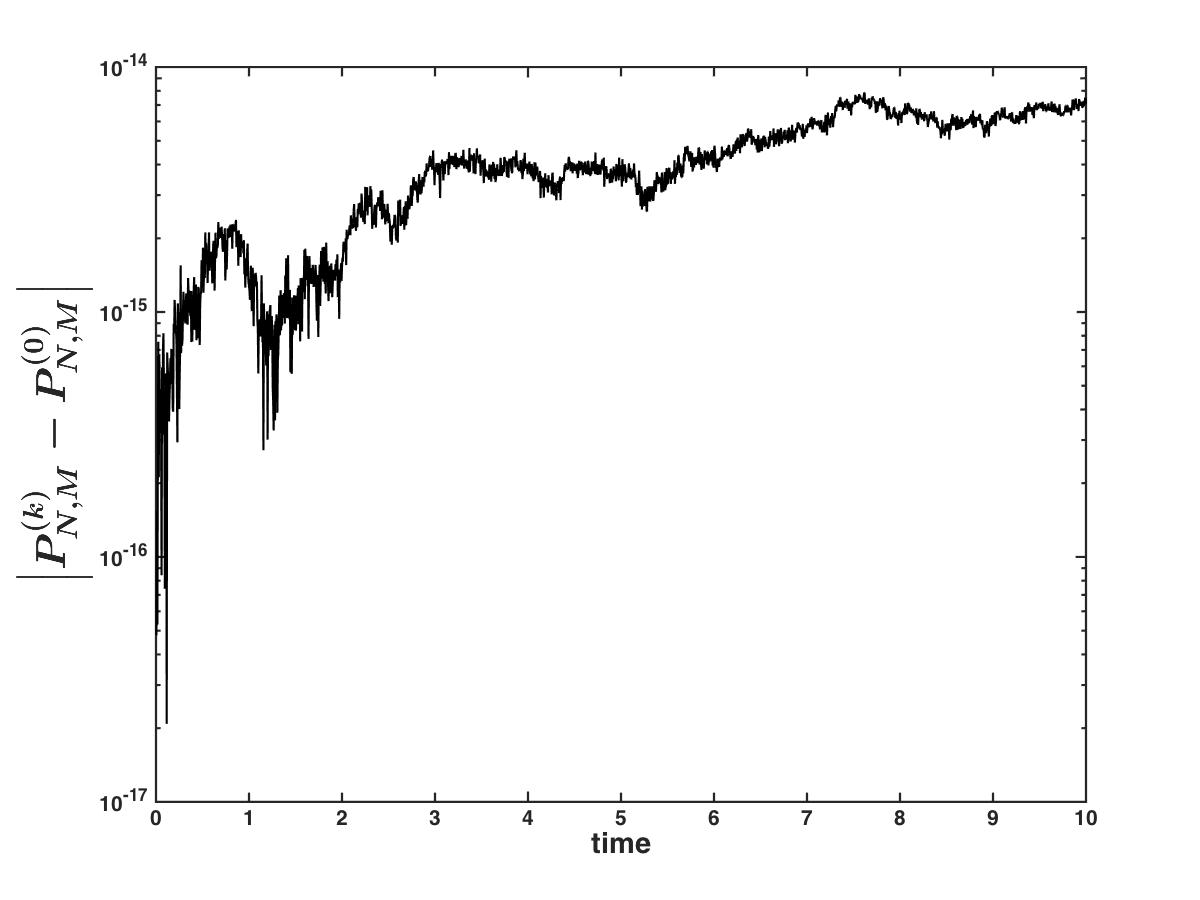}%{./BDF2BubboneMomentoT10dt0005.jpg} 
  }
  \centerline{
    \includegraphics[height=6.25cm]{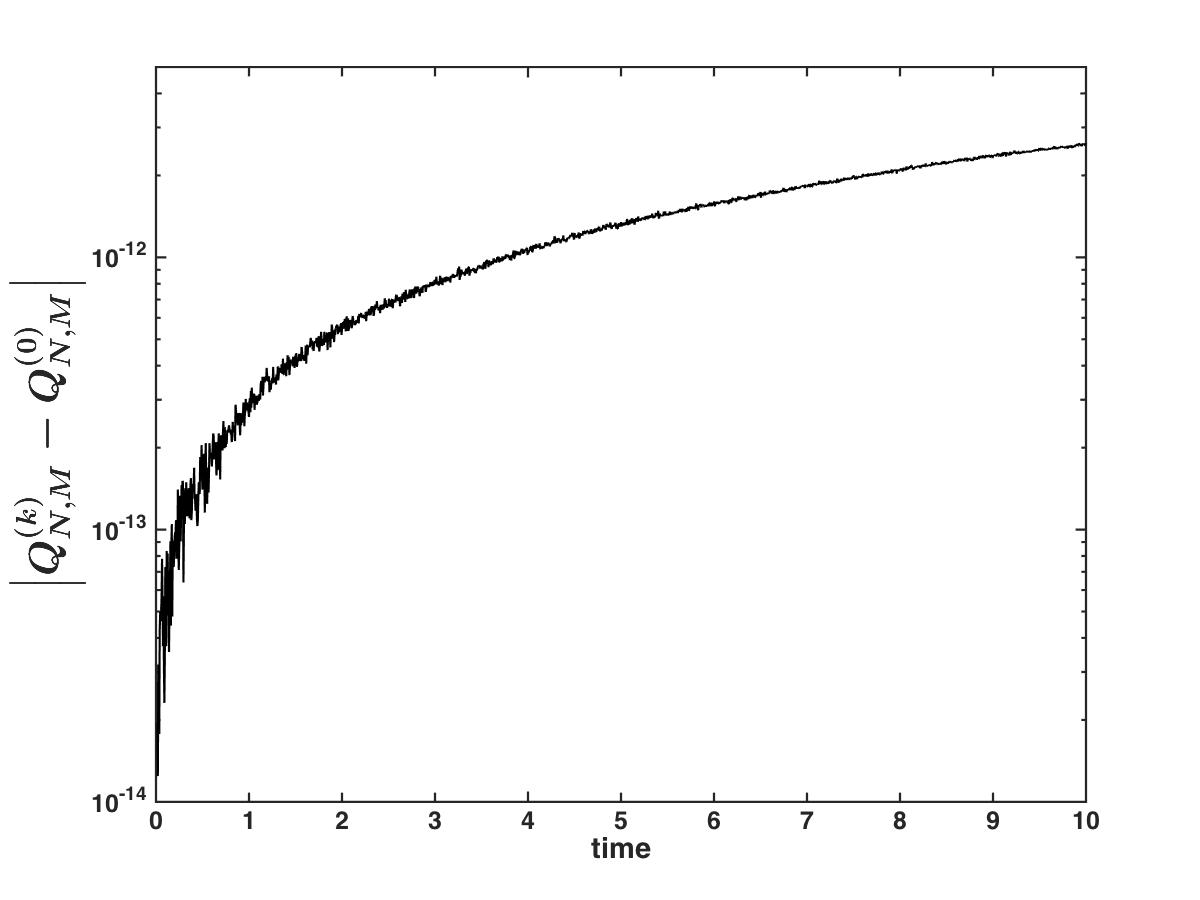}%{./BDF3BubboneMassaT10dt0005.jpg} 
    \includegraphics[height=6.25cm]{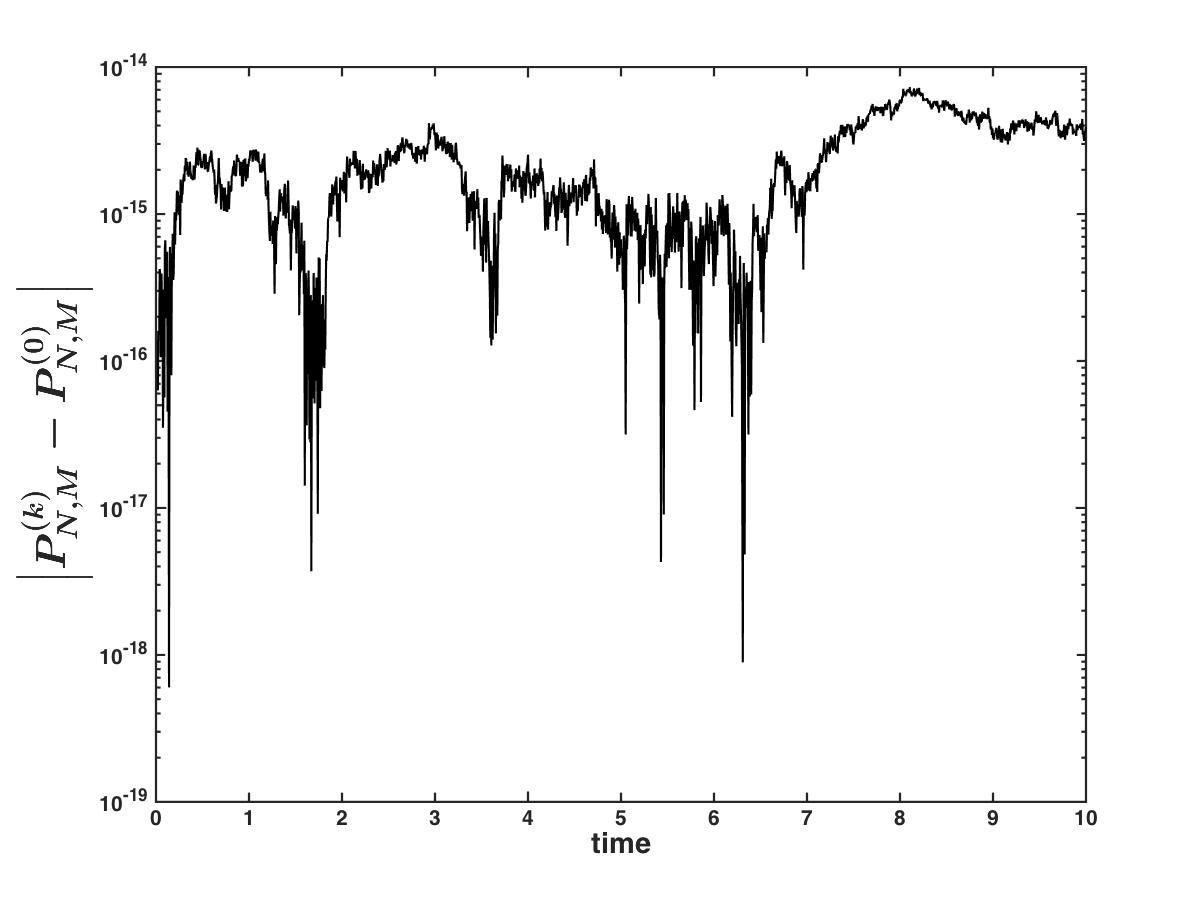}%{./BDF3BubboneMomentoT10dt0005.jpg} 
  }
  \caption{Two-stream instability test: conservation of the number of
    particles (left) and momentum (right), when applying the
    first-order one-step scheme (top), the second-order
    BDF scheme (center) and the third-order BDF
    scheme (bottom). The plots show the variation with
    respect to the initial value. All calculations are
    carried out by choosing $N=2^{5}$, $M=2^{7}$, $T=10$ and $\Delta
    t=5\cdot10^{-3}$.  
  }
  \label{fig3two-streamMassaMomento}
\end{figure}

\begin{figure}
  \centerline{\includegraphics[height=7.5cm]{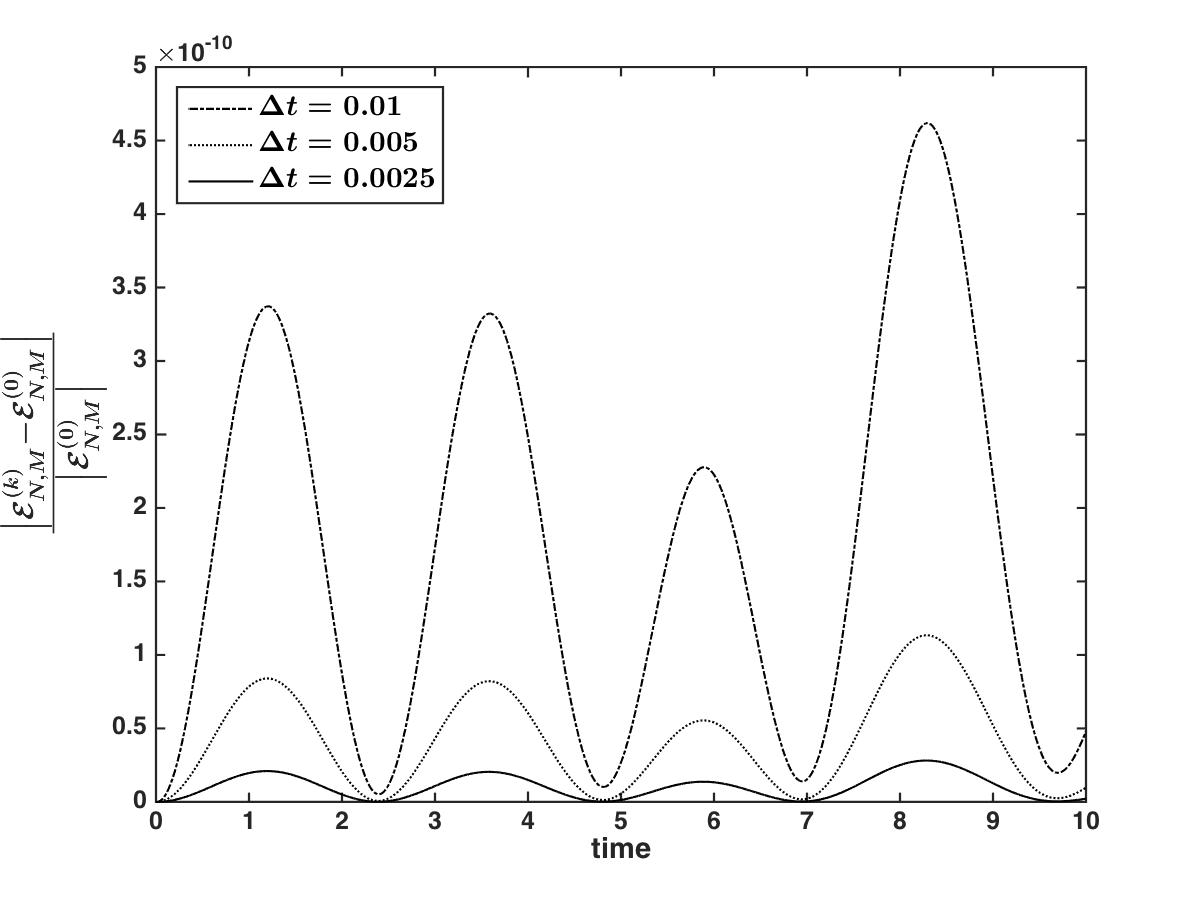}} %{./BDF2confrontoEnergia.jpg} }
  \centerline{\includegraphics[height=7.5cm]{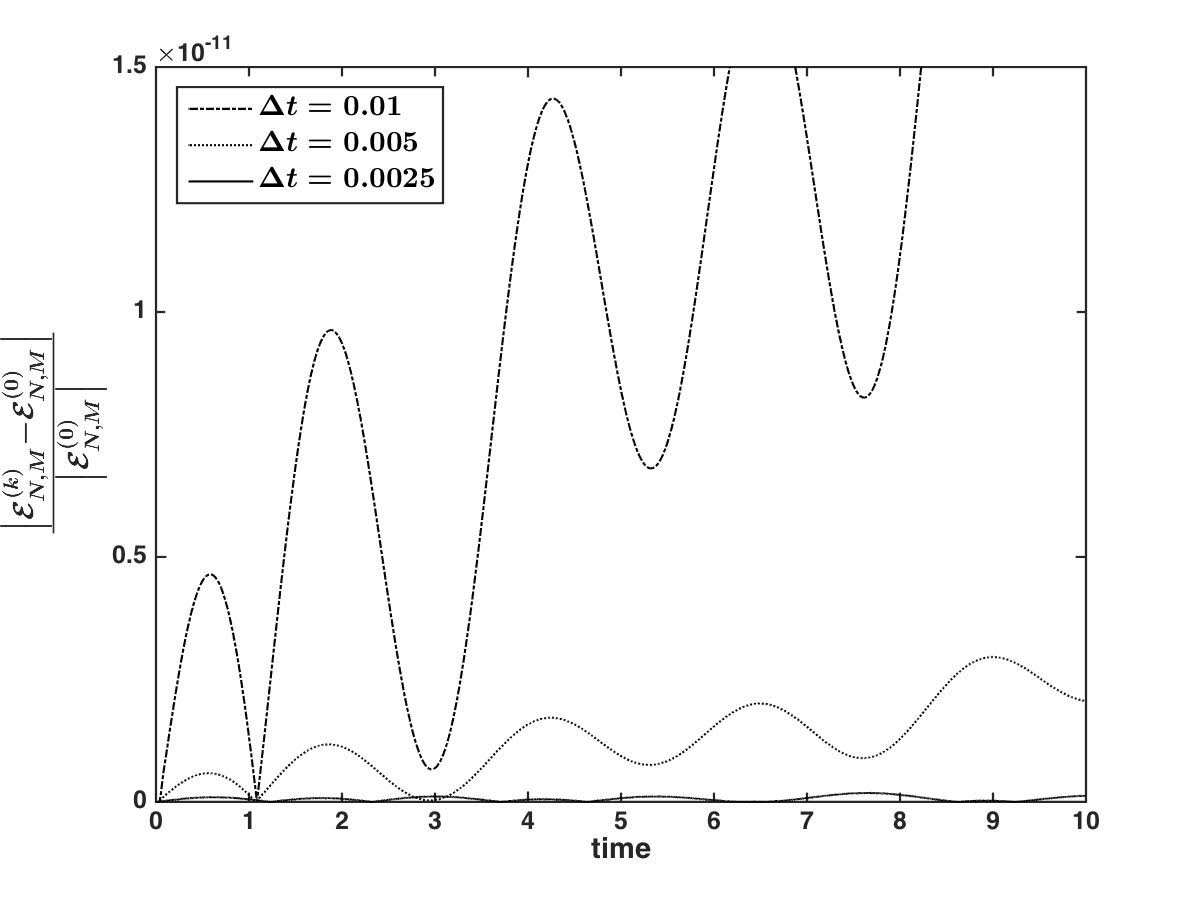}} %{./BDF3confrontoEnergia.jpg} }
  %%\centerline{\includegraphics[height=8cm]{./figures/JPG/new-fig4.jpg} }
  \caption{ Two-stream instability test: violation of total energy
    conservation using the second-order BDF scheme (top) and the third-order BDF scheme (bottom), for $N=2^{5}$, $M=2^{7}$, $T = 10$. 
    From the plots, it is clearly evident that such violation decays when the time step 
    diminishes.
    As predicted by the theoretical considerations of
    Section~\ref{sec:conservation:properties}, this decay is quadratic
    in the first case and cubic in the second case.  }
  \label{fig4two-streamConfrontoEnergia}
\end{figure}

%%%%%%%%%%%%%%%%%%%%%%%%%%%%%%%%%%%%%%%%%%%%%%%%%%%%%%%%%%
%%
To further validate our new schemes, we tested them on two standard
test cases of plasma physics: the two-stream instability and the
Landau damping (see next section).
To this end, in the two-stream instability problem, we set
$\Omega_{x}=[0,4\pi]$, $\Omega_{v}=[-5,5]$ in \eqref{eq:1D1V:V},
\eqref{eq:1D1V:Vci}, \eqref{eq:1D1V:EMaxwell},
\eqref{eq:1D1V:charge-density}.
The initial guess is given by:
\begin{align} \label{eq:1D1V:twostreamf}
  \bar{\fs}(\xs,\vs)=\frac{1}{2\alpha\sqrt{2\pi}}  
  \left[ 
    \exp { \left(-\frac{\vs-\beta}{\alpha\sqrt{2}}\right)^{2}}+
    \exp{ \left(-\frac{\vs+\beta}{\alpha\sqrt{2}}\right)^{2}}
  \right]
  \left[ 
    1 + \epsilon \cos\left( \kappa \xs \right)
  \right],
  %%\qquad\xs\in\Omega_{x},\,\,\vs\in\Omega_{v},
\end{align}
with $\alpha=1\slash{\sqrt 8}$, $\beta=1$, $\epsilon=10^{-3}$,
$\kappa=0.5$.
The exact solution is approximated by periodic functions in the
variables $\xs$ and $\vs$.
We integrate in time up to time $T=30$ using the second-order one-step
scheme \eqref{eq:1D1V:eqcloumod} with timestep $\Delta t=10^{-2}$.
This value is sufficiently small to guarantee stability, since the CFL
condition \eqref{eq:1D1V:CFL} requires $\Delta t$ to be proportional to
$1/\max\{N,M\}$.
The results of our simulations are presented in Figures
\ref{fig1two-streamApprox}, \ref{fig2two-streamImodoE},
\ref{fig3two-streamMassaMomento}  and
\ref{fig4two-streamConfrontoEnergia}.
In particular, in Figure~\ref{fig1two-streamApprox}, calculations are
carried out for different values of the discretization parameters $N$
and $M$.
The plots on the left show the interpolations of the initial solution
\eqref{eq:1D1V:twostreamf} with respect to the variable $\vs$ at
$\xs=0$.
Only in the top one there is a little disagreement, since the degrees
of freedom look not sufficient, which has, of course, a negative
reflection on the final solution.
The plots on the right show the corresponding numerical distribution
at the final time $T=30$.
The choice $N=M=2^{5}$ already gives reliable approximation results
but to completely eliminate the wiggles it is recommendable to
increase $M$ up to $2^{7}$.
Note, however, that the global number of degrees of freedom
$2^5\times 2^7=32\times128$ is rather low.

In Figure~\ref{fig2two-streamImodoE}, we plot the time evolution of
the ($\log$ of the) first Fourier mode of the electric field
$\Es^{(k)}_N$ in \eqref{eq:1D1V:FourierSeriesEnodi}, for different
values of the discretization parameters.
According to \eqref{eq:1D1V:coeffFourierSeries}, this is given by
$|\hat{\as}_{1}^{(k)}|$.
In particular, the plots show $|\hat{\as}_{1}^{(k)}|$ versus time,
when $N=2^{5}$, $M=2^{7}$, $\Delta t=10^{-2}$ and $T=100$.
These results are in agreement with the behavior expected from the
theory.
In particular, the slope of the numerical curves in the initial part
of the dynamics, where the two-streams instability starts developing,
matches well the slope predicted by the linear theory.
The stability of the numerical method is shown by the ``plateau'' up
to the final time $T=100$, which implies that the method is also
suitable for long-time integration.

To study the capability of the proposed schemes to preserve physical
invariants, we compute the variation with respect to the initial value
of the following quantities:
\begin{align} 
  \left| 
    \Qs^{(k)}_{N,M} - \Qs^{(0)}_{N,M}
  \right|, %% \qquad\quad k=0,1,\ldots, K,
  \label{eq:1D1V:varQ}
\end{align}
and
\begin{align}
  \left| 
    \Ps^{(k)}_{N,M} - \Ps^{(0)}_{N,M} 
  \right|, %% \qquad\quad k=0,1,\ldots, K,
  \label{eq:1D1V:varQ}
\end{align}
where the formulas for the  discrete number of particles
$\Qs^{(k)}_{N,M}$ and the discrete momentum $\Ps^{(k)}_{N,M}$ are
defined in \eqref{eq:1D1V:massadiscreta}
and~\eqref{eq:1D1V:momentodiscreto}, respectively.
The results of this study are given in
Figure~\ref{fig3two-streamMassaMomento}, for different time-marching
schemes.
The plots show (in a semi-$\log$ diagram) the variation versus time of
the number of particles and the momentum, with respect to their
initial value, when $N=2^{5}$, $M=2^{7}$, $T = 10$ and $\Delta t=5\cdot10^{-3}$.
In the first case (top), the results are excellent (i.e., within the
machine precision).
In the other cases, a weak growth in time is observed, probably due to
the accumulation of rounding errors.

To study the conservation of the total energy, we computed the
relative variation of the discrete energy with respect to the initial
value:
\begin{equation} \label{eq:1D1V:varEnergy}
  \frac{\left| \mathcal{\Es}_{N,M}^{(k)} - \mathcal{\Es}_{N,M}^{(0)}  \right|}{\left|\mathcal{\Es}_{N,M}^{(0)}\right|},  
  %% \qquad  k=0,1,\ldots, K,
\end{equation}
where $\mathcal{\Es}_{N,M}^{(k)}$ is defined in
\eqref{eq:1D1V:energiadiscreta}.
The results of Figure \ref{fig4two-streamConfrontoEnergia} show (in a
semi-$\log$ diagram) the behavior of the above quantities for
different values of the timestep $\Delta t$, for $N=2^5$, $M=2^7$ and
$T=10$.
Here, we implemented the second-order BDF scheme and the third-order
BDF scheme.
The energy is not perfectly preserved, but the discrepancy decays fast
by diminishing $\Delta t$, according to the accuracy of the method.
Indeed, these plots show that the decay rate for the first scheme is
quadratic, while that of the second scheme is cubic.
It has to be observed that this last method requires a more
restrictive condition on the timestep.
First of all, this is true because of the smaller domain of stability
of BDF high-order methods.
Secondly, because in the build-up of the method we trace back the
characteristic curves of several multiples of $\Delta t$ (see, for
instance, the second relation in~\eqref{eq:1D1V:appf}).
  
%%%%%%%%%%%%%%%%%%%%%%%%%%%%%%%%%%%%%%%%%%%%%%%%%%%%%%%%%%
\subsection{Landau damping}
%%%%%%%%%%%%%%%%%%%%%%%%%%%%%%%%%%%%%%%%%%%%%%%%%%%%%%%%%%
%%
In the following numerical tests, the proposed numerical schemes are
applied in order to capture the Landau damping phenomenon.
Landau damping is a classical kinetic effect in warm plasmas due to
the resonance of the particles with an initial wave perturbation.
In this classical and well-studied example, the continuous
filamentation process in velocity space occurs.

We initialize the electron Maxwellian distribution with a suitable
perturbation as follows:
\begin{align} \label{eq:1D1V:landau}
  &
  \fs(0,\xs,\vs)=\frac{1}{\sqrt{2\pi}}  
  \left[
    1 + \gamma \cos\left(\kappa   \xs   \right)
  \right] \exp (- \vs^{2}/2), 
  %% \qquad\,\,\xs\in\Omega_{x},\,\,\vs\in\Omega_{v},
\end{align}
where $\gamma$ is the size of the perturbation and $\kappa$ is the
wave-number. 
For this test, we set $\Omega_{x}=[0,4\pi]$ and $\Omega_{v}=[-10,10]$.
The size of $\Omega_v$ ensures that the values attained by $f$ at
$v=\pm10$ are negligible.

%%%%damping%%%%%%%%%%%%%%%%%%%%%%%%%%%%%%%%%%%%%%%%%%%%%%%%%%%%%%
\subsubsection{Linear Landau damping}
\begin{figure}
  \centerline{
    \includegraphics[height=6.25cm]{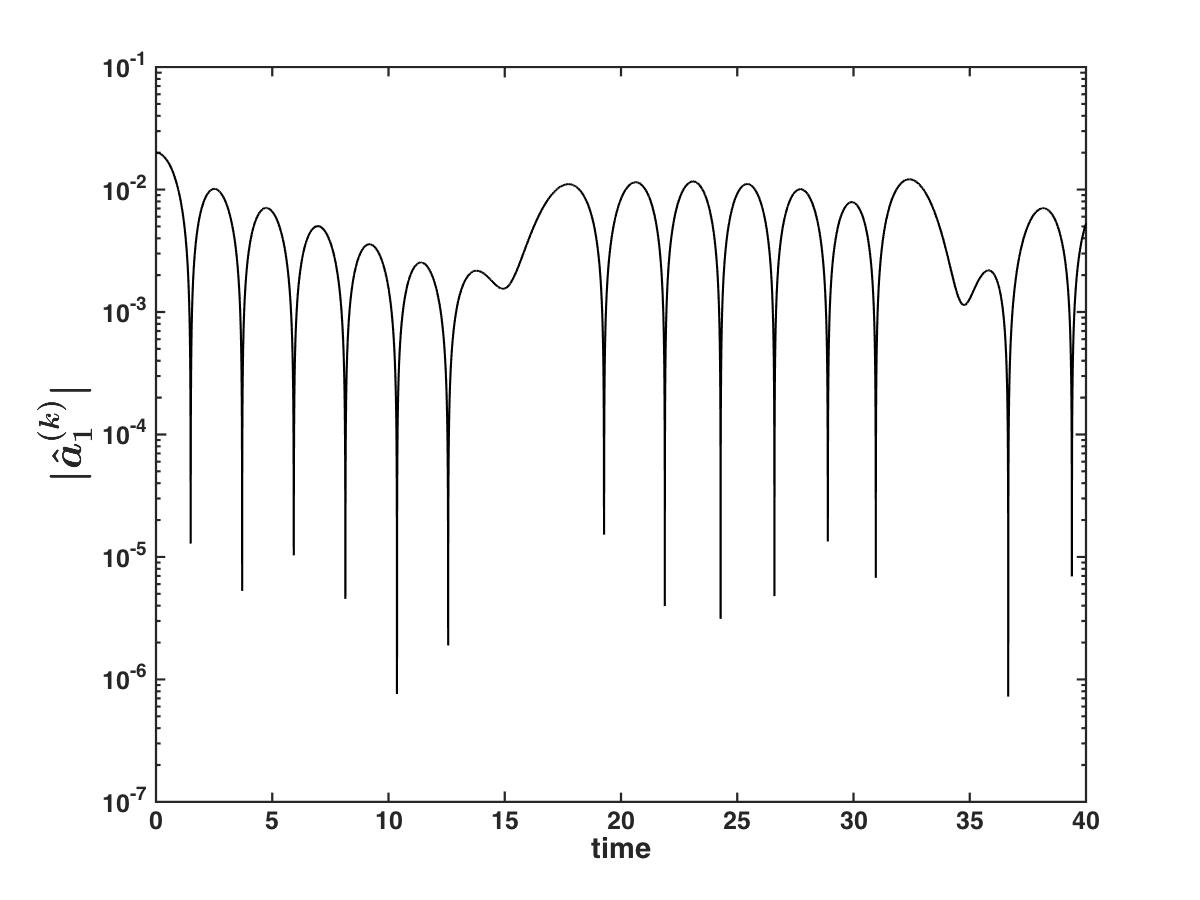}%{./LinearLDBDF2ImodoENM2^5.jpg} 
    \includegraphics[height=6.25cm]{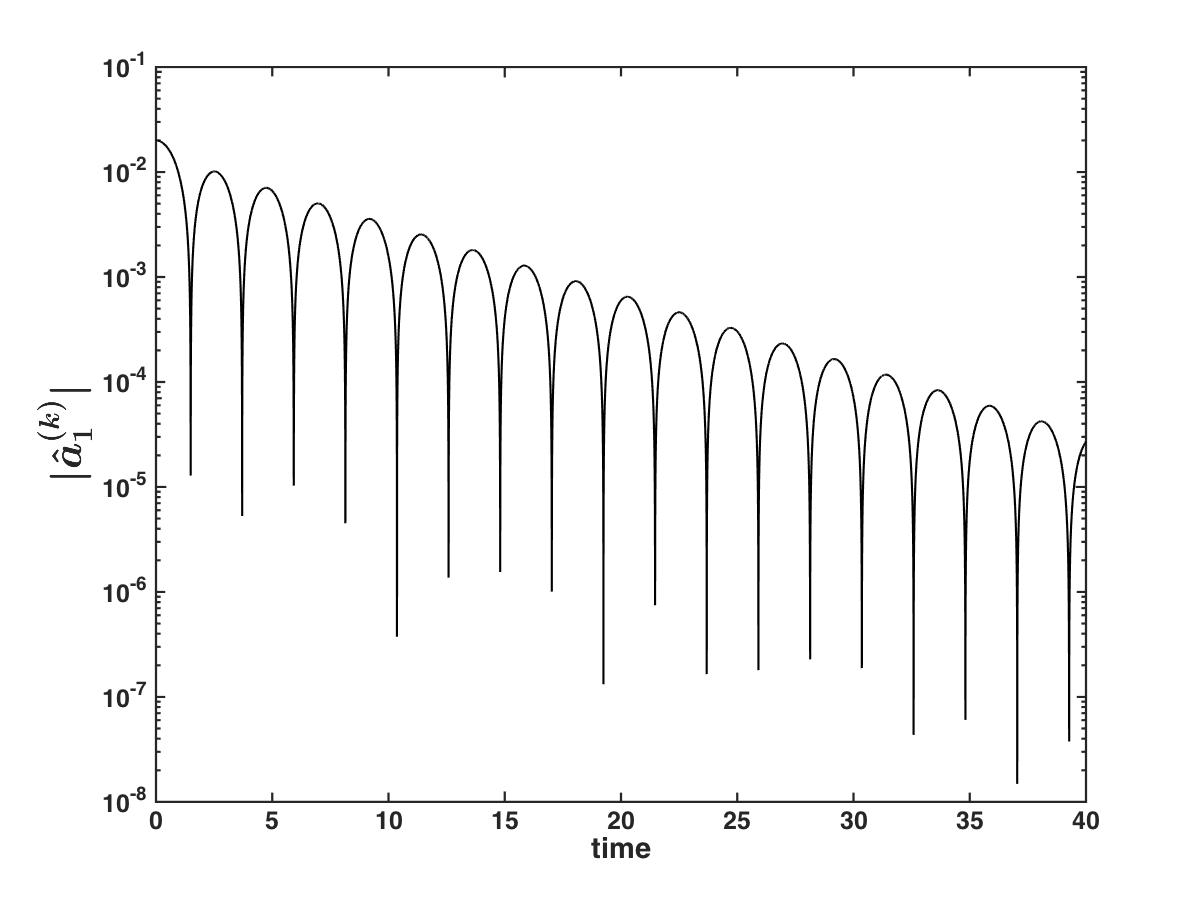}%{./LinearLDImodoEN2^5M2^7.jpg} 
  }
  \caption{ Linear Landau damping test: the first Fourier mode $|\hat
    a_{1}^{(k)}|$ of the electric field $|E^{(k)}_N|$ versus time, for
    the second-order BDF scheme with $T=40$, $\Delta t=2.5\cdot10^{-3}$
    and $N=M=2^{5}$ (left), $N=2^{5}$, $M=2^{7}$ (right). }
  \label{fig1LinearLD}
\end{figure}

%%%%%%%%%%%%%%%%%%%%%%%%%%%%%%%%%%%%%%%%%%%%%%%%%%%%%%%%%%
In this example, we set $\gamma=0.01$ and $\kappa=0.5$ in
\eqref{eq:1D1V:landau}.
Here, the perturbation is small and therefore the plasma behaves
according to the linear Landau theory.
The solution is computed up to time $T=40$ by using the second-order
BDF scheme in time with $\Delta t=2.5\cdot10^{-3}$ and $N=M=2^{5}$ (left),
$N=2^{5}$, $M=2^{7}$ (right).
Figure~\ref{fig1LinearLD} shows the behaviour in time of the first
Fourier mode of the electric field $\Es^{(k)}_N$ (see
$|\hat{\as}_{1}^{(k)}|$ in~\eqref{eq:1D1V:FourierSeriesEnodi}) in the
$\log$ scale.
The recurrence phenomenon starting at time $t\approx 12$ is clearly
visible on the left plot, which is due to an insufficient resolution
of the velocity domain.
This effect can be mitigated by increasing the accuracy of the
velocity approximation (we recall that we do not have any artificial
dissipation term in these schemes).
The plot on the right shows how the method performs when $M=2^7$
velocity degrees of freedom are used.
A similar behavior has been observed also for the other discretization
schemes proposed in this paper.

%%%%%%%%%%%%%%%%%%%%%%%%%%%%%%%%%%%%%%%%%%%%%%%%%%%%%%%%%%
\subsubsection{Nonlinear Landau damping}
%%%%%%%%%%%%%%%%%%%%%%%%%%%%%%%%%%%%%%%%%%%%%%%%%%%%%%%%%%
The initial distribution is again the function in
\eqref{eq:1D1V:landau}, but this time we set $\gamma=0.5$.
The other parameters are the same as in the linear Landau damping.
Therefore, a larger amplitude of the initial perturbation is used.
In this situation, the Landau linear theory does not hold, because the
nonlinear effects become relevant.
Nevertheless, several results obtained numerically are available in
the literature, since the nonlinear Landau damping is often used to
assess the performance of Vlasov-Poisson solvers.
  
\begin{figure}
  \centerline{
    \includegraphics[height=6.25cm]{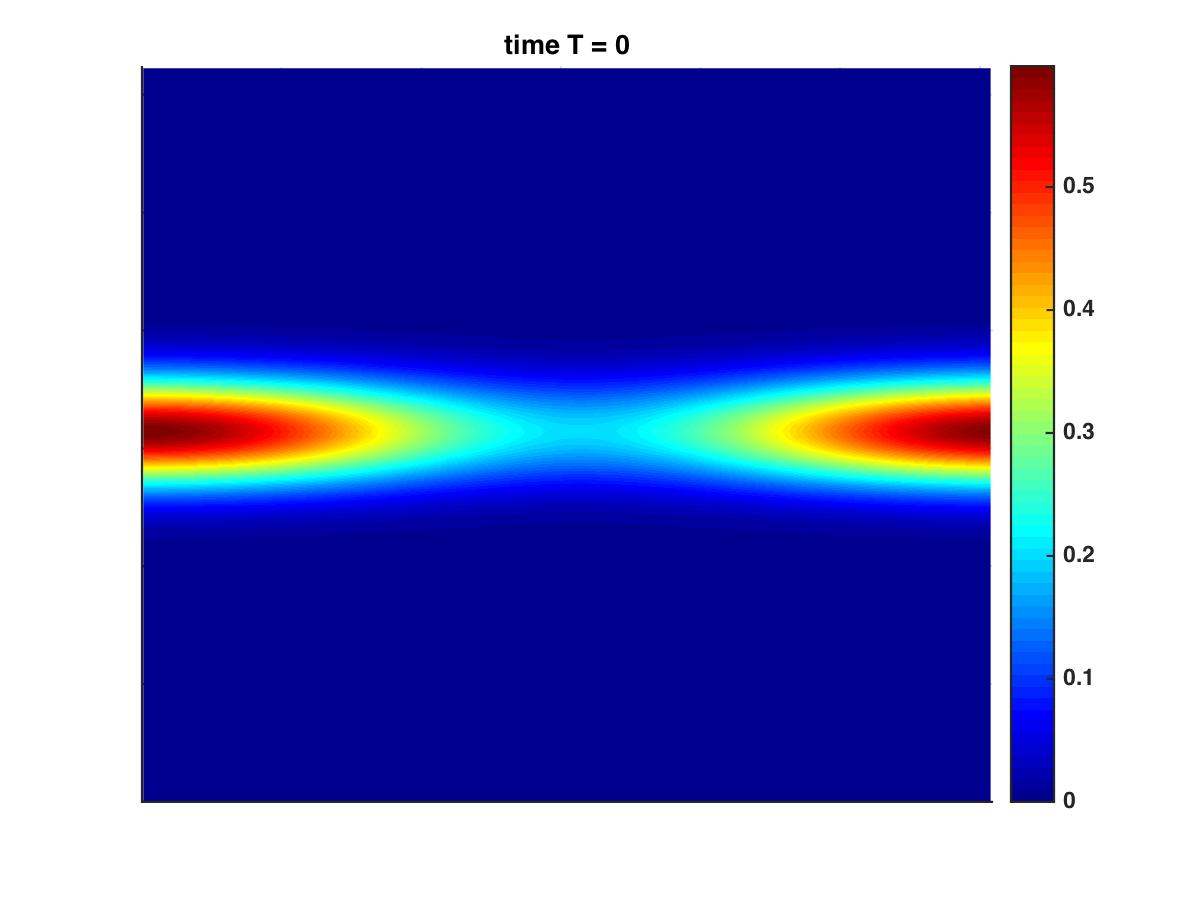}%{./NonlinearLDBDF2T0N2^5M2^7.jpg}      
    \includegraphics[height=6.25cm]{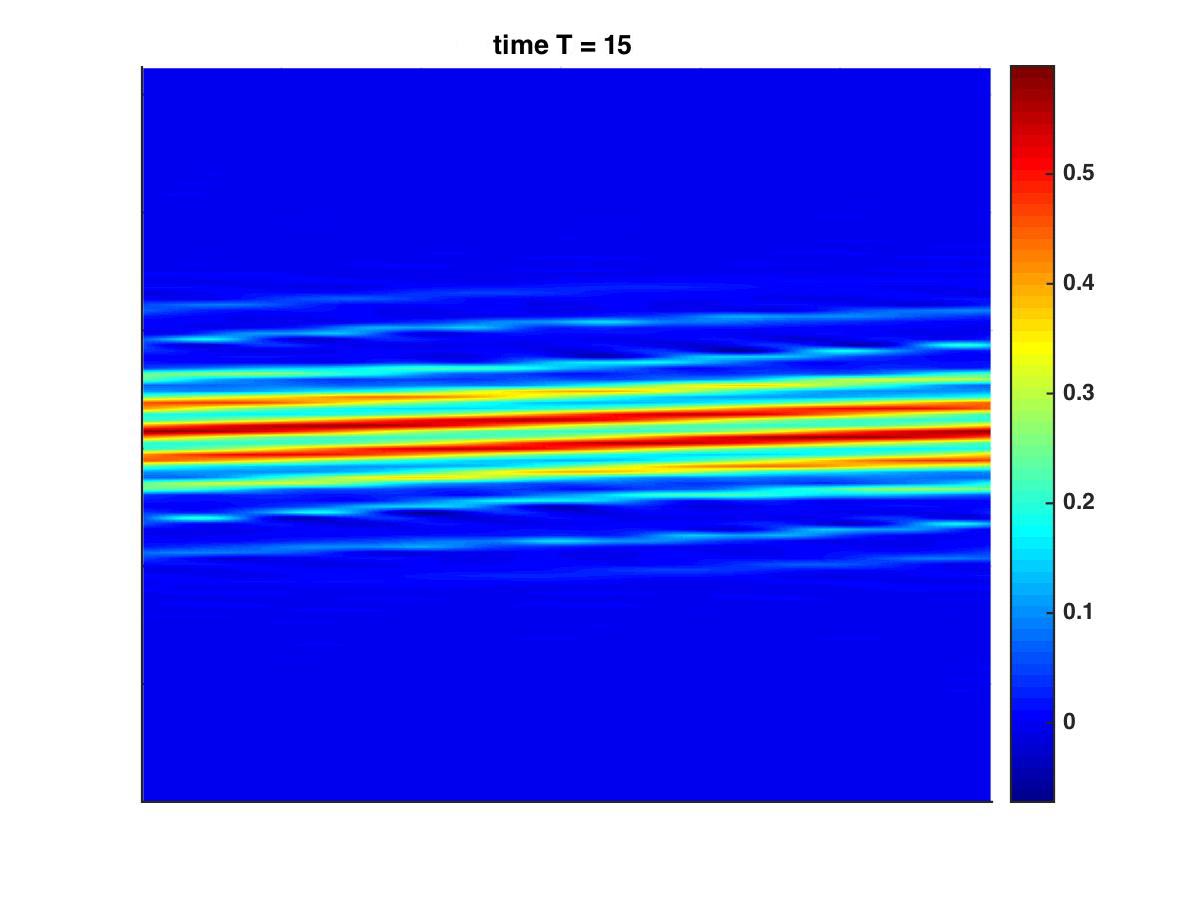}%{./NonlinearLDBDF2T15N2^5M2^7.jpg}      
  }
  \centerline{
      \includegraphics[height=6.25cm]{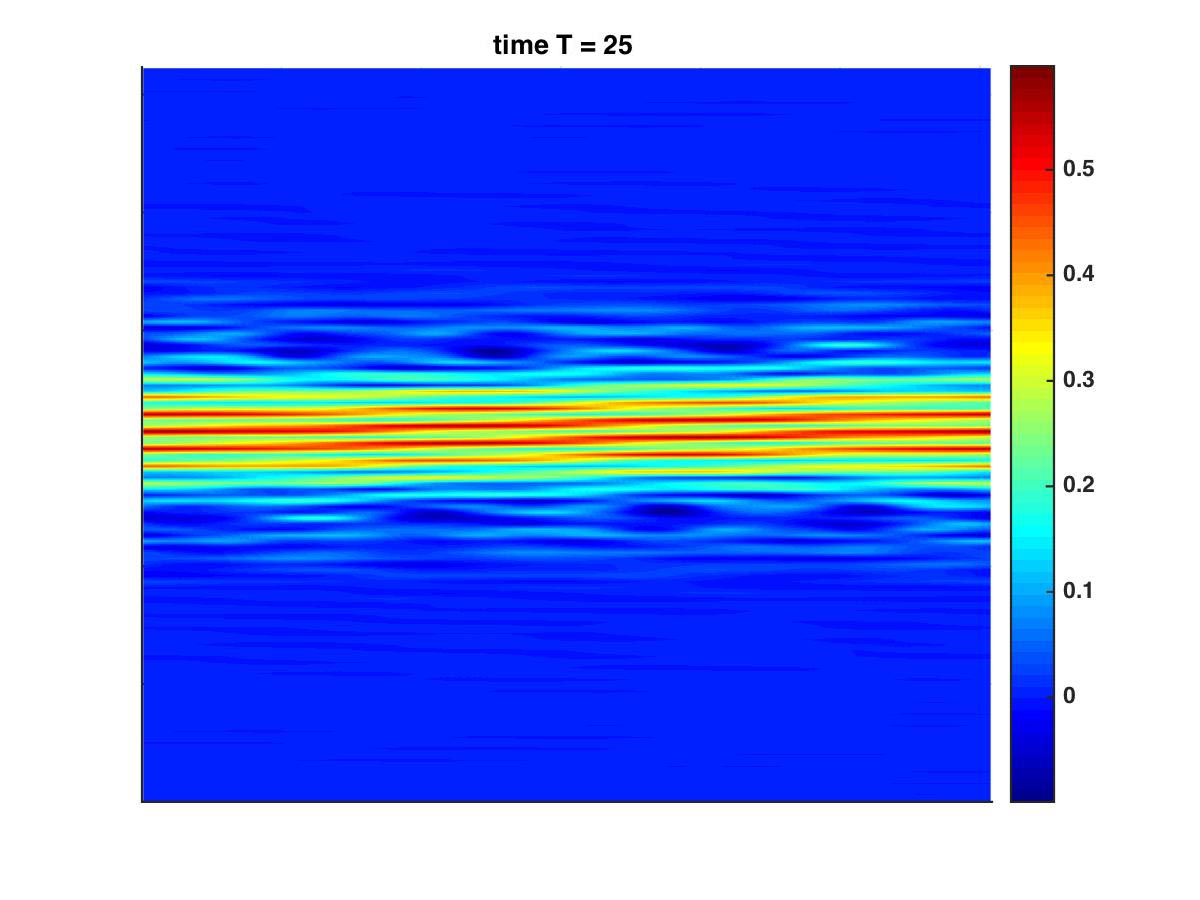}%{./NonlinearLDBDF2T25N2^5M2^7.jpg}      
      \includegraphics[height=6.25cm]{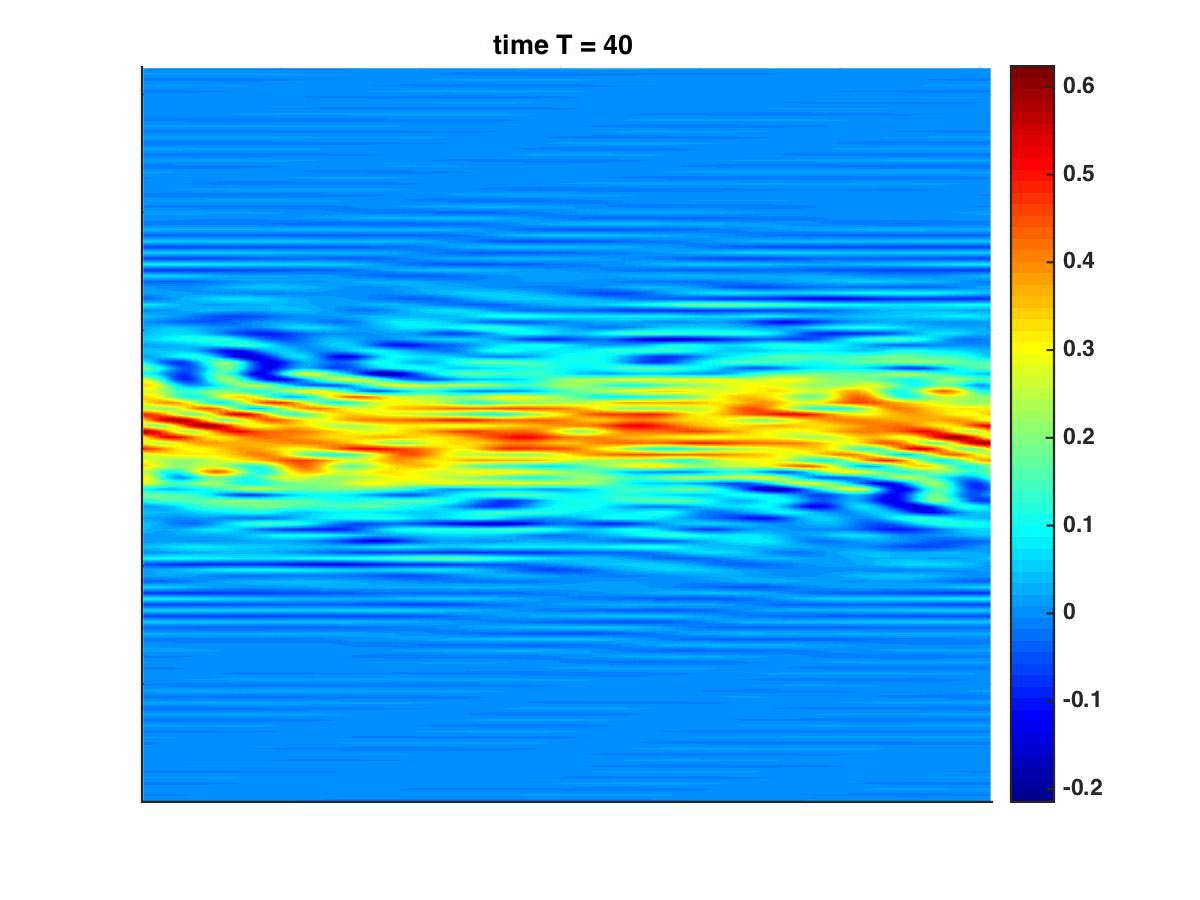}%{./NonlinearLDBDF2T40N2^5M2^7.jpg}     
    }
    \caption{ Nonlinear Landau damping test: approximated distribution
      functions obtained by using the second-order BDF scheme, with
      $N=2^{5}$, $M=2^{7}$ and $\Delta t=2.5\cdot10^{-3}$.
      Using the one-step second-order scheme with the same parameters gives
      exactly the same results.
  }
  \label{fig1NonlinearLD}
\end{figure}

\begin{figure}
  \centerline{
    \includegraphics[height=6.25cm]{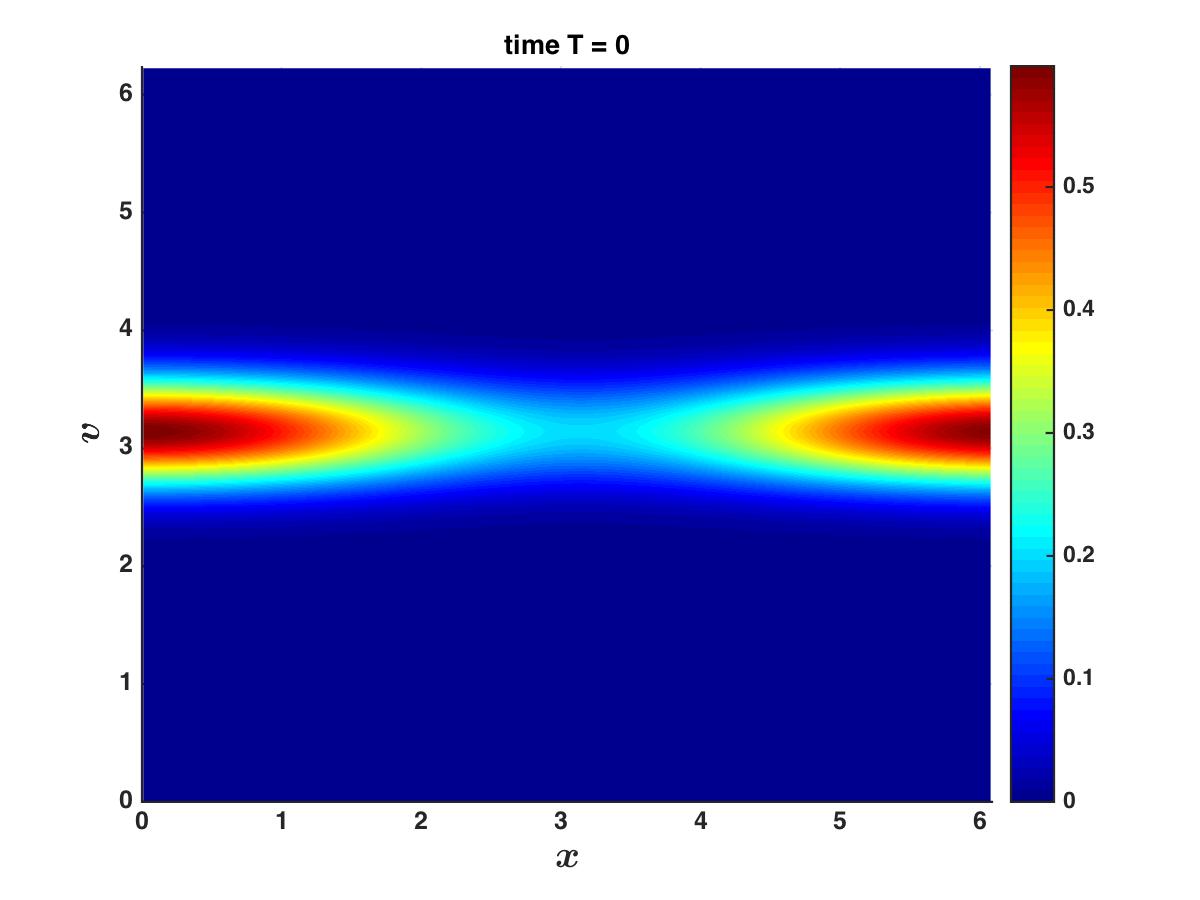}%{./NonlinearLDT0N2^5M2^7_bis.jpg}      
    \includegraphics[height=6.25cm]{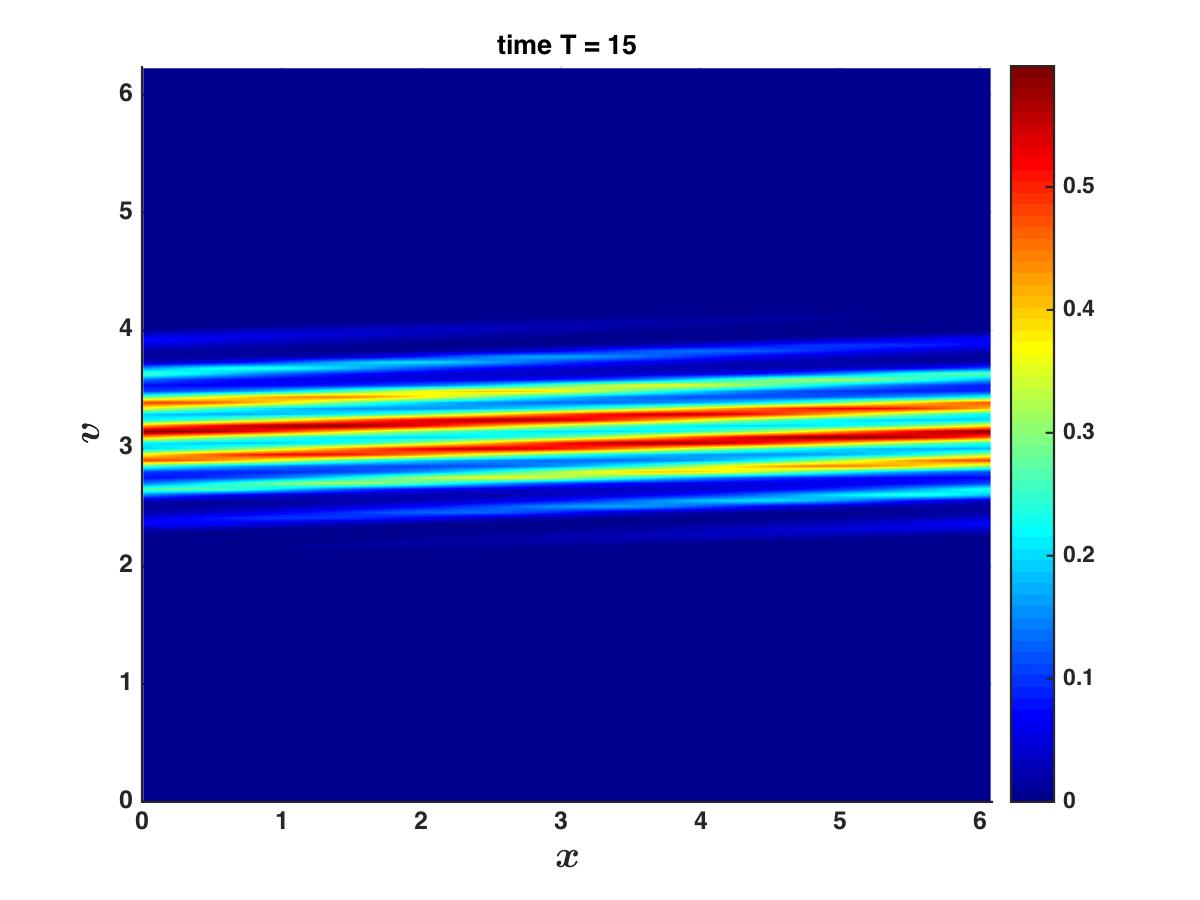}%{./NonlinearLDT15N2^5M2^7_bis.jpg}      
  }
  \centerline{
      \includegraphics[height=6.25cm]{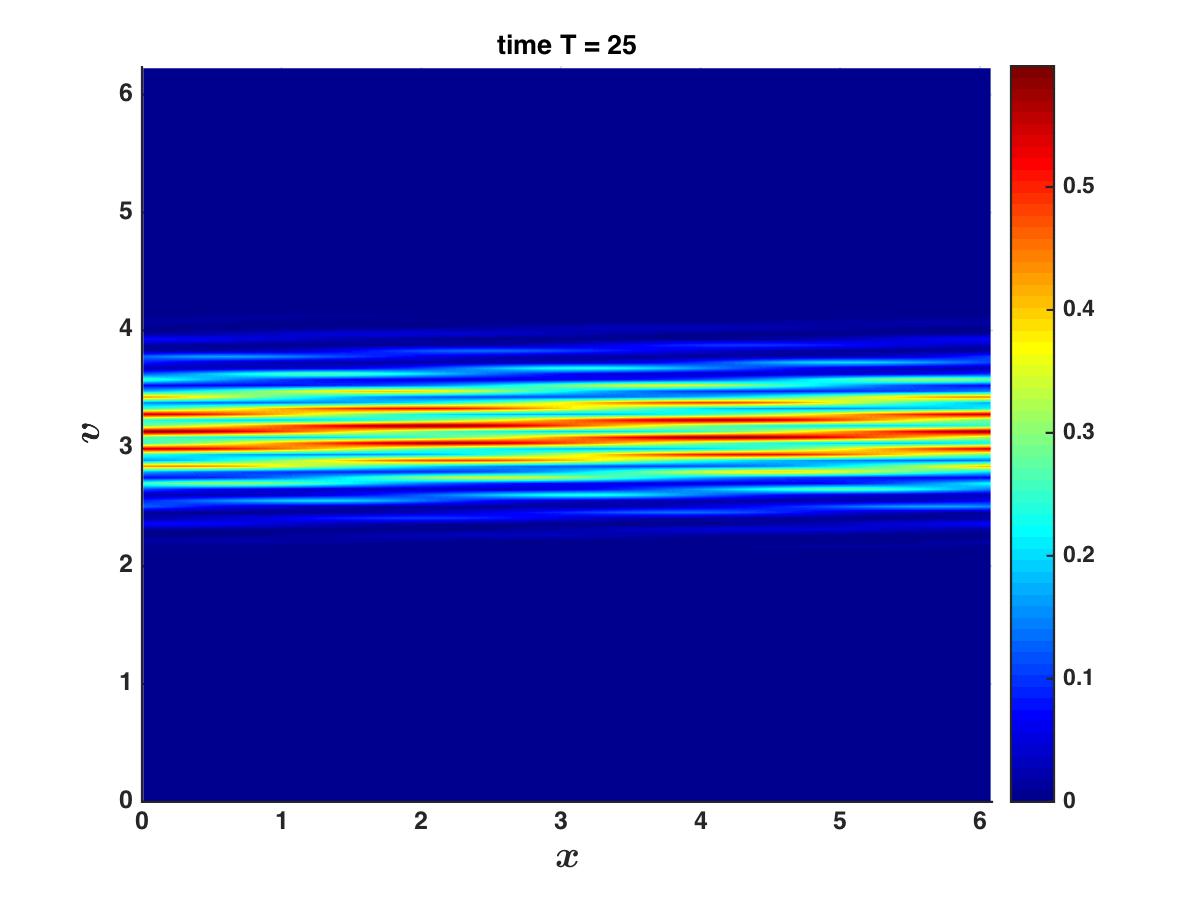}%{./NonlinearLDT25N2^5M2^7_bis.jpg}      
      \includegraphics[height=6.25cm]{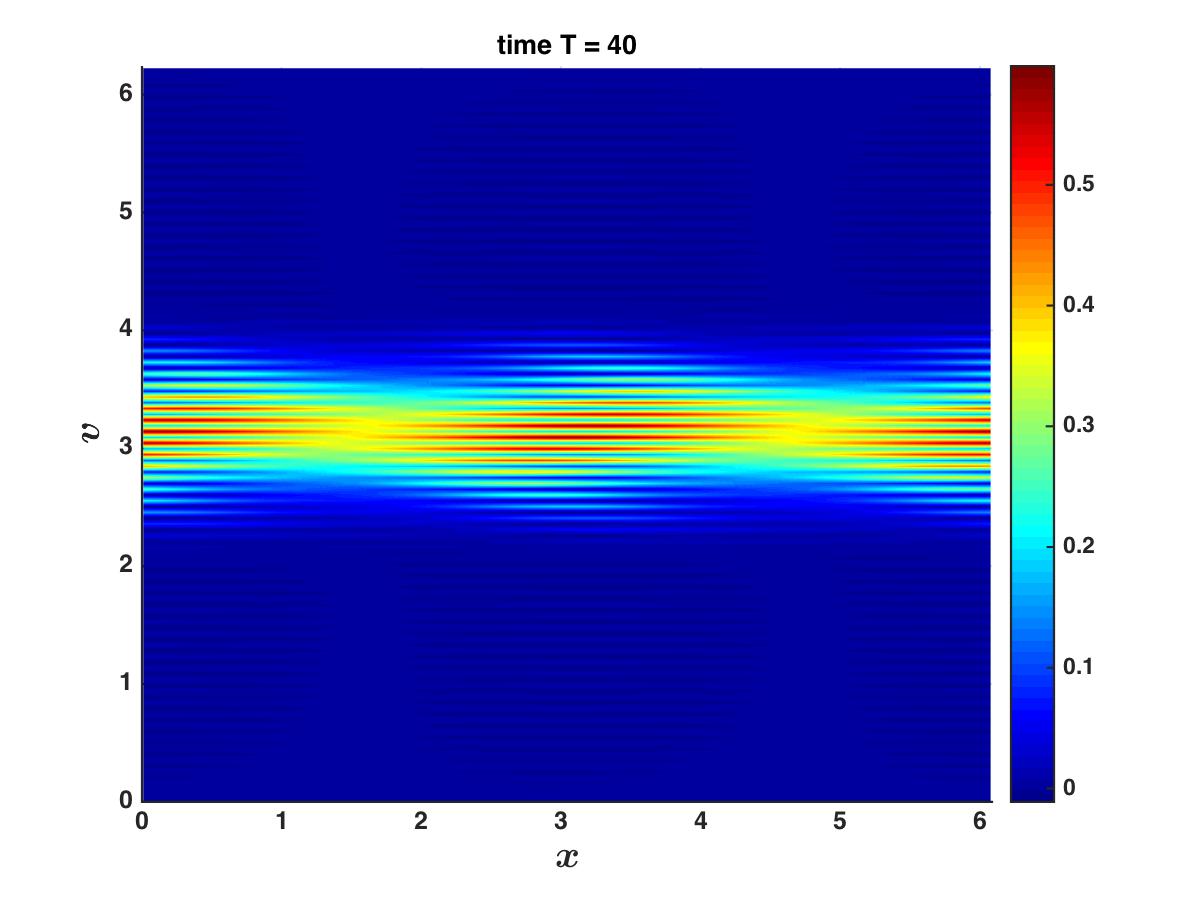}%{./NonlinearLDT40N2^5M2^7_bis.jpg}     
    }
    \caption{ Nonlinear Landau damping test: approximated distribution
      functions obtained by using the one-step second-order
      time-marching scheme, with $N=2^{5}$, $M=2^{7}$ and $\Delta
      t=5\cdot10^{-3}$. Note that the timestep is twice that of the
      calculation shown in Fig.~\ref{fig1NonlinearLD}.  }
  \label{fig1NonlinearLDMio}
\end{figure}

Figure \ref{fig1NonlinearLD} shows the plots at different times for
the computation relative to the second-order BDF scheme.
In this example, we work with $N=2^{5}$, $M=2^{7}$, $\Delta t=2.5\cdot10^{-3}$,
and $T=40$.
In these plots, the filamentation effect is clearly evident and it is
due to the fact that we do not have any explicit artificial
dissipation term in the method.
%5
The one-step second-order scheme provides identical results when is
applied with the same parameters.
However, the latter method has less restriction on the timestep than
the former one (see also the comments at the end of
Section~\ref{sec32}).
Therefore, we can run the same simulation with $\Delta t=5\cdot10^{-3}$.
The results are shown in Figure~\ref{fig1NonlinearLDMio}.
Filamentation is still visible, but less evident probably because of
some numerical diffusion due to the choice of a larger timestep.

\begin{figure}
  \centerline{
      \includegraphics[height=6.25cm]{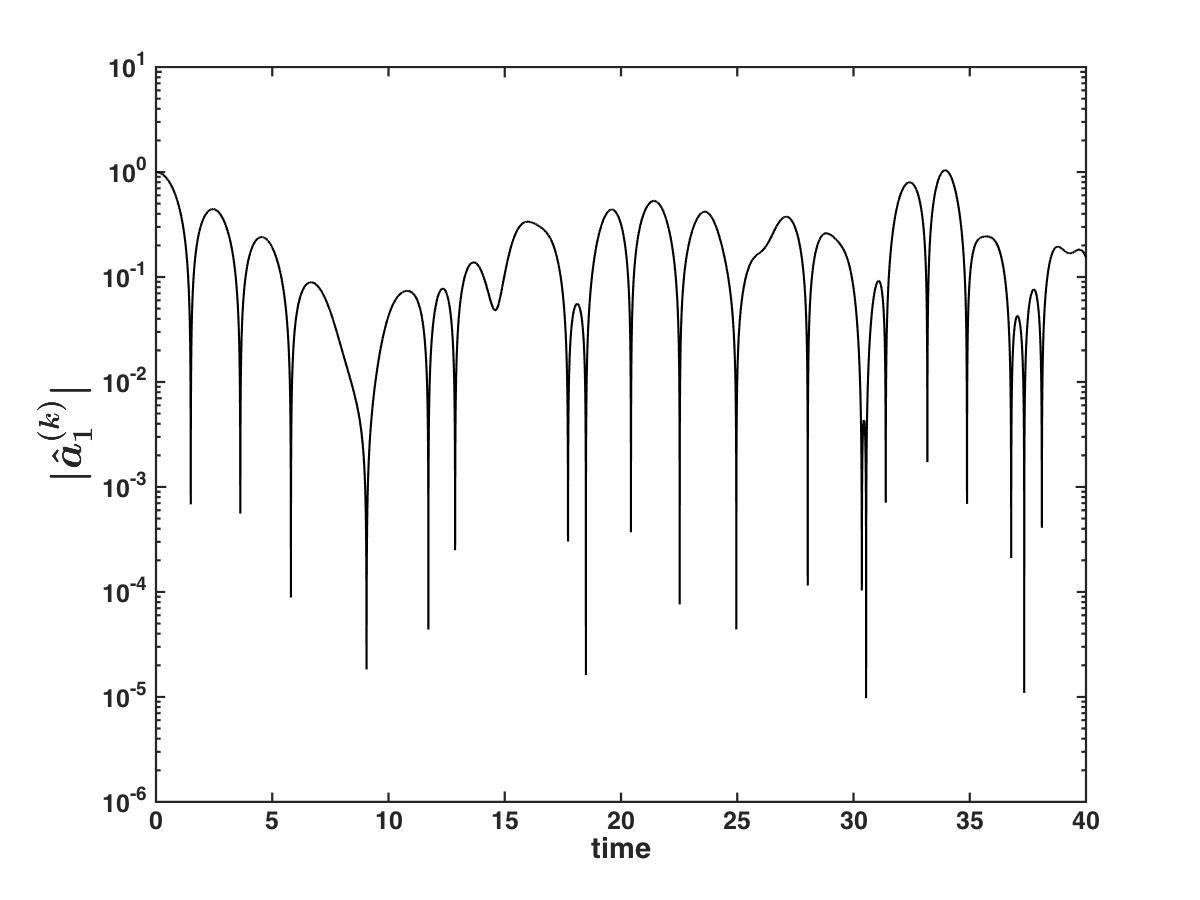}%{./NonLinearLDBDF2ImodoENM2^5.jpg} 
      \includegraphics[height=6.25cm]{plot27.jpg}%{./NonlinearLDDBDF2ImodoEN2^5M2^7.jpg} 
    }
    \caption{ Nonlinear Landau damping test: the first Fourier mode
      $|\hat a_{1}^{(k)}|$ of the electric field $|E^{(k)}_N|$ versus
      time, obtained by using the second-order BDF scheme, with $T=40$, $\Delta
      t=2.5\cdot10^{-3}$ and $N=M=2^{5}$ (left) and $N=2^{5}$, $M=2^{7}$
      (right).
      Using the one-step second-order scheme with the same parameters gives
      exactly the same results. 
    }
   \label{fig2NonlinearLD}
\end{figure}

Finally, in Figure~\ref{fig2NonlinearLD} we show the first Fourier
mode of the electric field $\Es^{(k)}_N$ in the log-scale computed
with the second-order BDF scheme for $\Delta t=2.5\cdot10^{-3}$, and using
$N=2^{5}$, $M=2^{7}$ on the right and $N=M=2^{5}$ on the left.
Again, the different behavior when more degrees of freedom are used
for the velocity representation is reflected by the comparison of the
corresponding curves.

\section{Conclusions}
\label{sec:conclusions}

In this work, a class of novel numerical methods for the system of equations of
Vlasov-Poisson has been designed, developed, and investigated
esperimentally.
These methods are based on a spectral approximation in the phase space in
a Semi-Lagrangian framework using a first- and a second-order accurate
approximation of the characteristics curves.
A single-step second-order method is thus obtained without resorting
to any splitting of the equations.
High-order time discretizations based on the  method-of-lines
approach are also proposed and studied, which are obtained by adopting
second-order and third-order multi-step Backward Differentiation  Formulas
(BDF).
Furthermore, conservation properties have been also investigated.
The performance of these methods  has been assessed by
thorugh a manufactured solution and standard benchmark problems as the
two stream instability and the Landau damping.

%\input figs.tex

%%%%%%%%%%%%%%%%%%%%%%%%%%%%%%%%%%%%%%%%%%%%%%%%%%%%%%%%%%
% Figure Two-stream instability
%%%%%%%%%%%%%%%%%%%%%%%%%%%%%%%%%%%%%%%%%%%%%%%%%%%%%%%%%%
%%
%% Fig. 1
%% 
%% 
%% Fig. 2
%% 
%%
%% Fig. 3
%%
%%
%% Fig. 4
%%
%%  
%%%%%%%%%%%%%%%%%%%%%%%%%%%%%%%%%%%%%%%%%%%%%%%%%%%%%%%%%% 
% Figure Linear Landau damping}
%%%%%%%%%%%%%%%%%%%%%%%%%%%%%%%%%%%%%%%%%%%%%%%%%%%%%%%%%% 
%% 
%% Fig. 5
%%

% LinearLDBDF2ImodoENM2^5.jpg

%%%%%%%%%%%%%%%%%%%%%%%%%%%%%%%%%%%%%%%%%%%%%%%%%%%%%%%%%% 
% Figure Non-Linear Landau damping N=2^5,M=2^7 BDF2
%%%%%%%%%%%%%%%%%%%%%%%%%%%%%%%%%%%%%%%%%%%%%%%%%%%%%%%%%%
%% 
%% Fig. 6
%% 

%%%%%%%%%%%%%%%%%%%%%%%%%%%%%%%%%%%%%%%%%%%%%%%%%%%%%%%%%% 
% Figure Non-Linear Landau damping N=2^5,M=2^7 I step II order
%%%%%%%%%%%%%%%%%%%%%%%%%%%%%%%%%%%%%%%%%%%%%%%%%%%%%%%%%%
%% 
%% Fig. 7
%% 

%%
%% Fig. 8
%%

\section*{Acknowledgements}
The second author was partially supported by the \emph{Short Term
  Mobility Program} of the Consiglio Nazionale delle Ricerche
(CNR-Italy), which partially funded a scientific visit to the Los
Alamos National Laboratory.
The third author was supported by the Laboratory Directed Research and
Development Program (LDRD), U.S. Department of Energy Office of
Science, Office of Fusion Energy Sciences, and the DOE Office of
Science Advanced Scientific Computing Research (ASCR) Program in
Applied Mathematics Research, under the auspices of the National
Nuclear Security Administration of the U.S. Department of Energy by
Los Alamos National Laboratory, operated by Los Alamos National
Security LLC under contract DE-AC52-06NA25396.

%\bibliographystyle{abbrv}
%% \bibliography{./TMP/vlasov.bib}

\end{document}